%% file: manuscript.tex
\documentclass[final,11pt,3p,doublespacing]{elsarticle}
\journal{journal for consideration.}
\usepackage{XCharter}
\usepackage{booktabs}
\emergencystretch=\maxdimen
\usepackage{titlesec}  % For more title formatting control
\usepackage{sectsty}
\sectionfont{\color{black}\large}
\subsectionfont{\color{black}}
\subsubsectionfont{\color{black}\itshape}
\usepackage{algorithm}
\usepackage{algpseudocode}
\usepackage{amsmath, amssymb}
\biboptions{sort&compress}
\usepackage[colorlinks]{hyperref}
\AtBeginDocument{%
	\hypersetup{
		plainpages = false, 
		bookmarksopen = true,
		bookmarksnumbered = true,
		breaklinks = true,
		linktocpage,
		colorlinks = true,
		linkcolor = purple,
		urlcolor  = black,
		citecolor = Magenta,
}}
\usepackage{multicol}
\usepackage{multirow}
\usepackage{imakeidx}
\makeindex
\usepackage{nomencl}
\makenomenclature
\usepackage[xindy]{glossaries}
\makeglossaries
\usepackage{graphics}
\usepackage{graphicx}
\usepackage{epsf,psfrag}
\usepackage{epstopdf}
\usepackage[para]{threeparttable}
\usepackage{subfigure}
\usepackage{epsfig}
\usepackage{pdflscape}
\usepackage{rotating}
\usepackage{lscape}
\usepackage{float}
\usepackage{stfloats}
\usepackage{textgreek}
\usepackage{longtable}
\usepackage{lscape}
\usepackage{capt-of}
\usepackage{tablefootnote}
\usepackage{footnote}
\usepackage{scrextend}
\usepackage{caption}
\usepackage[autopunct=true]{csquotes}
\usepackage{rotating}
\usepackage[normalem]{ulem}
\usepackage{resizegather}
\usepackage{siunitx}
\usepackage{cases}
\usepackage[nodisplayskipstretch]{setspace}
\usepackage{color,soul}
\usepackage[usenames,dvipsnames]{xcolor}
\usepackage{amsfonts,amsmath,amssymb,gensymb,stmaryrd}
\usepackage{latexsym,array}%,subeqn}
\usepackage{textcomp}

\newcommand{\RomanNumeralCaps}[1]
\linenumbers
\makesavenoteenv{tabular}
\makesavenoteenv{table}
\usepackage[left,displaymath, mathlines]{lineno}
\DeclareMathAlphabet{\mathpzc}{OT1}{pzc}{m}{it}
\def\fig{Fig.~}
\def\figs{Figs.~}
\def\eqn{Eq.~}
\def\eqns{Eqs.~}
\def\tab{Table~}
\def\tabs{Tables~}
\def\sect{Section~}

\def\tsc#1{\csdef{#1}{\textsc{\lowercase{#1}}\xspace}}
\tsc{WGM}
\tsc{QE}
\tsc{EP}
\tsc{PMS}
\tsc{BEC}
\tsc{DE}
\usepackage{easyReview} \setreviewsoff
\newcommand{\removeEq}[1]{\ifistoreview{\@\expandafter\removeColor{#1}\hspace{-0.6em}} \else {}\fi}
\graphicspath{{../}{Figures/}{Figures/Validation/}}
\DeclareGraphicsExtensions{.png, .PNG,.pdf,.PDF,.jpg,.JPG,.eps,.EPS}
\begin{document}
%--------------------------------------------
%\input{reviewer-1.tex}\newpage
%\setcounter{page}{1}
%--------------------------------------------
%\input{highlights.tex}\newpage
%--------------------------------------------
\setcounter{page}{1}
\begin{frontmatter} % for elsevier jrnls
% 
%\openup -0.5em % 1em for double spacing
%\title{An Improved Formulation of the Cell-Centered Nodal Integral Method for Burgers’ Equation}
%\title{An Improved Cell-Centered Nodal Integral Method (RCCNIM) for Non-Linear Burgers' Equation: Accurate Formulation, Efficient Implementation, and Validation}
\title{A Robust Cell-Centered Nodal Integral Method (RCCNIM) for Nonlinear Burgers' Equation: Accurate Formulation, Efficient Implementation, and Validation}
\author[labela]{Nadeem Ahmed}%\ead{jdhakar@ch.iitr.ac.in}
\emailauthor{nadeem.pd@ch.iitr.ac.in}{N. Ahmed}
\author[labela]{Ram Prakash {Bharti}\corref{coradd}}%\ead{rpbharti@iitr.ac.in}
\emailauthor{rpbharti@iitr.ac.in}{R.P. Bharti}
\address[labela]{Complex Fluid Dynamics and Microfluidics (CFDM) Lab, Department of Chemical Engineering, Indian Institute of Technology Roorkee, Roorkee - 247667, Uttarakhand, India}

%
%\address[labela]{Department of Chemical Engineering, Indian Institute of Technology Roorkee, Roorkee - 247667, Uttarakhand, INDIA}
%
\cortext[coradd]{\textit{Corresponding author. }}
%
%%%%%%%%%%%%%%%%%%%%%%%%%%%%%%%%%%%%%%%%%%%%%%%%%%%%%%%%%%%%%%%%%%%%%%%%%%%%%%%%%%%%%%
\begin{abstract}
\fontsize{11}{16pt}\selectfont
%----Text of abstract
%\add{text} \remove{text} \highlight{text}
%
%An improved formulation of the recently developed cell-centered nodal integral method (CCNIM) is presented for the numerical solution of Burgers’ equation. In the proposed approach, referred to as the improved CCNIM (RCCNIM), the nonlinear convection term is reformulated prior to discretization by approximating the convective velocity with the cell-averaged value from the previous time level, rather than the present-time values employed in the original CCNIM formulation. As a result, the coefficients of the resulting system of algebraic equations depend only on known quantities and are evaluated once per time step. This strategy significantly reduces the computational cost while maintaining accuracy comparable to that of the original CCNIM formulation. The effectiveness of the proposed formulation is assessed by solving the nonlinear Burgers’ equation using both the original CCNIM and the RCCNIM, followed by a systematic comparison of their numerical results. The comparisons demonstrate that the RCCNIM consistently outperforms the original CCNIM in terms of computational efficiency while maintaining comparable accuracy, thereby providing a strong basis for its extension to more complex fluid flow problems. 
An improved formulation of the recently developed cell-centered nodal integral method (MCCNIM) is proposed for the numerical solution of the nonlinear Burgers’ equation. The improved scheme, referred to as RCCNIM, reformulates the nonlinear convection term prior to discretization by evaluating the convective velocity using cell-averaged values from the previous-time level, instead of the present-time approximation employed in the original MCCNIM formulation. This reformulation leads to a fully algebraic system whose coefficients depend only on known quantities and are therefore evaluated once per time step. As a result, the proposed method significantly reduces computational cost while retaining the accuracy of the original MCCNIM. The performance of RCCNIM is assessed through systematic numerical comparisons with MCCNIM for the nonlinear Burgers’ equation. The results demonstrate that RCCNIM achieves comparable accuracy with improved computational efficiency, indicating its potential for extension to more complex nonlinear fluid flow problems.

%%%%
\end{abstract}
%%%%%%%%%%%%%%%%%%%%%%%%%%%%%%%%%%%%%%%%%%%%%%%%%%%%%%%%%%%%%%%%%%%%%%%%%%%%%%%%%%%%%%
%%%% 
%%%%%%%%%%%%%%%%%%%%%%%%%%%%%%%%%%%%%%%%%%%%%%%%%%%%%%%%%%%%%%%%%%%%%%%%%%%%%%%%%%%%%%
\end{frontmatter}
%%%%%%%%%%%%%%%%%%%%%%%%%%%%%%%%%%%%%%%%%%%%%%%%%%%%%%%%%%%%%%%%%%%%%%%%%%%%%%%%%%%%%%
%\fontsize{12}{14pt}\selectfont
%\linespread{1.6}
%\doublespacing
%\openup 1em % 1em for double spacing
%\setstretch{2} 
%%
%%%%%%%%%%%%%%%%%%%%%%%%%%%%%%%%%%%%%%%%%%%%%%%%%%%%%%%%%%%%%%%%%%%%%
\section{Introduction}
\label{sec:1}
%%%%%%%%%%%%%%%%%%%%%%%%%%%%%%%%%%%%%%%%%%%%%%%%%%%%%%%%%%%%%%%%%%%%%
%%%
Nodal Integral Methods (NIM) constitute a class of coarse-mesh techniques that have proven effective for solving a wide range of partial differential equations (PDEs) encountered in scientific and engineering applications \citep{05Dorning_1979,11Azmy_1983,06Hennart_1986,07Lawrence_1986}. 
%\citep{01Delp_1964,02Steinke_1973,03Burns_1975,04Lawrence_1980,05Dorning_1979,06Hennart_1986,07Lawrence_1986}. 
These methods offer an attractive alternative for large-scale simulations by achieving accuracy comparable to that of conventional fine-mesh methods at significantly lower computational cost. Their efficiency stems from their primary step known as the transverse integration procedure (TIP), which reduces the governing PDEs to a set of averaged ordinary differential equations (ODEs). The availability of analytical or semi-analytical solutions to these ODEs enables accurate numerical approximations even on very coarse spatial meshes \citep{11Azmy_1983,esser1993upwind,22Rizwan_uddin_1997}. To extend these advantages to complex fluid flow problems and improve numerical robustness, several variants of NIM have been developed in the literature. Notable examples include the Upwind Nodal Method (UNM) \citep{esser1993upwind}, the Cell Analytical-Numerical Method (CAN) \citep{20Elnawawy_1990}, the Nodal Expansion Method (NEM) \citep{lee2011nodal,zhou2016general}, the Modified Nodal Integral Method (MNIM) \citep{22Rizwan_uddin_1997,23Rizwan_uddin_1997,25Wescott_2001,wang2005modified}, and the Cell-Centered Nodal Integral Method (CCNIM) \citep{19Ahmed_2021,18Ahmed_2023,17Ahmed_2024,Ahmed_Diffusion_2024,ahmed2025efficient}, among others \citep{kumar2012pressure,Niteen2022nodal,gander2022new,jarrah2024nodal,Neeraj2024coupled}. 

Most of the aforementioned NIM variants are already in a mature and well-established stage, with the exception of the CCNIM, which remains under active development. To date, CCNIM has been successfully applied to diffusion \citep{19Ahmed_2021,Ahmed_Diffusion_2024}, convection–diffusion \citep{18Ahmed_2023,17Ahmed_2024}, and Burgers’ equations \citep{ahmed2025efficient}. Owing to several fundamental advantages over traditional NIM formulations \citep{22Rizwan_uddin_1997,23Rizwan_uddin_1997}, CCNIM has the potential to outperform other NIM variants for fluid flow problems in both flexibility and applicability \citep{ahmed2025efficient}. A key advantage is that the final discrete system is expressed in terms of cell-centered variables, rather than surface-averaged variables as in conventional NIM schemes \citep{23Rizwan_uddin_1997,wang2005modified}. This formulation naturally avoids staggered-grid structures in traditional NIM and leads to a more compact and intuitive discretization  \citep{ahmed2025efficient}. Furthermore, CCNIM can be readily applied to problems involving second-order time derivatives, such as wave equations, for which traditional NIM approaches are often unsuitable due to their simultaneous discretization of space and time \citep{19Ahmed_2021,18Ahmed_2023}. The CCNIM formulation allows the governing equations to be cast as a system of differential–algebraic equations (DAEs), enabling the use of higher-order time-integration schemes and improved temporal accuracy \citep{18Ahmed_2023}. From an implementation standpoint, CCNIM is comparatively simpler than traditional NIM, particularly in higher dimensions, where surface-averaged variables in traditional NIM act as staggered grids which require complex local coordinate constructions \citep{19Ahmed_2021,18Ahmed_2023,17Ahmed_2024,Ahmed_Diffusion_2024,ahmed2025efficient}. An additional predictive advantage of CCNIM is its natural compatibility with immersed boundary techniques, such as cut-cell approaches for handling irregular geometries, owing to the cell-centered nature of the discrete variables defined at cell centroids.     

Initially, CCNIM was developed for the diffusion equation \citep{19Ahmed_2021} as a simplified reformulation of methods already established in the literature for the neutron diffusion equation \citep{09Shober_1977,13Shober1978}, before being extended to fluid flow problems. Following its successful application to diffusion problems, the method was extended to the linear two-dimensional convection–diffusion equation, where detailed error analyses were performed in comparison with traditional NIM schemes \citep{zhou2016general} to demonstrate the accuracy of CCNIM for convection–diffusion problems \citep{18Ahmed_2023}. Prior to extending the method to the nonlinear problems such as Burgers' equations, it was recognized that the CCNIM formulation developed for the linear convection–diffusion equation may not be directly applicable to nonlinear problems. This limitation arises because CCNIM leads to a system of differential–algebraic equations (DAEs) \citep{19Ahmed_2021,18Ahmed_2023}, which impose algebraic constraints on the solution and therefore require carefully constructed consistent initial conditions, complicating their numerical integration. Additional limitations of the initial CCNIM formulation include its restriction to two- or higher-dimensional problems, rendering it inapplicable to one-dimensional cases, as well as the relatively complex implementation of Neumann boundary conditions \citep{19Ahmed_2021,18Ahmed_2023}. Following these observations, a modified version of CCNIM, termed MCCNIM, was developed for multidimensional linear convection–diffusion equations, in which the temporal domain is treated as an additional nodal dimension \citep{17Ahmed_2024}. This reformulation yields a fully algebraic system, rather than a system of DAEs, thereby overcoming the limitations of the original CCNIM \citep{18Ahmed_2023}. The performance of MCCNIM \citep{17Ahmed_2024} has been evaluated in the literature through comparisons with the original CCNIM \citep{18Ahmed_2023} and traditional NIM \citep{22Rizwan_uddin_1997} formulations. These studies demonstrate that the modified formulation preserves the solution accuracy of CCNIM, with both approaches exhibiting the same order of accuracy and yielding nearly identical numerical results. It has also been shown that CCNIM, when formulated as a system of DAEs, is particularly well suited for problems requiring enhanced temporal accuracy or involving higher-order time derivatives, such as wave equations. In contrast, MCCNIM offers greater flexibility for problems ranging from one to higher spatial dimensions and allows for a simpler implementation of Neumann boundary conditions. Moreover, it has been established that, by treating the temporal derivative explicitly, the MCCNIM framework, together with its flux definition, can be easily reformulated as a CCNIM approach. Overall, both methods have been shown to be equally effective, with their suitability primarily determined by the specific requirements of the problem under consideration. 

The MCCNIM was subsequently extended to multidimensional nonlinear Burgers’ equations, and its results were compared with those obtained using traditional NIM formulations \citep{23Rizwan_uddin_1997,25Wescott_2001} and other conventional numerical schemes \citep{ahmed2025efficient}. These comparisons demonstrated that MCCNIM is not only accurate but also computationally more efficient than traditional NIM when applied to nonlinear problems. In the existing MCCNIM formulation, the nonlinear convection term is linearized by replacing the convective velocity with the node-averaged velocity evaluated at the current time step. Consequently, the coefficients of the resulting discretized algebraic system depend on these current-time velocities and must be updated at every iteration within each time step, leading to increased computational cost. Drawing on the underlying idea of the existing M$^2$NIM approach \citep{25Wescott_2001}, this work proposes an improved formulation of CCNIM, hereafter referred to as RCCNIM, which preserves accuracy comparable to the existing MCCNIM for the Burgers’ problem while achieving a significant reduction in computational time. Rather than directly substituting the convective velocity at the current time step, the proposed approach augments the original governing equation by adding an identical convection term to both sides of the PDE, with the velocity evaluated at the previous time step. The analytical solution is then constructed using this added term, resulting in a formulation whose coefficients depend only on previously computed velocities and therefore need to be evaluated once per time step. This strategy yields a substantial reduction in computational expense compared to the existing MCCNIM. To validate the proposed scheme, we have solved various one- and two-dimensional Burgers' problems with available analytical or referenced benchmark solutions. Additionally, a comprehensive comparison of the proposed scheme with previously published results is presented to demonstrate its effectiveness.

The organization of the manuscript is outlined as follows. \sect\ref{sec:2} presents the formulation of the RCCNIM scheme for the two-dimensional, time-dependent Burgers’ equations. In \sect\ref{sec:3}, numerical experiments are provided to assess the accuracy and efficiency of the proposed method. Finally,  \sect\ref{sec:4} summarizes the key findings along with concluding remarks and discussions.

\section{RCCNIM Formulation}
\label{sec:2}
%%%%%%%%%%%%%%%%%%%%%%%%%%%%
%%
%
%
%
%\subsection{One-dimensional Burgers’ equation}
%\label{sec:2.1}
%
RCCNIM has been applied to the two-dimensional, time-dependent, Burgers' equation with sources
\begin{gather}
	\frac{\partial u(x,y,t)}{\partial t}+u\left(x,y,t\right)\frac{\partial u\left(x,y,t\right)}{\partial x}+v\left(x,y,t\right)\frac{\partial u\left(x,y,t\right)}{\partial y}=\frac{1}{Re}\left(\frac{\partial^2u\left(x,y,t\right)}{\partial x^2}+\frac{\partial^2u\left(x,y,t\right)}{\partial y^2}\right)+f_x(x,y,t) 
\label{eq:001}
\end{gather}
\begin{gather}
	\frac{\partial v(x,y,t)}{\partial t}+u\left(x,y,t\right)\frac{\partial v\left(x,y,t\right)}{\partial x}+v\left(x,y,t\right)\frac{\partial v\left(x,y,t\right)}{\partial y}=\frac{1}{Re}\left(\frac{\partial^2v\left(x,y,t\right)}{\partial x^2}+\frac{\partial^2v\left(x,y,t\right)}{\partial y^2}\right)+f_y(x,y,t) 
\label{eq:002}
\end{gather}
Here, $u(x,y,t)$ and $v(x,y,t)$ represent the velocity components in the $x$- and $y$-directions, respectively; $Re$ denotes the Reynolds number; and $f_x(x,y,t)$ and $f_y(x,y,t)$ represent the corresponding source terms.
\begin{figure}[!b]
\centering
\vfill

\begin{minipage}{0.6\linewidth}
	\centering
	\includegraphics[width=\linewidth]{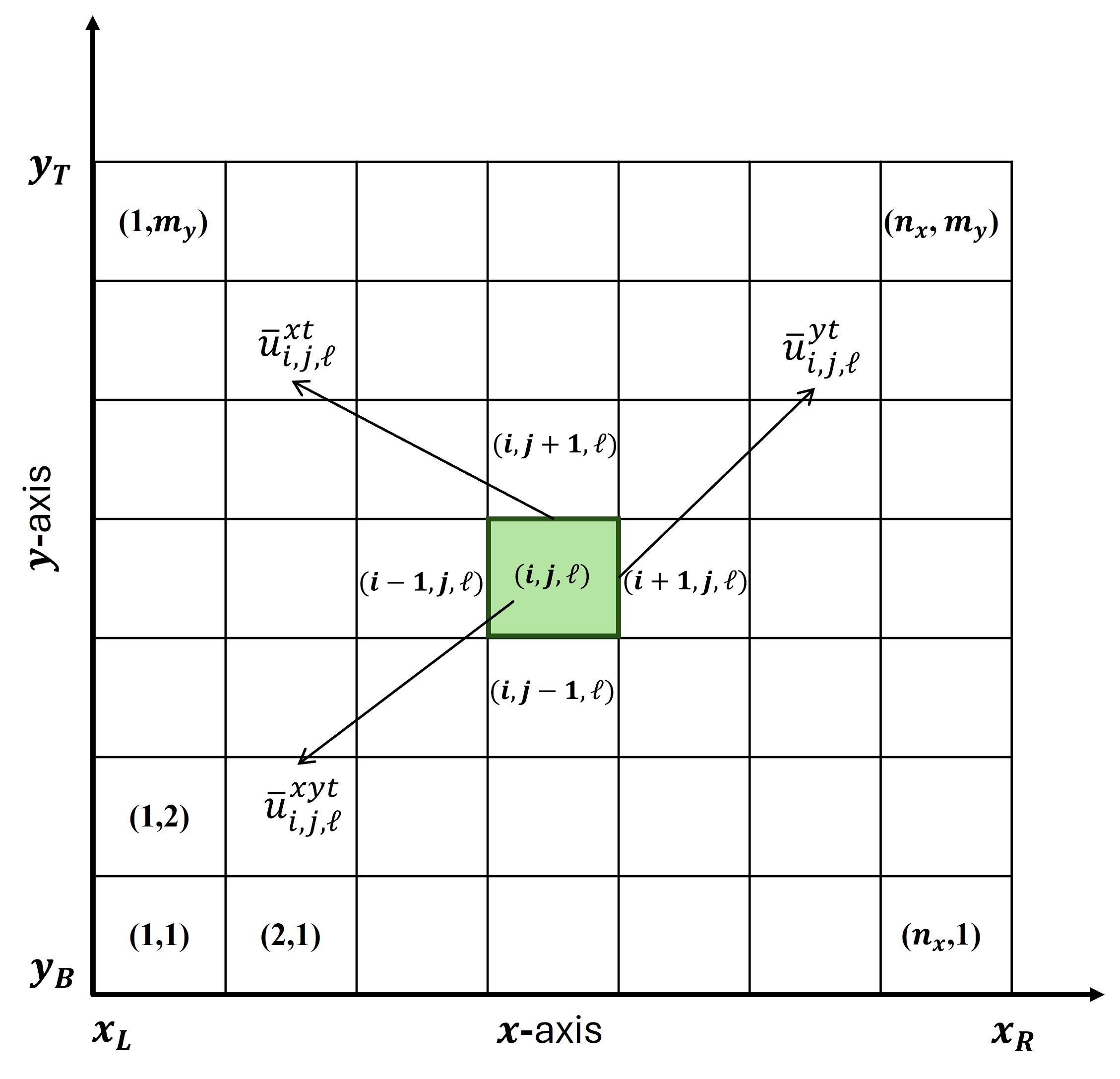}
	
	\vspace{0.2cm}
	{\Large (a)}
\end{minipage}

\vspace{0.5cm}

\begin{minipage}{0.8\linewidth}
	\centering
	\includegraphics[width=\linewidth]{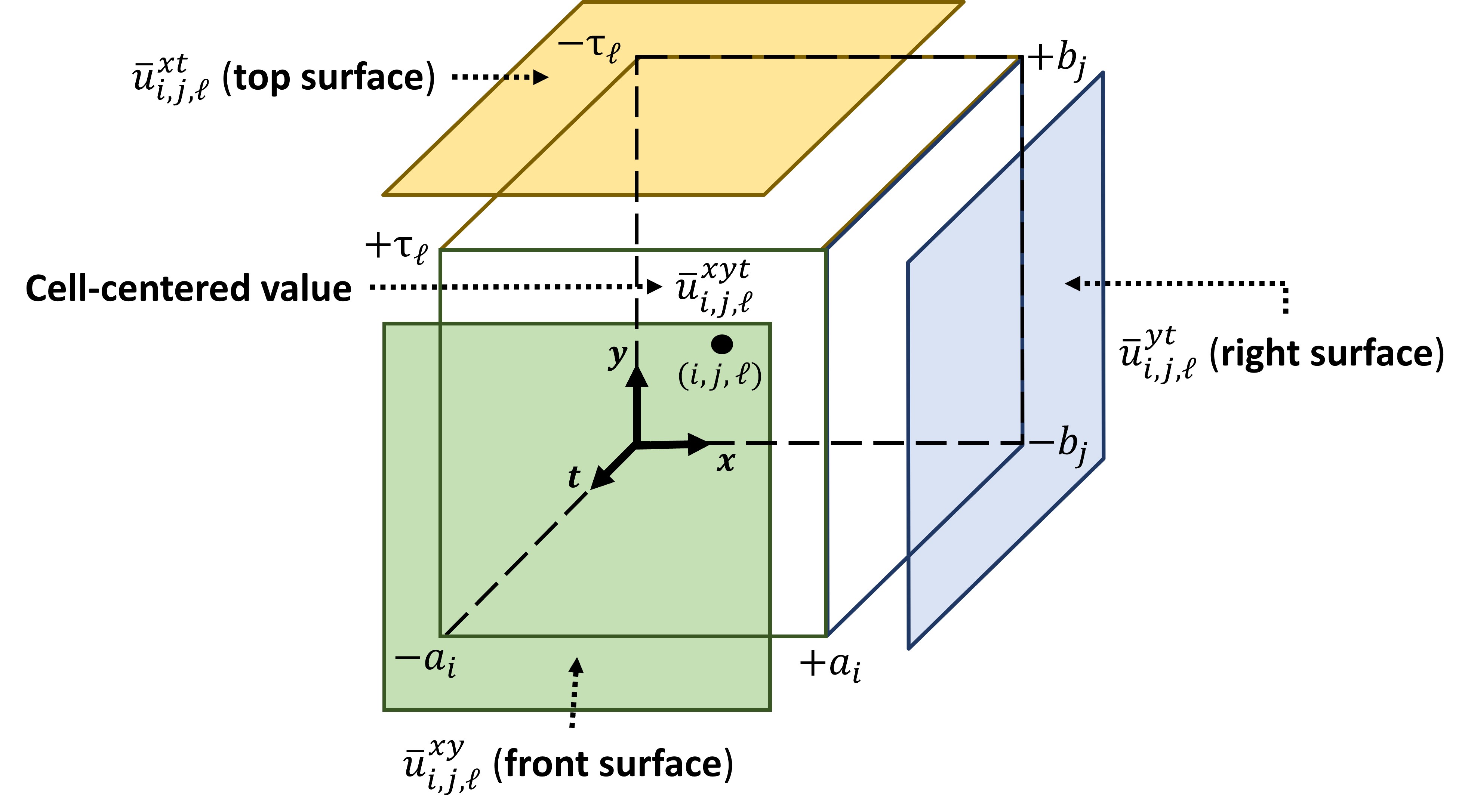}
	
	\vspace{0.2cm}
	{\Large (b)}
\end{minipage}

\vfill

\caption{Schematic representation of the computational grid for a two-dimensional RCCNIM problem. (a) Domain discretization. (b) Local coordinate framework and transverse-averaged velocities corresponding to cell $(i,j,\ell)$. \citep{ahmed2025efficient}}
\label{fig:011}

\end{figure}
RCCNIM is based on the same fundamental principle as the traditional NIM, with certain distinguishing features in its formulation. Like NIM, it is a semi-analytical numerical technique designed to solve partial differential equations (PDEs) efficiently and accurately, even on relatively coarse computational grids. Similar to conventional approaches such as the finite difference, finite volume, and finite element methods, the computational domain in RCCNIM is discretized into smaller elements referred to as nodes, rectangular or square nodes in one dimension and cuboid-like elements in two dimensions. However, RCCNIM is distinguished by its use of the transverse integration process (TIP), through which the governing PDEs are converted into a system of analytically solvable ordinary differential equations (ODEs) within each node.

In the two-dimensional formulation governed by \eqn\eqref{eq:001} and \eqn\eqref{eq:002}, the space--time domain is partitioned into finite cuboidal elements, hereafter referred to as nodes or cells, as depicted in panel (a) of \fig\ref{fig:011}. These nodes are indexed by $i$, $j$, and $l$ along the $x$-, $y$-, and temporal ($t$) directions, respectively. Each node $(i,j,l)$ is associated with dimensions $\Delta x \times \Delta y \times \Delta t$, where $\Delta x = 2a_i$, $\Delta y = 2b_j$, and $\Delta t = 2\tau_l$, with the origin located at the center of the node as illustrated in panel (b) of \fig\ref{fig:011}. Upon discretization, the transverse integration process is applied to each node, converting the governing PDEs into a corresponding set of ODEs, which form the basis for further analytical treatment within the nodal framework.

Prior to applying the transverse integration process (TIP) for RCCNIM to transform the system of PDEs into a set of ODEs, the original two-dimensional Burgers’ equations described by \eqn\eqref{eq:001} and \eqn\eqref{eq:002} are first recast in the local coordinate system, as presented below:
 \begin{gather}
	 \frac{\partial u}{\partial t}+u^p\frac{\partial u}{\partial x} + v^p\frac{\partial u}{\partial y} = \frac{1}{Re}\left[\frac{\partial^2u}{\partial x^2} + \frac{\partial^2u}{\partial y^2} \right]  - (u^0-u^p)\frac{\partial u}{\partial x} - (v^0-v^p)\frac{\partial u}{\partial y} + f_x(x,y,t) 
\label{eq:003}
\end{gather}
\begin{gather}
	\frac{\partial v}{\partial t}+u^p\frac{\partial v}{\partial x} + v^p\frac{\partial v}{\partial y} = \frac{1}{Re}\left[\frac{\partial^2v}{\partial x^2} + \frac{\partial^2v}{\partial y^2} \right]  - (u^0-u^p)\frac{\partial v}{\partial x} - (v^0-v^p)\frac{\partial v}{\partial y} + f_y(x,y,t) 
\label{eq:004}
\end{gather}
Here, ${\bar{u}}^p$ and ${\bar{v}}^p$ denote the cell-averaged convective velocities at the previous time level, whereas ${\bar{u}}^0$ and ${\bar{v}}^0$ represent the corresponding cell-averaged velocities at the current time step. In the earlier formulation of MCCNIM for nonlinear Burgers' equations, only the convective velocity at the current time level was incorporated. As a result, this velocity appeared within the coefficients of the resulting algebraic system, requiring their re-evaluation at every time step and consequently increasing the computational cost. In contrast, the present formulation, as given in \eqns\eqref{eq:003} and \eqref{eq:004}, introduces the convective term based on the cell-averaged velocity from the previous time step by adding it to both sides of the governing equations (\eqns\eqref{eq:001} and \eqref{eq:002}), leading to the modified forms in \eqns\eqref{eq:003} and \eqref{eq:004}. A similar strategy has been adopted in the surface-averaged NIM framework, commonly referred to as M$^2$NIM \citep{25Wescott_2001}. The motivation behind this reformulation is to express the final algebraic system in terms of coefficients that depend only on the cell-averaged velocities from the previous time level. This allows the coefficients to be evaluated just once per time step, thereby significantly reducing computational effort. The definitions of cell-averaged velocities are provided in \sect\ref{sec:2.0.1}.

The modified Burgers' equations presented in \eqns\eqref{eq:003} and \eqref{eq:004} are processed through transverse integration in each direction, resulting in a system of six ordinary differential equations, as discussed below. First, applying the transverse-integration operator 
$\frac{1}{4 a_i b_j} \int_{-a_i}^{+a_i} \int_{-b_j}^{+b_j} \, \mathrm{d}y \, \mathrm{d}x$ 
to \eqns\eqref{eq:003} and \eqref{eq:004}, respectively, leads to the first two ODEs
\begin{gather}
	\frac{\mathrm{d}{\bar{u}}^{xy}(t)}{\mathrm{d}t}={\bar{S}}_{1}^{xy}(t) 
\label{eq:005}
\end{gather}
\begin{gather}
	\frac{\mathrm{d}{\bar{v}}^{xy}(t)}{\mathrm{d}t}={\bar{S}}_{4}^{xy}(t)
\label{eq:006}
\end{gather}
where the cell-specific subscripts ($i,j,l$) on independent variables have been omitted and ${\bar{u}}^{xy}(t)$ and ${\bar{v}}^{xy}(t)$ represent the space averaged, time-dependent velocities which are defined as
\begin{gather}
	{\bar{u}}^{xy}\left(t\right)=\frac{1}{4 a_i b_j} \int_{-a_i}^{+a_i} \int_{-b_j}^{+b_j} u(x,y,t)\, \mathrm{d}y \, \mathrm{d}x 
\label{eq:007}
\end{gather}
\begin{gather}
	{\bar{v}}^{xy}\left(t\right)=\frac{1}{4 a_i b_j} \int_{-a_i}^{+a_i} \int_{-b_j}^{+b_j} v(x,y,t)\, \mathrm{d}y \, \mathrm{d}x 
\label{eq:008}
\end{gather}
In a similar fashion, the remaining transverse-averaged velocities, namely ${\bar{u}}^{xt}(y)$, ${\bar{v}}^{xt}(y)$, ${\bar{u}}^{yt}(x)$, and ${\bar{v}}^{yt}(x)$, can also be defined for later use.
Moreover, the quantities ${\bar{S}}_{1}^{xy}(t)$ and ${\bar{S}}_{4}^{xy}(t)$ represent pseudo-source terms that collect contributions which are not easily integrable and are given as
\begin{gather}
	\begin{split}
		{\bar{S}}_{1}^{xy}(t)=\frac{1}{4a_ib_j}\int_{-a_i}^{+a_i}\int_{-b_j}^{+b_j}\Bigg[\frac{1}{Re}\left(\frac{\partial^2{{u}}}{\partial x^2}+\frac{\partial^2{{u}}}{\partial y^2}\right)-{{u}}^0\frac{\partial{{u}}}{\partial x}-{{v}}_{}^0\frac{\partial{{u}}}{\partial y}+{{f}}_{x}\Bigg]\mathrm{d}y\mathrm{d}x
	\end{split}
\label{eq:009}
\end{gather}
\begin{gather}
	\begin{split}
		{\bar{S}}_{4}^{xy}(t)=\frac{1}{4a_ib_j}\int_{-a_i}^{+a_i}\int_{-b_j}^{+b_j}\Bigg[\frac{1}{Re}\left(\frac{\partial^2{{v}}}{\partial x^2}+\frac{\partial^2{{v}}}{\partial y^2}\right)-{{u}}^0\frac{\partial{{v}}}{\partial x}-{{v}}_{}^0\frac{\partial{{v}}}{\partial y}+{{f}}_{y}\Bigg]\mathrm{d}y\mathrm{d}x
	\end{split}
\label{eq:010}
\end{gather}
In a similar manner, applying transverse integration to \eqns\eqref{eq:003} and \eqref{eq:004} along the remaining directions yields four additional ordinary differential equations, given as follows:
\begin{gather}
	\frac{1}{Re}\frac{\mathrm{d}^2{\bar{u}}^{xt}\left(y\right)}{\mathrm{d}y^2}-{\bar{v}}_{}^p\frac{\mathrm{d}{\bar{u}}^{xt}(y)}{\mathrm{d}y}={\bar{S}}_{2}^{xt}(y) 
\label{eq:011}
\end{gather}
\begin{gather}
	\frac{1}{Re}\frac{\mathrm{d}^2{\bar{u}}^{yt}(x)}{\mathrm{d}x^2}-{\bar{u}}_{}^p\frac{\mathrm{d}{\bar{u}}^{yt}(x)}{\mathrm{d}x}={\bar{S}}_{3}^{yt}(x) 
\label{eq:012}
\end{gather}
\begin{gather}
	\frac{1}{Re}\frac{\mathrm{d}^2{\bar{v}}^{xt}\left(y\right)}{\mathrm{d}y^2}-{\bar{v}}_{}^p\frac{\mathrm{d}{\bar{v}}^{xt}(y)}{\mathrm{d}y}={\bar{S}}_{5}^{xt}(y)
\label{eq:013}
\end{gather}
\begin{gather}
	\frac{1}{Re}\frac{\mathrm{d}^2{\bar{v}}^{yt}(x)}{\mathrm{d}x^2}-{\bar{u}}_{}^p\frac{\mathrm{d}{\bar{v}}^{yt}(x)}{\mathrm{d}x}={\bar{S}}_{6}^{yt}(x)  
\label{eq:014}
\end{gather}
and the corresponding pseudo-source terms are given by:
\begin{gather}
	\begin{split}
		{\bar{S}}_{2}^{xt}(y)=\frac{1}{4a_i\tau_l}\int_{-\tau_l}^{+\tau_l}\int_{-a_i}^{+a_i}\Bigg[{\frac{\partial{{u}}^{}}{\partial t}}  + {{u}}_{}^0\frac{\partial{{u}}^{}}{\partial x} + (v^0-v^p)\frac{\partial {{u}}^{}}{\partial y}  - \frac{1}{Re}\frac{\partial^2{{u}}^{}}{\partial x^2}-{{f}}_{x}^{} \Bigg]\mathrm{d}x\mathrm{d}t 
	\end{split}
\label{eq:015}
\end{gather}
\begin{gather}
	\begin{split}
		{\bar{S}}_{3}^{yt}(x)=\frac{1}{4b_j\tau_l}\int_{-\tau_l}^{+\tau_l}\int_{-b_j}^{+b_j}\Bigg[{\frac{\partial{{u}}^{}}{\partial t}}  + (u^0-u^p)\frac{\partial{{u}}^{}}{\partial x} + v^0\frac{\partial {{u}}^{}}{\partial y}  - \frac{1}{Re}\frac{\partial^2{{u}}^{}}{\partial y^2}-{{f}}_{x}^{} \Bigg]\mathrm{d}y\mathrm{d}t 
	\end{split}
\label{eq:016}
\end{gather}
\begin{gather}
	\begin{split}
		{\bar{S}}_{5}^{xt}(y)=\frac{1}{4a_i\tau_l}\int_{-\tau_l}^{+\tau_l}\int_{-a_i}^{+a_i}\Bigg[{\frac{\partial{{v}}^{}}{\partial t}}  + {{u}}_{}^0\frac{\partial{{v}}^{}}{\partial x} + (v^0-v^p)\frac{\partial {{v}}^{}}{\partial y}  - \frac{1}{Re}\frac{\partial^2{{v}}^{}}{\partial x^2}-{{f}}_{y}^{} \Bigg]\mathrm{d}x\mathrm{d}t 
	\end{split}
\label{eq:017}
\end{gather}
\begin{gather}
	\begin{split}
		{\bar{S}}_{6}^{yt}(x)=\frac{1}{4b_j\tau_l}\int_{-\tau_l}^{+\tau_l}\int_{-b_j}^{+b_j}\Bigg[{\frac{\partial{{v}}^{}}{\partial t}}  + (u^0-u^p)\frac{\partial{{v}}^{}}{\partial x} + v^0\frac{\partial {{v}}^{}}{\partial y}  - \frac{1}{Re}\frac{\partial^2{{v}}^{}}{\partial y^2}-{{f}}_{y}^{} \Bigg]\mathrm{d}y\mathrm{d}t 
	\end{split}
\label{eq:018}
\end{gather}
It should be noted that no approximations have been introduced in these terms at this stage. All subsequent approximations are confined to the pseudo-source terms through a Legendre polynomial expansion, and the overall accuracy of the scheme depends on this treatment. In particular, the order of accuracy is governed by the level at which the expansion of the pseudo-source terms are truncated; for instance, second-order accuracy is obtained by truncating at the zeroth-order term, while higher-order accuracy can be achieved by retaining additional terms in the expansion. 

In the present study, the primary objective is to develop a computationally efficient scheme that reduces execution time compared to the existing approach, rather than pursuing higher-order accuracy. Accordingly, the pseudo-source terms are truncated at the zeroth-order level, which effectively yields constant right-hand sides in the resulting ordinary differential equations (\eqns\eqref{eq:005}-\eqref{eq:006}, \eqref{eq:011}-\eqref{eq:014}). The simplified form of the ODEs after truncation is given as follows:  
\begin{gather}
	\frac{\mathrm{d}{\bar{u}}^{xy}\left(t\right)}{\mathrm{d}t}={\bar{S}}_{1}^{xyt} 	
\label{eq:019}
\end{gather}
\begin{gather}
	\frac{1}{Re}\frac{\mathrm{d}^2{\bar{u}}^{xt}\left(y\right)}{\mathrm{d}y^2}-{{v}}_{}^p\frac{\mathrm{d}{\bar{u}}^{xt}(y)}{\mathrm{d}y}={\bar{S}}_{2}^{xyt} 
\label{eq:020}
\end{gather}
\begin{gather}
	\frac{1}{Re}\frac{\mathrm{d}^2{\bar{u}}^{yt}(x)}{\mathrm{d}x^2}-{{u}}_{}^p\frac{\mathrm{d}{\bar{u}}^{yt}(x)}{\mathrm{d}x}={\bar{S}}_{3}^{xyt} 
\label{eq:021}
\end{gather}
\begin{gather}
	\frac{\mathrm{d}{\bar{v}}^{xy}\left(t\right)}{\mathrm{d}t}={\bar{S}}_{4}^{xyt}
\label{eq:022}
\end{gather}
\begin{gather}
	\frac{1}{Re}\frac{\mathrm{d}^2{\bar{v}}^{xt}\left(y\right)}{\mathrm{d}y^2}-{{v}}_{}^p\frac{\mathrm{d}{\bar{v}}^{xt}(y)}{\mathrm{d}y}={\bar{S}}_{5}^{xyt}
\label{eq:023}
\end{gather}
\begin{gather}
	\frac{1}{Re}\frac{\mathrm{d}^2{\bar{v}}^{yt}(x)}{\mathrm{d}x^2}-{{u}}_{}^p\frac{\mathrm{d}{\bar{v}}^{yt}(x)}{\mathrm{d}x}={\bar{S}}_{6}^{xyt}  
\label{eq:024}
\end{gather}
It is worth noting that, in \eqns\eqref{eq:019}-\eqref{eq:024}, the right-hand sides reduce to constant terms as a result of truncating the expansion at the zeroth order. Consequently, the resulting ordinary differential equations can now be solved analytically. Solving \eqn\eqref{eq:019} and applying the nodal boundary condition at $t={+\tau}_{\ell}$, namely ${\bar{u}}^{xy}(t)={\bar{u}}^{xy}\left({+\tau}_{\ell}\right)={\bar{u}}_{i,j,\ell}^{xy}$, we obtain
\begin{gather}
	{\bar{u}}^{xy}\left(t\right)={\bar{S}}_{1i,j,\ell}^{xyt}\left(t-\tau_{\ell}\right)+{\bar{u}}_{i,j,\ell}^{xy}
\label{eq:025}
\end{gather}
To compute the cell-centered value, \eqn\eqref{eq:025} is further integrated in the remaining direction (i.e., the $t$-direction) by applying the operator $\frac{1}{2{\tau}_{\ell}}\int_{-{\tau}_{\ell}}^{+{\tau}_{\ell}}\mathrm{d}t$. This yields
\begin{gather}
	{\bar{u}}_{i,j,\ell}^{xyt}=-\tau_{\ell}{\bar{S}}_{1i,j,\ell}^{xyt}+{\bar{u}}_{i,j,\ell}^{xy} 
\label{eq:026}
\end{gather}
Here, ${\bar{u}}_{i,j,\ell}^{xyt}$ denotes the cell-centered variable. It is important to note that the pseudo source terms are already treated as cell-centered quantities, since, upon truncation to zeroth order, they reduce to constants. Consequently, ${\bar{S}}_{1i,j,\ell}^{xyt}$ is also a cell-centered quantity. Thus, in \eqn\eqref{eq:026}, the only remaining surface-averaged term is ${\bar{u}}_{i,j,\ell}^{xy}$. In order to derive a fully cell-centered scheme, all surface-averaged quantities must be eliminated from the final system of algebraic equations. Therefore, \eqn\eqref{eq:026} is rearranged to express the surface-averaged value in terms of the cell-centered variables.
\begin{gather}
	{\bar{u}}_{i,j,\ell}^{xy}={\bar{u}}_{i,j,\ell}^{xyt}+\tau_{\ell}{\bar{S}}_{1i,j,\ell}^{xyt} 
\label{eq:027}
\end{gather}
Similarly, the expression for ${\bar{v}}_{i,j,\ell}^{xy}$ is obtained by analytically solving the corresponding ODE, \eqn\eqref{eq:022}, and is given by 
\begin{gather}
\begin{split}
		{\bar{v}}_{i,j,\ell}^{xy}={\bar{v}}_{i,j,\ell}^{xyt}+\tau_{\ell}{\bar{S}}_{4i,j,\ell}^{xyt}
\end{split}
\label{eq:028}
\end{gather}
Up to this point, only two ODEs have been utilized, while the solutions of the remaining four ODEs (i.e., \eqns\eqref{eq:020}, \eqref{eq:021}, \eqref{eq:023}, and \eqref{eq:024}) are yet to be obtained. Since all four ODEs share an identical structure, only the solution of one representative ODE (i.e., \eqn\eqref{eq:021}) is presented here; the solutions of the others can be derived analogously. To solve \eqn\eqref{eq:021}, the surface-averaged flux $\left({\bar{J}}^{yt}_x\left(x\right)=-\frac{1}{Re}\frac{\mathrm{d}{\bar{u}}^{yt}(x)}{\mathrm{d}x}\right)$ and surface-averaged variable (${\bar{u}}^{yt}(x)$) at the interface of the adjacent node are employed as local nodal boundary conditions. Specifically, solving \eqn\eqref{eq:021} for the node $(i,j,\ell)$ and applying the nodal boundary conditions at the surface $x = {+a}_i$, namely the averaged velocity ${\bar{u}}^{yt}(x) = {\bar{u}}^{yt}(+a_i)$ and the averaged flux ${\bar{J}}^{yt}_x(x) = {\bar{J}}^{yt}_x(+a_i)$, yields
\begin{gather}
	\begin{split}
\left[{\bar{u}}^{yt}\left(x\right)\right]_{i,j,\ell}={\bar{u}}^{yt}\left({+a}_i\right)-\frac{-1+e^{{\left(x-a_i\right)\bar{u}}_{i,j,\ell}^pRe}}{\ {\bar{u}}_{i,j,\ell}^p}{\bar{J}}^{yt}_x\left(+a_i\right) \\
+\frac{\left(-1+e^{{\left(x-a_i\right)\bar{u}}_{i,j,\ell}^pRe}-(x-a_i){\bar{u}}_{i,j,\ell}^pRe\right)}{\left({\bar{u}}_{i,j,\ell}^p\right)^2Re}{\bar{S}}_{3i,j,\ell}^{xyt} 
	\end{split}
\label{eq:029}
\end{gather}
It is to be noted that the subscript $x$ in ${\bar{J}}^{yt}_x$ denotes the flux associated with the $x$-momentum equation. The same ODE (i.e., \eqn\eqref{eq:021}) is again solved for the neighboring node $(i+1,j,\ell)$ by applying the nodal boundary conditions at the interface $x = {-a}_{i+1}$. Specifically, the averaged variable and flux are prescribed as ${\bar{u}}^{yt}(x) = {\bar{u}}^{yt}({-a}_{i+1})$ and ${\bar{J}}^{yt}_x(x) = {\bar{J}}^{yt}_x({-a}_{i+1})$, respectively. This yields
\begin{gather}
	\begin{split}
    \left[{\bar{u}}^{yt}\left(x\right)\right]_{i+1,j,\ell}={\bar{u}}^{yt}\left({-a}_{i+1}\right)-\frac{-1+e^{\left(x+a_{i+1}\right){\bar{u}}_{i+1,j,\ell}^pRe}}{\ {\bar{u}}_{i+1,j,\ell}^p}{\bar{J}}^{yt}_x\left(-a_{i+1}\right) \\ +\frac{\left(-1+e^{\left(x+a_{i+1}\right){\bar{u}}_{i+1,j,\ell}^pRe}-{(x+a_{i+1})\bar{u}}_{i+1,j,\ell}^pRe\right)}{\left({\bar{u}}_{i+1,j,\ell}^p\right)^2Re}{\bar{S}}_{3i+1,j,\ell}^{xyt} 
	\end{split}
\label{eq:030}
\end{gather}
To compute the cell-centered values, \eqns\eqref{eq:029} and \eqref{eq:030} are further integrated in the remaining direction (i.e., the $x$-direction). This is achieved by applying the operators $\frac{1}{2{a}_{i}}\int_{-{a}_{i}}^{+{a}_{i}}\mathrm{d}x$ to \eqn\eqref{eq:029} and $\frac{1}{2{a}_{i+1}}\int_{-{a}_{i+1}}^{+{a}_{i+1}}\mathrm{d}x$ to \eqn\eqref{eq:030}, respectively. This yields
\begin{gather}
	\begin{split}
	{\bar{u}}_{i,j,\ell}^{xyt}&={\bar{u}}^{yt}\left({+a}_i\right)-\frac{\left(1-e^{-{Ru}_{i,j,\ell}}-2a_i{\bar{u}}_{i,j,\ell}^pRe\right)}{2a_i\left({\bar{u}}_{i,j,\ell}^p\right)^2{Re}}{\bar{J}}^{yt}_x\left(+a_i\right)\\&+\frac{\left(1-e^{-{Ru}_{i,j,\ell}}+2a_i{\bar{u}}_{i,j,\ell}^pRe\left(-1+a_i{\bar{u}}_{i,j,\ell}^pRe\right)\right)}{2a_i\left({\bar{u}}_{i,j,\ell}^p\right)^3{Re}^2}{\bar{S}}_{3i,j,\ell}^{xyt} 
	\end{split}
\label{eq:031}
\end{gather}
\begin{gather}
	\begin{split}
		{\bar{u}}_{i+1,j,\ell}^{xyt}&={\bar{u}}^{yt}\left({-a}_{i+1}\right)+\frac{\left(1-e^{{Ru}_{i+1,j,\ell}}-2a_{i+1}{\bar{u}}_{i+1,j,\ell}^pRe\right)}{2a_{i+1}\left({\bar{u}}_{i+1,j,\ell}^{p}\right)^2{Re}}{\bar{J}}^{yt}_x\left(-a_{i+1}\right)\\ &+\frac{\left(-1+e^{{Ru}_{i+1,j,\ell}}-2a_{i+1}{\bar{u}}_{i+1,j,\ell}^pRe\left(1+a_{i+1}{\bar{u}}_{i+1,j,\ell}^pRe\right)\right)}{2a_{i+1}\left({\bar{u}}_{i+1,j,\ell}^p\right)^3{Re}^2}{\bar{S}}_{3i+1,j,\ell}^{xyt} 
	\end{split}
\label{eq:032}
\end{gather}
It is to be noted that in \eqns\eqref{eq:031} and \eqref{eq:032}, the exponential terms are expressed in terms of the corresponding local Reynolds number. For example, at the node $(i,j,\ell)$, the local Reynolds number is defined as ${Ru}_{i,j,\ell}=2a_i{\bar{u}}_{i,j,\ell}^{p}Re$. In addition, \eqns\eqref{eq:018} and \eqref{eq:019} may be reformulated into the following equivalent forms:
\begin{gather}
	{\bar{J}}^{yt}_x\left(+a_i\right)=A_{31}\left({\bar{u}}_{i,j,\ell}^{xyt}-{\bar{u}}^{yt}\left({+a}_i\right)\right)+A_{32}{\bar{S}}_{3i,j,\ell}^{xyt}  
\label{eq:033}
\end{gather}
\begin{gather}
{\bar{J}}^{yt}_x\left(-a_{i+1}\right)=A_{51,i+1}\left({\bar{u}}_{i+1,j,\ell}^{xyt}-{\bar{u}}^{yt}\left({-a}_{i+1}\right)\right)+A_{52,i+1}{\bar{S}}_{3i+1,j,\ell}^{xyt} 
\label{eq:034}
\end{gather}
where $A_{31}$, $A_{32}$, $A_{51}$ and $A_{52}$ are the coefficients that are dependent on the various parameters  such as $a_i$, ${Reu}_{i,j,\ell}$, ${\bar{u}}_{i,j,\ell}^{p}$. The explicit definitions of these coefficients are presented in \ref{app:A}. Next, continuity conditions are imposed at the shared interface between the adjacent nodes $(i,j,\ell)$ and $(i+1,j,\ell)$. In particular, continuity of the surface-averaged velocity requires
\begin{gather}
	{\bar{u}}^{yt}\left({+a}_i\right)={\bar{u}}^{yt}\left({-a}_{i+1}\right)={\bar{u}}_{i,j,\ell}^{yt} 
\label{eq:035}
\end{gather}
while continuity of the surface-averaged flux is expressed as
\begin{gather}
{\bar{J}}^{yt}_x\left(+a_i\right)={\bar{J}}^{yt}_x\left(-a_{i+1}\right)={\bar{J}}_{xi,j,\ell}^{yt} 
\label{eq:036}
\end{gather}
By enforcing these interface conditions and making use of \eqns\eqref{eq:033} and \eqref{eq:034}, an explicit expression for the interfacial surface-averaged velocity ${\bar{u}}_{i,j,\ell}^{yt}$ is obtained as
\begin{gather}
	{\bar{u}}_{i,j,\ell}^{yt}=\frac{A_{32}{\bar{S}}_{3i,j,\ell}^{xyt}-\ A_{52,i+1}{\bar{S}}_{3i+1,j,\ell}^{xyt}+\ A_{31}{\bar{u}}_{i,j,\ell}^{xyt}-\ A_{51,i+1}{\bar{u}}_{i+1,j,\ell}^{xyt}}{A_{31}-A_{51,i+1}} 
\label{eq:037}
\end{gather}
Substituting the expression for ${\bar{u}}_{i,j,\ell}^{yt}$ from \eqn\eqref{eq:037} into either of the flux relations, \eqns\eqref{eq:033} or \eqref{eq:034}, leads to an identical expression for ${\bar{J}}_{i,j,\ell}^{yt}$, thereby confirming that the continuity conditions are satisfied. The resulting expression for ${\bar{J}}_{i,j,\ell}^{yt}$ is given by
\begin{gather}
	{\bar{J}}_{xi,j,\ell}^{yt}=\frac{-A_{32}A_{51,i+1}{\bar{S}}_{3i,j,\ell}^{xyt}+\ {A_{31}A}_{52,i+1}{\bar{S}}_{3i+1,j,\ell}^{xyt}+\ A_{31}A_{51,i+1}({\bar{u}}_{i+1,j,\ell}^{xyt}-{\bar{u}}_{i,j,\ell}^{xyt})}{A_{31}-A_{51,i+1}}  
\label{eq:038}
\end{gather}
The remaining three ODEs, namely \eqns\eqref{eq:020}, \eqref{eq:023}, and \eqref{eq:024}, can be solved in an analogous manner to obtain the corresponding surface-averaged variables, which are given by
\begin{gather}
	\begin{split}
			{\bar{u}}_{i,j,\ell}^{xt}=\frac{B_{32}{\bar{S}}_{2i,j,\ell}^{xyt}-\ B_{52,j+1}{\bar{S}}_{2i,j+1,\ell}^{xyt}+\ B_{31}{\bar{u}}_{i,j,\ell}^{xyt}-\ B_{51,j+1}{\bar{u}}_{i,j+1,\ell}^{xyt}}{B_{31}-B_{51,j+1}} 
	\end{split}
\label{eq:039}
\end{gather}
\begin{gather}
\begin{split}
	{\bar{J}}_{xi,j,\ell}^{xt}=\frac{-B_{32}B_{51,j+1}{\bar{S}}_{2i,j,\ell}^{xyt}+\ {B_{31}B}_{52,j+1}{\bar{S}}_{2i,j+1,\ell}^{xyt}+\ B_{31}B_{51,j+1}({\bar{u}}_{i,j+1,\ell}^{xyt}-{\bar{u}}_{i,j,\ell}^{xyt})}{B_{31}-B_{51,j+1}} 
\end{split}
\label{eq:040}
\end{gather}
\begin{gather}
	\begin{split}
			{\bar{v}}_{i,j,\ell}^{xt}=\frac{B_{32}{\bar{S}}_{5i,j,\ell}^{xyt}-\ B_{52,j+1}{\bar{S}}_{5i,j+1,\ell}^{xyt}+\ B_{31}{\bar{v}}_{i,j,\ell}^{xyt}-\ B_{51,j+1}{\bar{v}}_{i,j+1,\ell}^{xyt}}{B_{31}-B_{51,j+1}}
	\end{split}
\label{eq:041}
\end{gather}
\begin{gather}
	\begin{split}
			{\bar{J}}_{yi,j,\ell}^{xt}=\frac{-B_{32}B_{51,j+1}{\bar{S}}_{5i,j,\ell}^{xyt}+\ {B_{31}B}_{52,j+1}{\bar{S}}_{5i,j+1,\ell}^{xyt}+\ B_{31}B_{51,j+1}({\bar{v}}_{i,j+1,\ell}^{xyt}-{\bar{v}}_{i,j,\ell}^{xyt})}{B_{31}-B_{51,j+1}}
	\end{split}
\label{eq:042}
\end{gather}
\begin{gather}
	\begin{split}
			{\bar{v}}_{i,j,\ell}^{yt}=\frac{A_{32}{\bar{S}}_{6i,j,\ell}^{xyt}-\ A_{52,i+1}{\bar{S}}_{6i+1,j,\ell}^{xyt}+\ A_{31}{\bar{v}}_{i,j,\ell}^{xyt}-\ A_{51,i+1}{\bar{v}}_{i+1,j,\ell}^{xyt}}{A_{31}-A_{51,i+1}}
	\end{split}
\label{eq:043}
\end{gather}
\begin{gather}
	\begin{split}
			{\bar{J}}_{yi,j,\ell}^{yt}=\frac{-A_{32}A_{51,i+1}{\bar{S}}_{6i,j,\ell}^{xyt}+\ {A_{31}A}_{52,i+1}{\bar{S}}_{6i+1,j,\ell}^{xyt}+\ A_{31}A_{51,i+1}({\bar{v}}_{i+1,j,\ell}^{xyt}-{\bar{v}}_{i,j,\ell}^{xyt})}{A_{31}-A_{51,i+1}}
	\end{split}
\label{eq:044}
\end{gather}
At this stage, solutions to all the ODEs have been obtained in terms of surface-averaged variables expressed through cell-centered quantities. To formulate a fully cell-centered scheme, it is necessary to eliminate all surface-averaged variables. This is accomplished in the final step by deriving appropriate constraint equations based on the definitions of the pseudo-source terms introduced earlier in \eqns\eqref{eq:009}-\eqref{eq:010}, \eqref{eq:015}-\eqref{eq:018} thereby closing the system of algebraic equations.
After performing averaging in all directions, equivalently, integrating each equation over the remaining coordinate directions, the pseudo-source terms defined in \eqns\eqref{eq:009}-\eqref{eq:010} and \eqref{eq:015}-\eqref{eq:018} can be expressed as
\begin{gather}
	\begin{split}
		{\bar{S}}_{1i,j,\ell}^{xyt}=\frac{1}{2a_i}\int_{-a_i}^{+a_i}\left(\frac{1}{Re}\frac{\mathrm{d}^2{\bar{u}}^{yt}\left(x\right)}{\mathrm{d}x^2}-{\bar{u}}_{i,j,\ell}^0\frac{\mathrm{d}{\bar{u}}^{yt}\left(x\right)}{\mathrm{d}x}\right)\mathrm{d}x+\frac{1}{2b_j}\int_{-b_j}^{+b_j}\left(\frac{1}{Re}\frac{\mathrm{d}^2{\bar{u}}^{xt}\left(y\right)}{\mathrm{d}y^2}-{\bar{v}}_{i,j,\ell}^0\frac{\mathrm{d}{\bar{u}}^{xt}\left(y\right)}{\mathrm{d}y}\right)\mathrm{d}y+{\bar{f}}_{xi,j,\ell}^{xyt}
	\end{split}
\label{eq:045}
\end{gather}
\begin{gather}
	\begin{split}
		{\bar{S}}_{2i,j,\ell}^{xyt}=\frac{1}{2\tau_{\ell}}\int_{-\tau_{\ell}}^{+\tau_{\ell}}{\frac{\mathrm{d}{\bar{u}}^{xy}\left(t\right)}{\mathrm{d}t}\mathrm{d}t}-\frac{1}{2a_i}\int_{-a_i}^{+a_i}\left(\frac{1}{Re}\frac{\mathrm{d}^2{\bar{u}}^{yt}\left(x\right)}{\mathrm{d}x^2}-{\bar{u}}_{i,j,\ell}^0\frac{\mathrm{d}{\bar{u}}^{yt}\left(x\right)}{\mathrm{d}x}\right)\mathrm{d}x+  ({\bar{v}}_{i,j,\ell}^0-{\bar{v}}_{i,j,\ell}^p)\frac{1}{2b_j}\int_{-b_j}^{+b_j}\frac{\mathrm{d} {\bar{u}}^{xt}\left(y\right)}{\mathrm{d} y}\mathrm{d}y - {\bar{f}}_{xi,j,\ell}^{xyt} 
	\end{split}
\label{eq:046}
\end{gather}
\begin{gather}
	\begin{split}
		{\bar{S}}_{3i,j,\ell}^{xyt}=\frac{1}{2\tau_{\ell}}\int_{-\tau_{\ell}}^{+\tau_{\ell}}{\frac{\mathrm{d}{\bar{u}}^{xy}\left(t\right)}{\mathrm{d}t}\mathrm{d}t}-\frac{1}{2b_j}\int_{-b_j}^{+b_j}\left(\frac{1}{Re}\frac{\mathrm{d}^2{\bar{u}}^{xt}\left(y\right)}{\mathrm{d}y^2}-{\bar{v}}_{i,j,\ell}^0\frac{\mathrm{d}{\bar{u}}^{xt}\left(y\right)}{\mathrm{d}y}\right)\mathrm{d}y + ({\bar{u}}_{i,j,\ell}^0-{\bar{u}}_{i,j,\ell}^p)\frac{1}{2a_i}\int_{-a_i}^{+a_i}\frac{\mathrm{d} {\bar{u}}^{yt}\left(x\right)}{\mathrm{d} x}\mathrm{d}x - {\bar{f}}_{xi,j,\ell}^{xyt}
	\end{split}
\label{eq:047}
\end{gather}
\begin{gather}
	\begin{split}
		{\bar{S}}_{4i,j,\ell}^{xyt}=\frac{1}{2a_i}\int_{-a_i}^{+a_i}\left(\frac{1}{Re}\frac{\mathrm{d}^2{\bar{v}}^{yt}\left(x\right)}{\mathrm{d}x^2}-{\bar{u}}_{i,j,\ell}^0\frac{\mathrm{d}{\bar{v}}^{yt}\left(x\right)}{\mathrm{d}x}\right)\mathrm{d}x+\frac{1}{2b_j}\int_{-b_j}^{+b_j}\left(\frac{1}{Re}\frac{\mathrm{d}^2{\bar{v}}^{xt}\left(y\right)}{\mathrm{d}y^2}-{\bar{v}}_{i,j,\ell}^0\frac{\mathrm{d}{\bar{v}}^{xt}\left(y\right)}{\mathrm{d}y}\right)\mathrm{d}y+{\bar{f}}_{yi,j,\ell}^{xyt} 
	\end{split}
\label{eq:048}
\end{gather}
\begin{gather}
	\begin{split}
		{\bar{S}}_{5i,j,\ell}^{xyt}=\frac{1}{2\tau_{\ell}}\int_{-\tau_{\ell}}^{+\tau_{\ell}}{\frac{\mathrm{d}{\bar{v}}^{xy}\left(t\right)}{\mathrm{d}t}\mathrm{d}t}-\frac{1}{2a_i}\int_{-a_i}^{+a_i}\left(\frac{1}{Re}\frac{\mathrm{d}^2{\bar{v}}^{yt}\left(x\right)}{\mathrm{d}x^2}-{\bar{u}}_{i,j,\ell}^0\frac{\mathrm{d}{\bar{v}}^{yt}\left(x\right)}{\mathrm{d}x}\right)\mathrm{d}x + ({\bar{v}}_{i,j,\ell}^0-{\bar{v}}_{i,j,\ell}^p)\frac{1}{2b_j}\int_{-b_j}^{+b_j}\frac{\mathrm{d} {\bar{v}}^{xt}\left(y\right)}{\mathrm{d} y}\mathrm{d}y - {\bar{f}}_{yi,j,\ell}^{xyt}
	\end{split}
\label{eq:049}
\end{gather}
\begin{gather}
	\begin{split}
		{\bar{S}}_{6i,j,\ell}^{xyt}=\frac{1}{2\tau_{\ell}}\int_{-\tau_{\ell}}^{+\tau_{\ell}}{\frac{\mathrm{d}{\bar{v}}^{xy}\left(t\right)}{\mathrm{d}t}\mathrm{d}t}-\frac{1}{2b_j}\int_{-b_j}^{+b_j}\left(\frac{1}{Re}\frac{\mathrm{d}^2{\bar{v}}^{xt}\left(y\right)}{\mathrm{d}y^2}-{\bar{v}}_{i,j,\ell}^0\frac{\mathrm{d}{\bar{v}}^{xt}\left(y\right)}{\mathrm{d}y}\right)\mathrm{d}y + ({\bar{u}}_{i,j,\ell}^0-{\bar{u}}_{i,j,\ell}^p)\frac{1}{2a_i}\int_{-a_i}^{+a_i}\frac{\mathrm{d} {\bar{v}}^{yt}\left(x\right)}{\mathrm{d} x}\mathrm{d}x - {\bar{f}}_{yi,j,\ell}^{xyt}
	\end{split}
\label{eq:050}
\end{gather}
By employing the definitions of the averaged fluxes, namely ${\bar{J}}^{yt}_x\left(x\right)=-\frac{1}{Re}\frac{\mathrm{d}{\bar{u}}^{yt}(x)}{\mathrm{d}x}$, ${\bar{J}}^{xt}_x\left(y\right)=-\frac{1}{Re}\frac{\mathrm{d}{\bar{u}}^{xt}(y)}{\mathrm{d}y}$, ${\bar{J}}^{yt}_y\left(x\right)=-\frac{1}{Re}\frac{\mathrm{d}{\bar{v}}^{yt}(x)}{\mathrm{d}x}$, and ${\bar{J}}^{xt}_y\left(y\right)=-\frac{1}{Re}\frac{\mathrm{d}{\bar{v}}^{xt}(y)}{\mathrm{d}y}$ and subsequently rearranging the pseudo-source terms, a set of six constraint equations can be derived in terms of the surface-averaged quantities as follows:
\begin{gather}
	\begin{split}
		{\bar{S}}_{1i,j,\ell}^{xyt}=-\frac{{\bar{J}}_{xi,j,\ell}^{yt}-{\bar{J}}_{xi-1,j,\ell}^{yt}}{2a_i}-{\bar{u}}_{i,j,\ell}^0\frac{{\bar{u}}_{i,j,\ell}^{yt}-{\bar{u}}_{i-1,j,\ell}^{yt}}{2a_i}-\frac{{\bar{J}}_{xi,j,\ell}^{xt}-{\bar{J}}_{xi,j-1,\ell}^{xt}}{2b_j}-{\bar{v}}_{i,j,\ell}^0\frac{{\bar{u}}_{i,j,\ell}^{xt}-{\bar{u}}_{i,j-1,\ell}^{xt}}{2b_j}+{\bar{f}}_{xi,j,\ell}^{xyt} 	
	\end{split}
\label{eq:051}
\end{gather}
\begin{gather}
	\begin{split}
		{\bar{S}}_{2i,j,\ell}^{xyt}=\frac{{\bar{u}}_{i,j,\ell}^{xy}-{\bar{u}}_{i,j,\ell-1}^{xy}}{2\tau_{\ell}}+\frac{{\bar{J}}_{xi,j,\ell}^{yt}-{\bar{J}}_{xi-1,j,\ell}^{yt}}{2a_i}+{\bar{u}}_{i,j,\ell}^0\frac{{\bar{u}}_{i,j,\ell}^{yt}-{\bar{u}}_{i-1,j,\ell}^{yt}}{2a_i} + ({\bar{v}}_{i,j,\ell}^0-{\bar{v}}_{i,j,\ell}^p)\frac{{\bar{u}}_{i,j,\ell}^{xt}-{\bar{u}}_{i,j-1,\ell}^{xt}}{2b_j} -{\bar{f}}_{xi,j,\ell}^{xyt}
	\end{split}
\label{eq:052}
\end{gather}
\begin{gather}
	\begin{split}
		{\bar{S}}_{3i,j,\ell}^{xyt}=\frac{{\bar{u}}_{i,j,\ell}^{xy}-{\bar{u}}_{i,j,\ell-1}^{xy}}{2\tau_{\ell}}+\frac{{\bar{J}}_{xi,j,\ell}^{xt}-{\bar{J}}_{xi,j-1,\ell}^{xt}}{2b_j}+{\bar{v}}_{i,j,\ell}^0\frac{{\bar{u}}_{i,j,\ell}^{xt}-{\bar{u}}_{i,j-1,\ell}^{xt}}{2b_j} + ({\bar{u}}_{i,j,\ell}^0-{\bar{u}}_{i,j,\ell}^p)\frac{{\bar{u}}_{i,j,\ell}^{yt}-{\bar{u}}_{i-1,j,\ell}^{yt}}{2a_i} - {\bar{f}}_{xi,j,\ell}^{xyt}
	\end{split}
\label{eq:053}
\end{gather}
\begin{gather}
	\begin{split}
		{\bar{S}}_{4i,j,\ell}^{xyt}=-\frac{{\bar{J}}_{yi,j,\ell}^{yt}-{\bar{J}}_{yi-1,j,\ell}^{yt}}{2a_i}-{\bar{u}}_{i,j,\ell}^0\frac{{\bar{v}}_{i,j,\ell}^{yt}-{\bar{v}}_{i-1,j,\ell}^{yt}}{2a_i}-\frac{{\bar{J}}_{yi,j,\ell}^{xt}-{\bar{J}}_{yi,j-1,\ell}^{xt}}{2b_j}-{\bar{v}}_{i,j,\ell}^0\frac{{\bar{v}}_{i,j,\ell}^{xt}-{\bar{v}}_{i,j-1,\ell}^{xt}}{2b_j}+{\bar{f}}_{yi,j,\ell}^{xyt}
	\end{split}
\label{eq:054}
\end{gather}
\begin{gather}
	\begin{split}
		{\bar{S}}_{5i,j,\ell}^{xyt}=\frac{{\bar{v}}_{i,j,\ell}^{xy}-{\bar{v}}_{i,j,\ell-1}^{xy}}{2\tau_{\ell}}+\frac{{\bar{J}}_{yi,j,\ell}^{yt}-{\bar{J}}_{yi-1,j,\ell}^{yt}}{2a_i}+{\bar{u}}_{i,j,\ell}^0\frac{{\bar{v}}_{i,j,\ell}^{yt}-{\bar{v}}_{i-1,j,\ell}^{yt}}{2a_i} + ({\bar{v}}_{i,j,\ell}^0-{\bar{v}}_{i,j,\ell}^p)\frac{{\bar{v}}_{i,j,\ell}^{xt}-{\bar{v}}_{i,j-1,\ell}^{xt}}{2b_j} -{\bar{f}}_{yi,j,\ell}^{xyt}
	\end{split}
\label{eq:055}
\end{gather}
\begin{gather}
	\begin{split}
		{\bar{S}}_{6i,j,\ell}^{xyt}=\frac{{\bar{v}}_{i,j,\ell}^{xy}-{\bar{v}}_{i,j,\ell-1}^{xy}}{2\tau_{\ell}}+\frac{{\bar{J}}_{yi,j,\ell}^{xt}-{\bar{J}}_{yi,j-1,\ell}^{xt}}{2b_j}+{\bar{v}}_{i,j,\ell}^0\frac{{\bar{v}}_{i,j,\ell}^{xt}-{\bar{v}}_{i,j-1,\ell}^{xt}}{2b_j} + ({\bar{u}}_{i,j,\ell}^0-{\bar{u}}_{i,j,\ell}^p)\frac{{\bar{v}}_{i,j,\ell}^{yt}-{\bar{v}}_{i-1,j,\ell}^{yt}}{2a_i} -{\bar{f}}_{yi,j,\ell}^{xyt}
	\end{split}
\label{eq:056}
\end{gather}
In addition, the original Burgers’ equations (\eqns\ref{eq:003} and \ref{eq:004}) are averaged over the control volume corresponding to node ($i,j,\ell$) by applying the operator $\left(\frac{1}{8a_ib_j\tau_\ell}\int_{-\tau_\ell}^{+\tau_\ell}\int_{-b_j}^{+b_j}\int_{-a_i}^{+a_i} \mathrm{d}t \mathrm{d}y \mathrm{d}x \right)$.
Upon invoking the definitions of the pseudo-source terms (\eqns\eqref{eq:019}--\eqref{eq:024}), this procedure leads to the formulation of two additional constraint relations, which can be expressed as follows.
\begin{gather}
	{\bar{S}}_{1i,j,\ell}^{xyt}={\bar{S}}_{2i,j,\ell}^{xyt}+{\bar{S}}_{3i,j,\ell}^{xyt} - ({\bar{u}}_{i,j,\ell}^0-{\bar{u}}_{i,j,\ell}^p)\frac{{\bar{u}}_{i,j,\ell}^{yt}-{\bar{u}}_{i-1,j,\ell}^{yt}}{2a_i} - ({\bar{v}}_{i,j,\ell}^0-{\bar{v}}_{i,j,\ell}^p)\frac{{\bar{u}}_{i,j,\ell}^{xt}-{\bar{u}}_{i,j-1,\ell}^{xt}}{2b_j}+{\bar{f}}_{xi,j,\ell}^{xyt} 
\label{eq:057}
\end{gather}
\begin{gather}
	{\bar{S}}_{4i,j,\ell}^{xyt}={\bar{S}}_{5i,j,\ell}^{xyt}+{\bar{S}}_{6i,j,\ell}^{xyt} - ({\bar{u}}_{i,j,\ell}^0-{\bar{u}}_{i,j,\ell}^p)\frac{{\bar{v}}_{i,j,\ell}^{yt}-{\bar{v}}_{i-1,j,\ell}^{yt}}{2a_i} - ({\bar{v}}_{i,j,\ell}^0-{\bar{v}}_{i,j,\ell}^p)\frac{{\bar{v}}_{i,j,\ell}^{xt}-{\bar{v}}_{i,j-1,\ell}^{xt}}{2b_j} + {\bar{f}}_{yi,j,\ell}^{xyt}  	
\label{eq:058}
\end{gather}
In the earlier formulation of the CCNIM scheme for the Burgers’ equations (MCCNIM), the surface-averaged quantities were directly substituted into the constraint equations to obtain the final system of algebraic equations. However, extending this procedure to the present RCCNIM framework proved to be considerably more involved. To address this difficulty, a modified strategy is adopted in the current work, wherein each term in the constraint equations is treated individually rather than substituting the surface-averaged expressions in a single step. The corresponding coefficients are then systematically constructed from these individually evaluated contributions. This reformulation offers several advantages. First, it provides a unified and straightforward pathway for developing both MCCNIM and RCCNIM schemes, as the individual terms are common to both formulations; moreover, it facilitates a more natural extension to the Navier–Stokes equations. Second, the resulting algebraic expressions are comparatively compact, reducing both the complexity and the number of coefficients that must be evaluated at each iteration. Finally, numerical observations indicate improved stability characteristics: in certain cases where the coefficients obtained via the previous approach fail to converge, the present formulation exhibits robust convergence behavior, highlighting the effectiveness of the proposed treatment.
In the set of four constraint equations (\eqns\ref{eq:051}-\ref{eq:053} and \ref{eq:057}) obtained from the $x$-momentum equation, the constituent terms are identified and treated individually, as outlined below:
\begin{gather}
	\begin{split}
		\frac{{\bar{J}}_{xi,j,\ell}^{yt}-{\bar{J}}_{xi-1,j,\ell}^{yt}}{2a_i} =  C_{31}{\bar{S}}_{3i,j,\ell}^{xyt}+C_{32}{\bar{S}}_{3i-1,j,\ell}^{xyt}+C_{33}{\bar{S}}_{3i+1,j,\ell}^{xyt}+C_{34}{\bar{u}}_{i,j,\ell}^{xyt}+C_{35}{\bar{u}}_{i-1,j,\ell}^{xyt}+C_{36}{\bar{u}}_{i+1,j,\ell}^{xyt}
	\end{split}
\label{eq:059}
\end{gather}
\begin{gather}
	\begin{split}
		\frac{{\bar{u}}_{i,j,\ell}^{yt}-{\bar{u}}_{i-1,j,\ell}^{yt}}{2a_i} =  C_{51}{\bar{S}}_{3i,j,\ell}^{xyt}+C_{52}{\bar{S}}_{3i-1,j,\ell}^{xyt}+C_{53}{\bar{S}}_{3i+1,j,\ell}^{xyt}+C_{54}{\bar{u}}_{i,j,\ell}^{xyt}+C_{55}{\bar{u}}_{i-1,j,\ell}^{xyt}+C_{56}{\bar{u}}_{i+1,j,\ell}^{xyt}
	\end{split}
\label{eq:060}
\end{gather}
\begin{gather}
	\begin{split}
		\frac{{\bar{J}}_{xi,j,\ell}^{xt}-{\bar{J}}_{xi,j-1,\ell}^{xt}}{2b_j} =  C_{71}{\bar{S}}_{2i,j,\ell}^{xyt}+C_{72}{\bar{S}}_{2i,j-1,\ell}^{xyt}+C_{73}{\bar{S}}_{2i,j+1,\ell}^{xyt}+C_{74}{\bar{u}}_{i,j,\ell}^{xyt}+C_{75}{\bar{u}}_{i,j-1,\ell}^{xyt}+C_{76}{\bar{u}}_{i,j+1,\ell}^{xyt}
	\end{split}
\label{eq:061}
\end{gather}
\begin{gather}
	\begin{split}
		\frac{{\bar{u}}_{i,j,\ell}^{xt}-{\bar{u}}_{i,j-1,\ell}^{xt}}{2b_j} =  C_{91}{\bar{S}}_{2i,j,\ell}^{xyt}+C_{92}{\bar{S}}_{2i,j-1,\ell}^{xyt}+C_{93}{\bar{S}}_{2i,j+1,\ell}^{xyt}+C_{94}{\bar{u}}_{i,j,\ell}^{xyt}+C_{95}{\bar{u}}_{i,j-1,\ell}^{xyt}+C_{96}{\bar{u}}_{i,j+1,\ell}^{xyt}
	\end{split}
\label{eq:062}
\end{gather}
The explicit definitions of all the $C$-coefficients are provided in \ref{app:A}. Furthermore, by substituting the expressions obtained from \eqns\eqref{eq:059}-\eqref{eq:062} into the constraint relations (\eqns\eqref{eq:051}-\eqref{eq:053} and \ref{eq:057}), a system of four algebraic equations corresponding to the $x$-momentum equation is derived at each node, expressed entirely in terms of cell-centered variables. An analogous procedure can be followed to derive the corresponding four algebraic equations for the $y$-momentum equation. 

It is worth noting that all coefficients appearing in the final system of algebraic equations depend solely on the velocity components from the previous time level (i.e., ${\bar{u}}_{i,j,k}^p$ and ${\bar{v}}_{i,j,k}^p$). As a result, within the present RCCNIM framework, these coefficients are required to be computed only once per time step, thereby reducing the overall computational effort in comparison to the previously developed MCCNIM scheme for the Burgers’ equation. 
\subsection{Approximation of nonlinear convective velocities}
\label{sec:2.0.1}
The formulation of the convective velocities ${\bar{u}}_{i,j,\ell}^p$, ${\bar{v}}_{i,j,\ell}^p$, ${\bar{u}}_{i,j,\ell}^0$ and ${\bar{v}}_{i,j,\ell}^0$ within the RCCNIM scheme for the two-dimensional case is adopted in a manner consistent with the MCCNIM methodology \citep{ahmed2025efficient}, and can be written as:
\begin{gather}
	\begin{split}
		{\bar{u}}_{i,j,\ell}^p=\frac{{\bar{u}}_{i,j,\ell-1}^{yt}+{\bar{u}}_{i-1,j,\ell-1}^{yt}+{\bar{u}}_{i,j,\ell-1}^{xt}+{\bar{u}}_{i,j-1,\ell-1}^{xt}}{4}
	\end{split}
\label{eq:063}
\end{gather}
\begin{gather}
	\begin{split}
	{\bar{v}}_{i,j,\ell}^p=\frac{{\bar{v}}_{i,j,\ell-1}^{yt}+{\bar{v}}_{i-1,j,\ell-1}^{yt}+{\bar{v}}_{i,j,\ell-1}^{xt}+{\bar{v}}_{i,j-1,\ell-1}^{xt}}{4}
\end{split}
\label{eq:064}
\end{gather}
\begin{gather}
	\begin{split}
		{\bar{u}}_{i,j,\ell}^0=\frac{{\bar{u}}_{i,j,\ell}^{yt}+{\bar{u}}_{i-1,j,\ell}^{yt}+{\bar{u}}_{i,j,\ell}^{xt}+{\bar{u}}_{i,j-1,\ell}^{xt}}{4}
	\end{split}
\label{eq:065}
\end{gather}
\begin{gather}
	\begin{split}
	{\bar{v}}_{i,j,\ell}^0=\frac{{\bar{v}}_{i,j,\ell}^{yt}+{\bar{v}}_{i-1,j,\ell}^{yt}+{\bar{v}}_{i,j,\ell}^{xt}+{\bar{v}}_{i,j-1,\ell}^{xt}}{4}
\end{split}
\label{eq:066}
\end{gather}
It has been established in earlier studies that the above definitions of the convective velocities yield reasonably accurate numerical solutions. However, further improvement can be achieved by reformulating these definitions in accordance with the CCNIM framework. In CCNIM, the pseudo source terms are retained as additional discrete unknowns rather than being eliminated from the final algebraic system. Consequently, their influence must be consistently incorporated into the formulation of the node-averaged convective velocities, namely ${\bar{u}}_{i,j,\ell}^p$, ${\bar{v}}_{i,j,\ell}^p$, ${\bar{u}}_{i,j,\ell}^0$ and ${\bar{v}}_{i,j,\ell}^0$.

To account for this, the expressions for ${\bar{u}}^{yt}$, ${\bar{u}}^{xt}$, ${\bar{v}}^{yt}$, and ${\bar{v}}^{xt}$, as given in \eqns\eqref{eq:039}–\eqref{eq:044}, are substituted into \eqns\eqref{eq:063}–\eqref{eq:066}, resulting in modified definitions of the convective velocities. These revised formulations have been demonstrated in the literature \citep{ahmed2025efficient} to provide enhanced accuracy and computational efficiency. Accordingly, in the present study, only the modified forms of the convective velocities are employed, both for solving the selected benchmark problems and for comparative assessment. In particular, substituting the relations from \eqns\eqref{eq:039}–\eqref{eq:044} into \eqn\eqref{eq:063} yields the following modified expression:
\begin{gather}
	\begin{split}
		{\bar{u}}_{i,j,\ell}^p & =F_{71}{\bar{S}}_{3i-1,j,\ell-1}^{xyt}+F_{72}{\bar{S}}_{3i,j,\ell-1}^{xyt}+F_{73}{\bar{S}}_{3i+1,j,\ell-1}^{xyt}+F_{74}{\bar{u}}_{i-1,j,\ell-1}^{xyt}+F_{75}{\bar{u}}_{i,j,\ell-1}^{xyt}+F_{76}{\bar{u}}_{i+1,j,\ell-1}^{xyt}\\&+G_{71}{\bar{S}}_{2i,j-1,\ell-1}^{xyt}+G_{72}{\bar{S}}_{2i,j,\ell-1}^{xyt}+G_{73}{\bar{S}}_{2i,j+1,\ell-1}^{xyt}+G_{74}{\bar{u}}_{i,j-1,\ell-1}^{xyt}+G_{75}{\bar{u}}_{i,j,\ell-1}^{xyt}+G_{76}{\bar{u}}_{i,j+1,\ell-1}^{xyt} 
	\end{split}
\label{eq:067}
\end{gather}
Details of all $F$- and $G$ - type coefficients are summarized in \ref{app:A}. The remaining velocity components (${\bar{v}}_{i,j,\ell}^p$, ${\bar{u}}_{i,j,\ell}^0$, and ${\bar{v}}_{i,j,\ell}^0$) can be derived in a similar manner.
\section{Results and discussion}
\label{sec:3}
%
%To demonstrate the efficiency of the proposed RCCNIM method, it is compared with the previously developed MCCNIM scheme. Several test problems based on the 2D Burgers’ equation are considered, and the execution times of both methods are evaluated. The MCCNIM scheme has already been established as accurate and efficient compared to traditional NIM approaches. Therefore, instead of revisiting earlier comparisons, the focus here is on a direct comparison between the two schemes.  
In this section, the performance of the proposed RCCNIM scheme is evaluated through four benchmark problems based on the Burgers’ equation, including one one-dimensional case and three two-dimensional cases, all admitting exact analytical solutions. The one-dimensional example is primarily included for validation purposes, as it serves as a standard test case widely used in the literature for assessing the accuracy of nodal schemes. The subsequent two-dimensional test cases are designed to examine different aspects of the method. The first of these provides a comparative assessment between the MCCNIM and RCCNIM schemes in terms of accuracy and computational efficiency. The second example demonstrates the capability of the proposed approach to accurately resolve problems involving mixed boundary conditions. The final test case is specifically chosen to highlight the inherent upwinding property of the RCCNIM scheme. For all test problems, whether one- or two-dimensional, the nonlinear algebraic systems arising from the discretization are solved iteratively using a Picard-based approach.

We compute RMS error (in two-dimension) using,
\begin{gather}
	\begin{split}
		\text{RMS}=\sqrt{\frac{\sum_{i=1}^{n_x}\sum_{j=1}^{n_y}\left|{\bar{u}}_{i,j,\ell}^{xyt}-{\bar{u}}_{i,j,\ell}^{exact}\right|^2}{n_xn_y}} 
	\end{split}
\label{eq:068}
\end{gather}
Here, $n_x$ and $n_y$ represent the number of discretization points along the $x$- and $y$-directions, respectively, while $\ell$ denotes the temporal index associated with the final time step in the two-dimensional setting. All numerical computations were carried out using MATLAB R2024b as the development and execution platform. The simulations were performed on a system equipped with an Intel Core i9-12900K processor, featuring 16 cores (8 performance and 8 efficiency cores) and 24 threads, operating at base frequencies of 3.2 GHz (P-cores) and 2.4 GHz (E-cores), with a maximum turbo frequency reaching up to 5.2 GHz.

\subsection{Numerical tests}
\label{sec:3.1}
\subsubsection{Example 1: One-Dimensional Propagating Shock Problem}
\label{sec:3.1.2}
The exact solution of the one-dimensional Burgers’ equation defined in the domain $x \in [-2,2]$, corresponding to the propagating shock problem \citep{23Rizwan_uddin_1997}, is given by
\begin{gather}
	\begin{split}
		u\left(x,t\right)=\frac{1}{2}\left[1-\tanh\left(\frac{xRe}{4}-\frac{tRe}{8}\right)\right]  
	\end{split}
\label{eq:069}
\end{gather}
where the initial condition at $t=0$ and the boundary conditions at both the left and right ends are prescribed using the exact solution given in \eqn\eqref{eq:069}. 

This problem is selected to assess the accuracy of the proposed scheme, as it has been widely used in the literature, including the works of Decker \citep{decker1993block} and Rizwan-uddin \citep{23Rizwan_uddin_1997, 25Wescott_2001}, to examine the characteristics of traditional NIM. Furthermore, the same benchmark problem was employed in the previously developed MCCNIM formulation for the Burgers’ equation, where a detailed analysis was presented. Consequently, a comparison with error results reported in the literature is provided here to validate the performance of the current scheme.  
\begin{table}[!t]
	\begin{center}\renewcommand{\arraystretch}{1.7}
		\caption{Comparison of RMS errors between RCCNIM and existing methods—CN-4PU, CNIM, MNIM \citep{23Rizwan_uddin_1997}, and MCCNIM \citep{ahmed2025efficient}, with $\Delta t=0.1$ and a grid of 20 nodes.} \label{tab:1}
		\resizebox{0.6\textwidth}{!}{		
			\begin{tabular}{|l|c|c|c|c|}
				\hline
		Method & \multicolumn{4}{c|}{RMS errors ($\times 10^{-2}$)} \\ \cline{2-5}
		 & \multicolumn{2}{c|}{$Re = 50$} & \multicolumn{2}{c|}{$Re = 100$} \\ \cline{2-5}
		& $t = 1.0$ & $t = 3.0$ & $t = 1.0$ & $t = 3.0$ \\ \hline
		CN-4PU & 3.652 & 3.524 & 4.546 & 4.367 \\ \hline
		CNIM   & 3.091 & 3.068 & 6.558 & 6.376 \\ \hline
		MNIM ($\bar{u}^x$) & 1.338 & 1.382 & 1.172 & 1.185 \\ \hline
		MNIM ($\bar{u}^t$) & 2.721 & 2.775 & 2.666 & 2.678 \\ \hline
		MNIM ($\bar{u}^x$ and $\bar{u}^t$) & 2.127 & 2.175 & 2.041 & 2.053 \\ \hline
		MCCNIM ($\bar{u}^{xt}$) & 1.649 & 1.722 & 1.648 & 1.688 \\ \hline
		RCCNIM ($\bar{u}^{xt}$) & 1.481 & 1.565 & 1.492 & 1.531 \\ \hline
        \end{tabular}}
	\end{center}
\end{table}
%
%\begin{table}[!t]
%	\begin{center}\renewcommand{\arraystretch}{1.5}
%        \caption{RMS error comparison of the developed RCCNIM scheme with MNIM, M$^2$NIM \citep{25Wescott_2001}, and MCCNIM \citep{ahmed2025efficient} for $Re=100$ and $t=2.0$.}\label{tab:2}
%		\resizebox{0.6\textwidth}{!}{		
%			\begin{tabular}{|c|c|c|c|c|c|}
%				\hline
%		$\Delta x$ & $\Delta t$ & \multicolumn{4}{c|}{RMS errors ($\times 10^{-2}$)} \\ \cline{3-6}
%		&  & MNIM & M{$^2$}NIM & MCCNIM & RCCNIM \\ \hline		
%		0.2 & 0.02  & 4.88 & 4.75 & 2.61 & 2.59 \\ \cline{2-6}
%		& 0.01  & 4.61 & 4.56 & 2.67 & 2.66 \\ \cline{2-6}
%		& 0.005 & 4.48 & 4.45 & 2.69 & 2.69 \\ \hline
%		0.1 & 0.02  & 2.15 & 1.94 & 1.44 & 1.41 \\ \cline{2-6}
%		& 0.01  & 2.02 & 1.95 & 1.49 & 1.48 \\ \cline{2-6}
%		& 0.005 & 1.95 & 1.92 & 1.52 & 1.52 \\ \hline
%		0.05 & 0.02  & 0.759 & 0.601 & 0.615 & 0.589 \\ \cline{2-6}
%		& 0.01  & 0.757 & 0.679 & 0.644 & 0.630 \\ \cline{2-6}
%		& 0.005 & 0.752 & 0.730 & 0.653 & 0.646 \\ \hline
%		\end{tabular}}
%	\end{center}
%\end{table}
%
As observed from \tabs\ref{tab:1} and \ref{tab:2}, the RCCNIM scheme achieves accuracy that is comparable to, or better than, both the traditional NIM and the previously developed MCCNIM scheme. A comparison of computational time is not presented for this one-dimensional case, as the differences are negligible due to the relatively low overall execution time. Therefore, to provide a more meaningful assessment of computational efficiency, two-dimensional test cases are considered in the subsequent examples.
%
%\begin{table}[!t]
%	\begin{center}\renewcommand{\arraystretch}{1.5}
%		\caption{RMS error comparison of the developed RCCNIM scheme with MNIM, M$^2$NIM \citep{25Wescott_2001}, and MCCNIM \citep{ahmed2025efficient} for $Re=200$ and $t=2.0$.} \label{tab:3}
%		\resizebox{0.6\textwidth}{!}{		
%			\begin{tabular}{|c|c|c|c|c|c|}
%				\hline
%				%
%		$\Delta x$ & $\Delta t$ & \multicolumn{4}{c|}{RMS errors ($\times 10^{-2}$)} \\ \cline{3-6} 
%		 &  & MNIM & M{$^2$}NIM & MCCNIM & RCCNIM \\ \hline		
%		0.2 & 0.02 & 6.36 & 6.59 & 3.01 & 2.99 \\ \cline{2-6}
%		& 0.01 & 6.20 & 6.12 & 3.08 & 3.08 \\ \cline{2-6}
%		& 0.005 & 5.90 & 5.87 & 3.11 & 3.11 \\ \hline
%		0.1 & 0.02 & 3.77 & 3.45 & 1.79 & 1.77 \\ \cline{2-6}
%		& 0.01 & 3.43 & 3.33 & 1.89 & 1.88 \\ \cline{2-6}
%		& 0.005 & 3.24 & 3.20 & 1.93 & 1.93 \\ \hline
%		0.05 & 0.02 & 1.66 & 1.29 & 0.92 & 0.89 \\ \cline{2-6}
%		& 0.01 & 1.51 & 1.33 & 1.03 & 1.01 \\ \cline{2-6}
%		& 0.005 & 1.42 & 1.37 & 1.07 & 1.06 \\ \hline
%		\end{tabular}}
%	\end{center}
%\end{table}
%
\begin{table}[!t]
\begin{center}
\renewcommand{\arraystretch}{1.5}
\caption{RMS error comparison of the developed RCCNIM scheme with MNIM, M$^2$NIM \citep{25Wescott_2001}, and MCCNIM \citep{ahmed2025efficient} for $Re=100$ and $Re=200$ at $t=2.0$.}
\label{tab:2}
\resizebox{1\textwidth}{!}{
\begin{tabular}{|c|c|c c c c | c c c c|}
\hline
& & \multicolumn{8}{c|}{RMS errors ($\times 10^{-2}$)} \\ \cline{3-10}

$\Delta x$ & $\Delta t$ & \multicolumn{4}{c|}{$Re=100$} & \multicolumn{4}{c|}{$Re=200$} \\ \cline{3-10}

& & MNIM & M$^2$NIM & MCCNIM & RCCNIM & MNIM & M$^2$NIM & MCCNIM & RCCNIM \\ \hline

0.2 & 0.02  & 4.88 & 4.75 & 2.61 & 2.59 & 6.36 & 6.59 & 3.01 & 2.99 \\ %\cline{2-10}
& 0.01  & 4.61 & 4.56 & 2.67 & 2.66 & 6.20 & 6.12 & 3.08 & 3.08 \\ %\cline{2-10}
& 0.005 & 4.48 & 4.45 & 2.69 & 2.69 & 5.90 & 5.87 & 3.11 & 3.11 \\ \hline

0.1 & 0.02  & 2.15 & 1.94 & 1.44 & 1.41 & 3.77 & 3.45 & 1.79 & 1.77 \\ %\cline{2-10}
& 0.01  & 2.02 & 1.95 & 1.49 & 1.48 & 3.43 & 3.33 & 1.89 & 1.88 \\ %\cline{2-10}
& 0.005 & 1.95 & 1.92 & 1.52 & 1.52 & 3.24 & 3.20 & 1.93 & 1.93 \\ \hline

0.05 & 0.02  & 0.759 & 0.601 & 0.615 & 0.589 & 1.66 & 1.29 & 0.92 & 0.89 \\ %\cline{2-10}
& 0.01  & 0.757 & 0.679 & 0.644 & 0.630 & 1.51 & 1.33 & 1.03 & 1.01 \\ %\cline{2-10}
& 0.005 & 0.752 & 0.730 & 0.653 & 0.646 & 1.42 & 1.37 & 1.07 & 1.06 \\ \hline

\end{tabular}}
\end{center}
\end{table}
The numerical results presented in \tabs\ref{tab:1}-\ref{tab:2} demonstrate that the RCCNIM scheme achieves accuracy that is comparable to, and in several cases surpasses, both the traditional NIM and the previously developed MCCNIM scheme. This assessment is supported by a detailed comparison of root-mean-square (RMS) errors reported in \tab\ref{tab:2}, where RCCNIM is evaluated against MCCNIM, MNIM, and M$^2$NIM across a range of Reynolds numbers, spatial discretizations, and time-step sizes available in the literature \citep{25Wescott_2001}. The results consistently indicate that, even for the nonlinear Burgers’ equation, the RCCNIM formulation provides accurate and reliable solutions, remaining competitive with existing NIM-based approaches for fluid flow applications.

A comparison of computational time is not included for this one-dimensional problem, as the overall execution time is relatively small and the differences are not sufficiently significant to yield meaningful insights. Instead, computational efficiency is examined in subsequent two-dimensional test cases, which constitute the primary focus of the present study.
\subsubsection{Example 2: Two-dimensional Burgers' equations with Dirichlet boundary conditions: Propagating shock problem}
\label{sec:3.1.3}
This two-dimensional propagating shock problem \citep{25Wescott_2001} extends the one-dimensional case discussed earlier (Example 1). The exact solution exhibits a similar structure and can be expressed as follows:
\begin{gather}
	\begin{split}
		u\left(x,y,t\right)=\frac{1}{2}\left[1-\tanh\left(\frac{xRe}{4}+\frac{yRe}{4}-\frac{tRe}{4}\right)\right] 
	\end{split}
\label{eq:070}
\end{gather}
\begin{gather}
	\begin{split}
		v\left(x,y,t\right)=\frac{1}{2}\left[1-\tanh\left(\frac{xRe}{4}+\frac{yRe}{4}-\frac{tRe}{4}\right)\right] 
	\end{split}
\label{eq:071}
\end{gather}
The computational domain is defined over $[-2,2]\times[-2,2]$, and the corresponding initial and boundary conditions are obtained from the exact solution given in \eqn\eqref{eq:070} and \eqn\eqref{eq:071}.
\begin{figure}[!htbp]%
	\centering
	%\subfigure[~]{\includegraphics[width=0.5\linewidth]{Figures/fig02_a.jpg}\label{fig02a}}%\\%
	%\subfigure[~]{\includegraphics[width=0.5\linewidth]{Figures/fig02_b.jpg}\label{fig02b}}\\
	%\subfigure[~]{\includegraphics[width=0.5\linewidth]{Figures/fig02_c.jpg}\label{fig02c}}%\\%
	%\subfigure[~]{\includegraphics[width=0.5\linewidth]{Figures/fig02_d.jpg}\label{fig02d}}\\
	%\subfigure[~]{\includegraphics[width=0.5\linewidth]{Figures/fig02_e.jpg}\label{fig02e}}%\\%
	%\subfigure[~]{\includegraphics[width=0.5\linewidth]{Figures/fig02_f.jpg}\label{fig02f}}\\
	%
    \subfigure[~]{\includegraphics[width=0.5\linewidth]{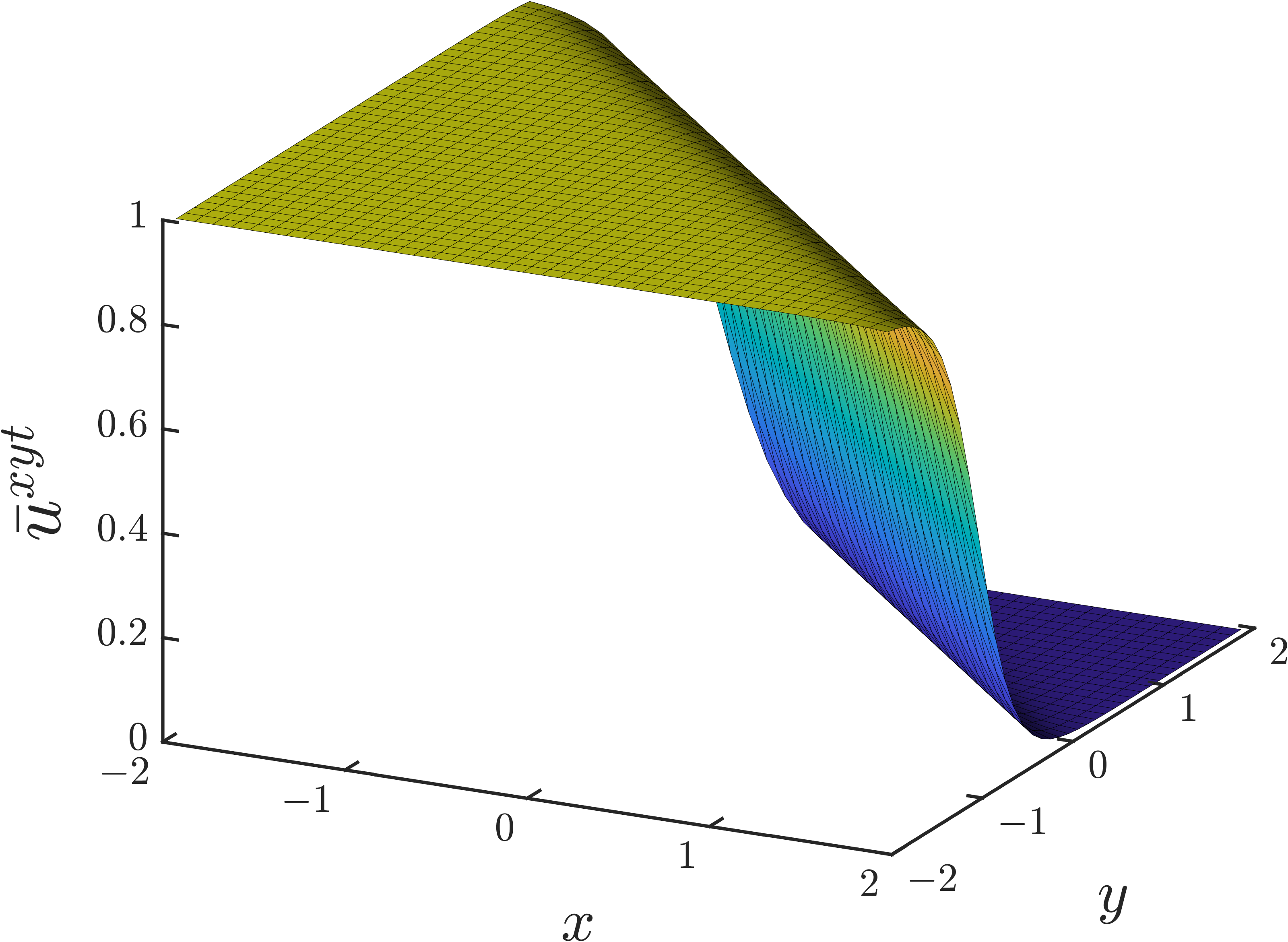}\label{fig02a}}%\\%
	\subfigure[~]{\includegraphics[width=0.5\linewidth]{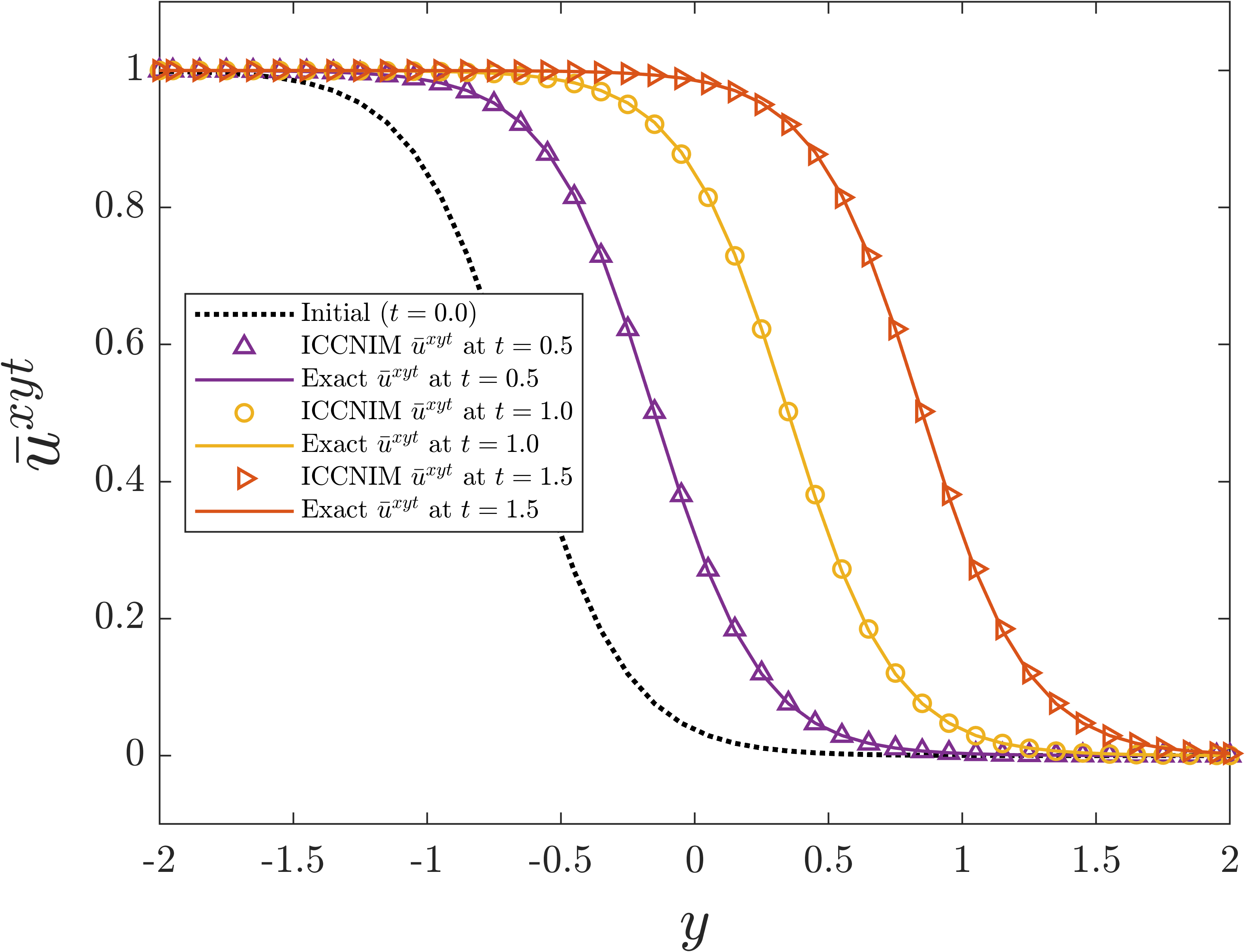}\label{fig02b}}\\
	\subfigure[~]{\includegraphics[width=0.5\linewidth]{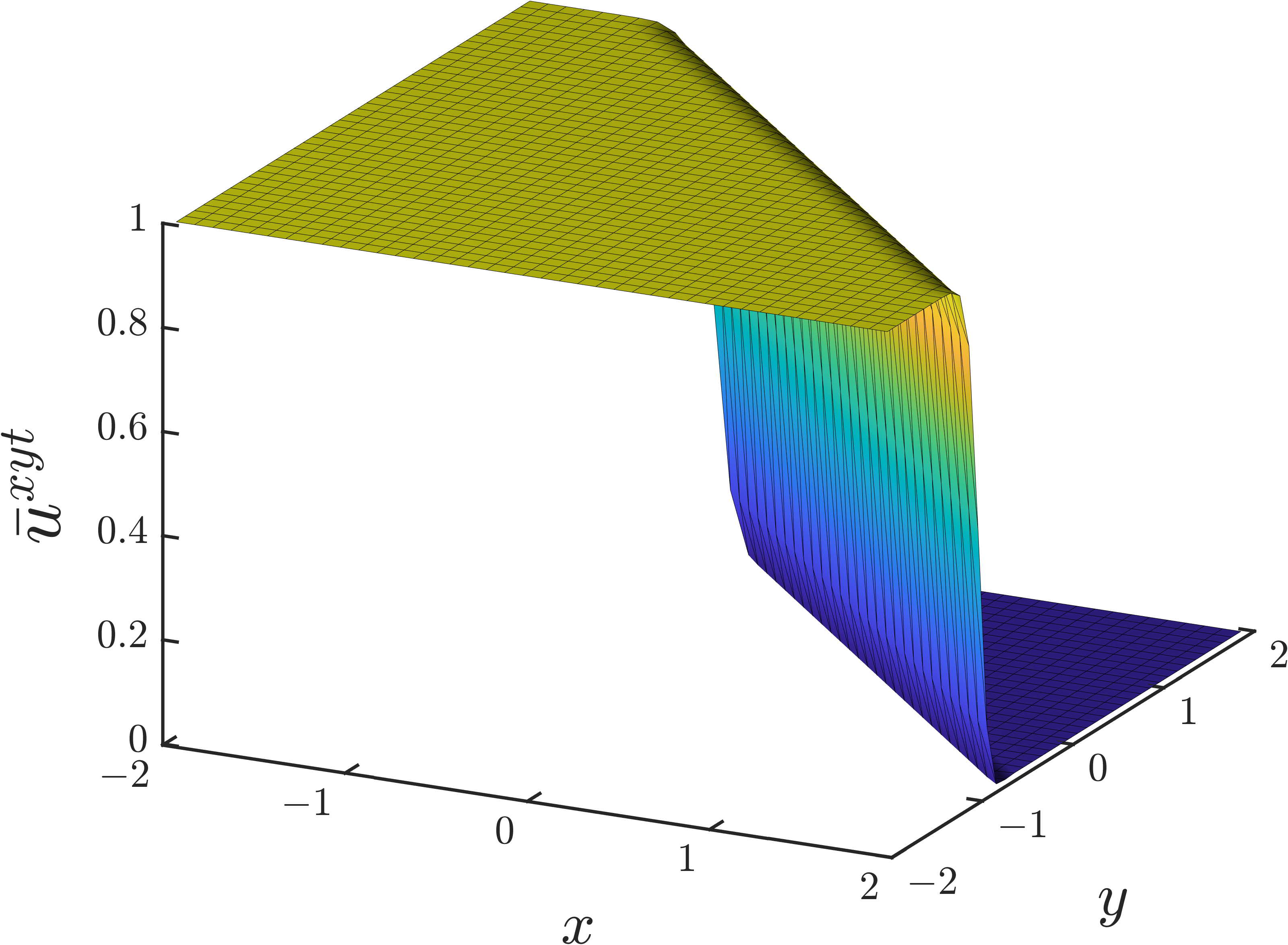}\label{fig02c}}%\\%
	\subfigure[~]{\includegraphics[width=0.5\linewidth]{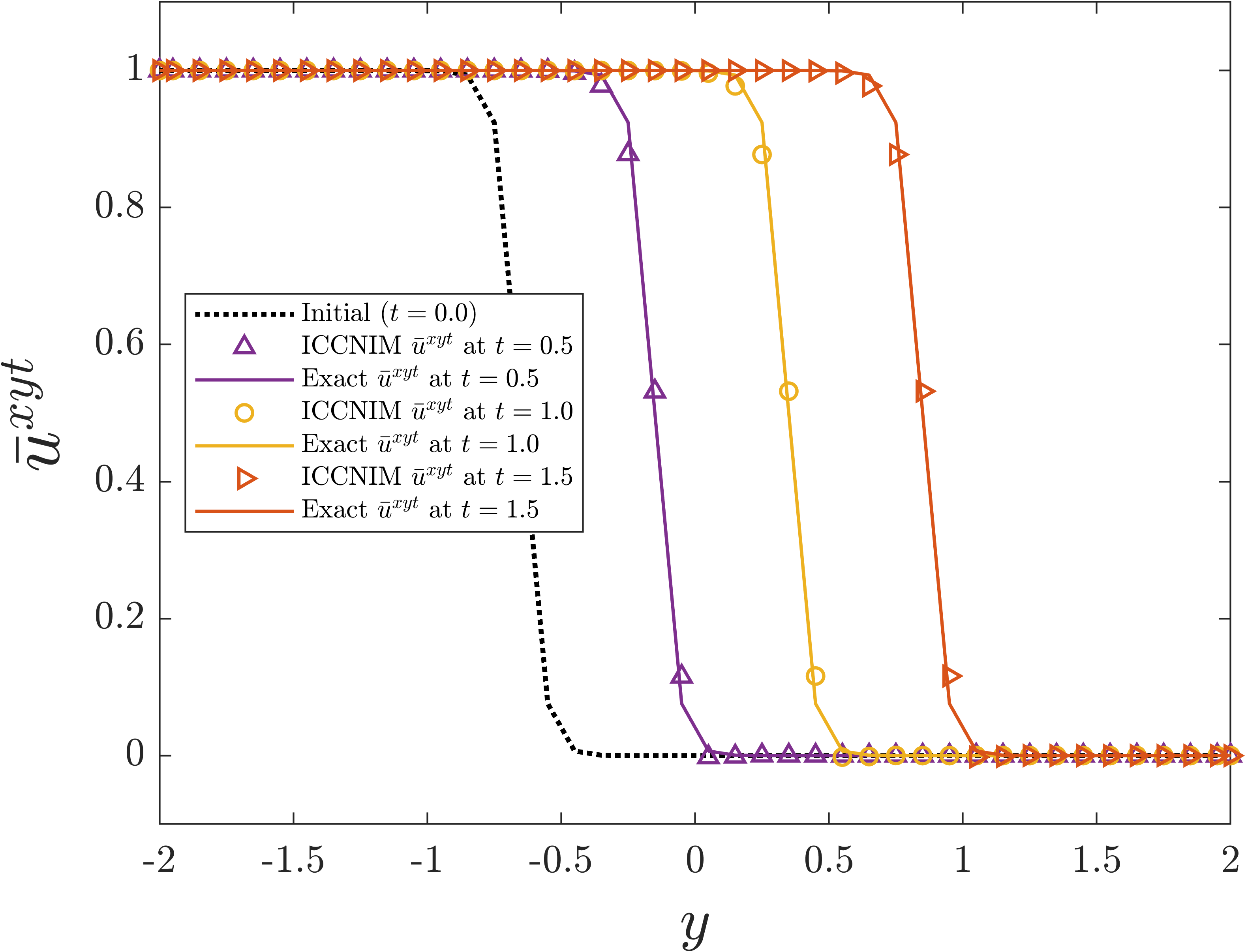}\label{fig02d}}\\
	\subfigure[~]{\includegraphics[width=0.5\linewidth]{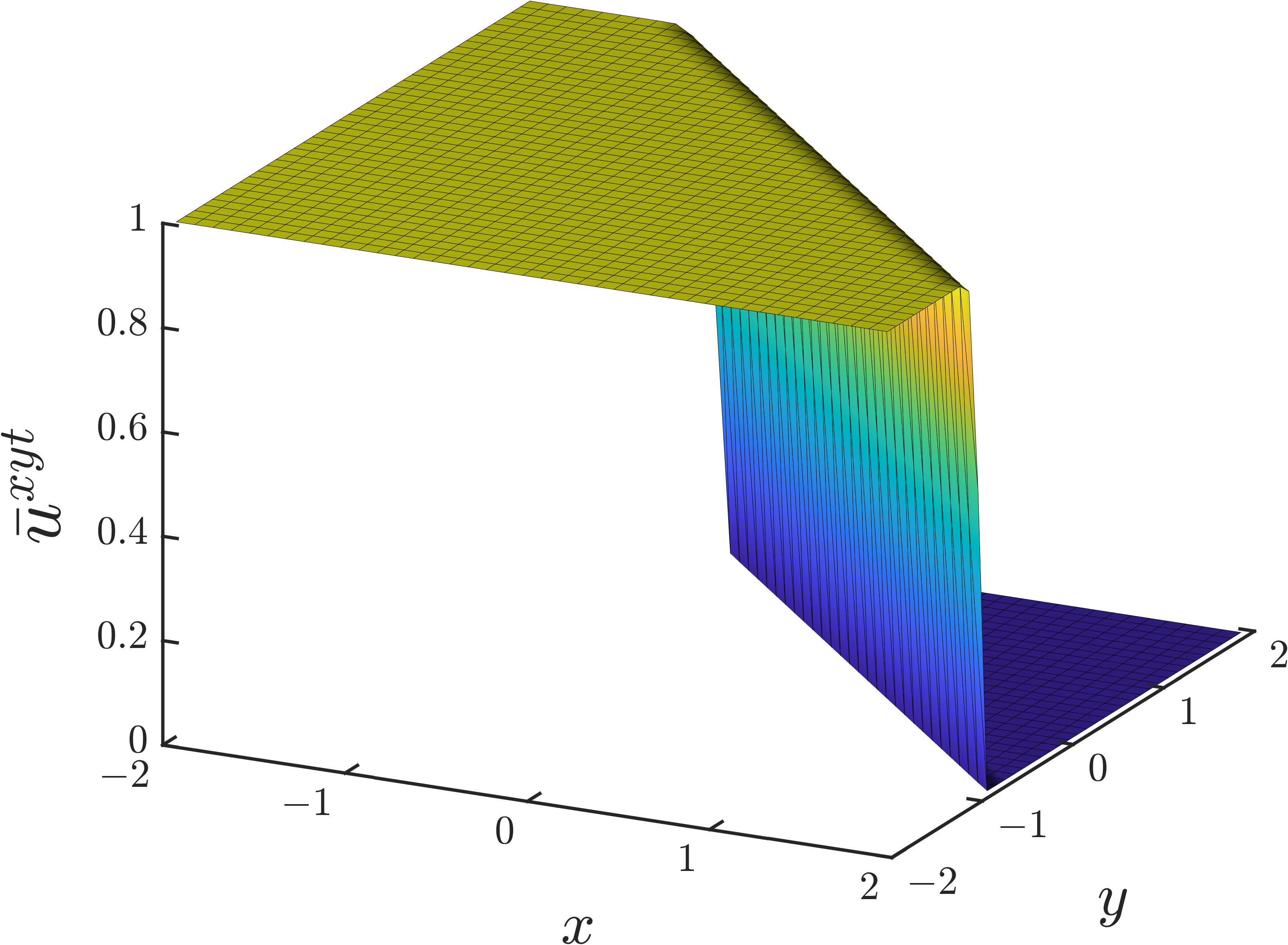}\label{fig02e}}%\\%
	\subfigure[~]{\includegraphics[width=0.5\linewidth]{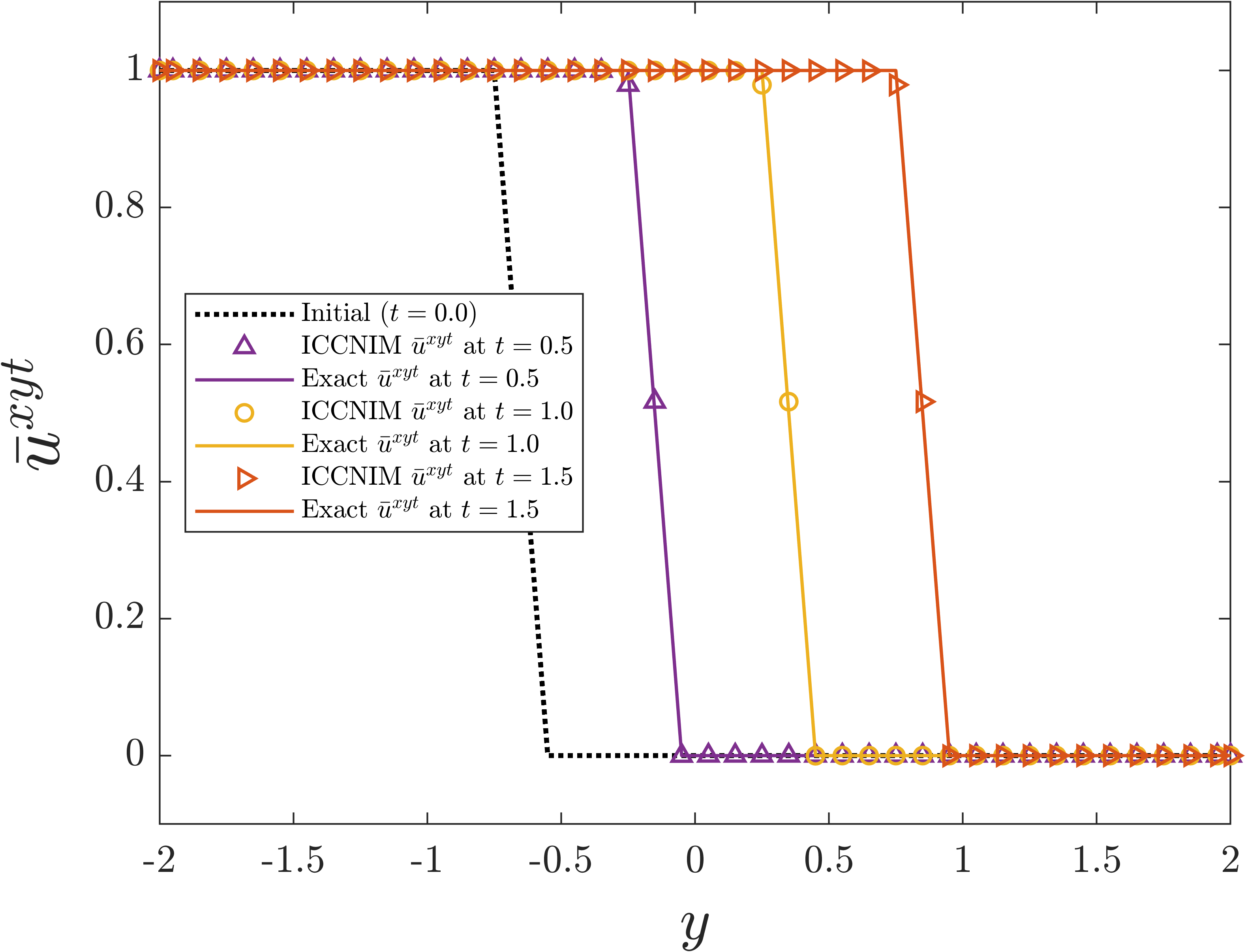}\label{fig02f}}\\
    \caption{Numerical solutions obtained using RCCNIM for different Reynolds numbers.
    (a, c, e) Surface plots depicting the overall solution structure at $t=2$ for
    $Re=10$, $Re=50$, and $Re=15000$, computed with a uniform grid size $\Delta x=0.05$
    and corresponding time step $\Delta t=0.001$. (b, d, f) Velocity profiles ($\bar{u}^{xyt}$) at $x=0.7$ at three different time levels, compared with the analytical solution for the respective Reynolds numbers.}%
	\label{fig:1}%
\end{figure}
\fig\ref{fig:1} presents the numerical results obtained using the RCCNIM scheme. The left column, comprising subfigures (a), (c), and (e), shows surface plots of the complete solution at time $t=2$ for three different Reynolds numbers, namely $Re=10$, $Re=50$, and $Re=15000$. These results are computed using a uniform grid size $\Delta x=0.05$ and corresponding time step $\Delta t=0.001$. The surface plots clearly demonstrate that, across all Reynolds numbers, the solution remains smooth and free from spurious oscillations, even as the wavefront becomes steeper at higher Reynolds numbers. It is noteworthy that, even at a high Reynolds number ($Re = 15000$), the wavefront retains a high degree of smoothness despite the use of a coarse spatial discretization ($\Delta x = 0.1$). This behavior confirms the stability and reliability of the proposed scheme for two-dimensional problems. For a detailed examination of the shock development, cross-sectional slices are depicted in the right column (subfigures (b), (d), and (f)) of \fig\ref{fig:1}, obtained by evaluating the solution along the line $x=0.7$ at multiple time levels. These plots facilitate a direct comparison between the numerical and analytical solutions. The results indicate a strong agreement between the two, even when relatively coarse grids are employed, despite the difficulty of resolving a sharply inclined shock front. In particular, \fig\ref{fig:1} highlights that the scheme maintains high accuracy even at a large Reynolds number ($Re=15000$), achieving close correspondence with the exact solution on a grid as coarse as $40 \times 40$.

\begin{table}[!t]
\begin{center}\renewcommand{\arraystretch}{2}
\caption{Comparison of RMS errors, iteration counts, and CPU runtimes for MCCNIM and RCCNIM at different Reynolds number ($Re$) and grid resolutions at $t=2$.}
 \label{tab:4}
\resizebox{1\textwidth}{!}{
\begin{tabular}{|c|c|c|c | c|c | c|c | c|c|}
\hline
$Re$ & $\Delta x$ & $\Delta t$ & \multicolumn{2}{c|}{RMS Errors ($\times 10^{-2}$)} & \multicolumn{2}{c|}{Iterations} & \multicolumn{2}{c|}{CPU Runtime (CR)} & Speedup  \\ \cline{4-9}
& &  & MCCNIM & RCCNIM & MCCNIM & RCCNIM & MCCNIM ($CR_M$) & RCCNIM ($CR_R$) & = $CR_M/CR_R$ \\ \hline

%\multicolumn{9}{|c|}{$Re=50$} \\ \hline
50 & 0.1    & 0.05   & 1.509 & 1.558 & 782   & 844   & 1.204 & 0.387 & 3.113 \\\cline{2-10} %\hline
& 0.05   & 0.025  & 0.804 & 0.814 & 1525  & 1605  & 9.544 & 3.318 & 2.877 \\ \cline{2-10}%\hline
& 0.025  & 0.0125 & 0.417 & 0.418 & 3360  & 3361  & 88.257 & 30.787 & 2.866 \\ \cline{2-10}%\hline
& 0.0125 & 0.005  & 0.169 & 0.169 & 12706 & 12355 & 1635.502 & 705.869 & 2.317 \\ \hline

%\multicolumn{9}{|c|}{$Re=500$} \\ \hline
500 &   0.1    & 0.05   & 2.060 & 2.179 & 999   & 1139  & 1.513 & 0.473 & 3.202 \\ \cline{2-10} %\hline
&   0.05   & 0.025  & 1.425  & 1.506 & 2118  & 2320  & 13.699 & 4.574 & 2.995 \\ \cline{2-10}%\hline
&   0.025  & 0.0125 & 0.937  & 0.998 & 4240  & 4559  & 123.716 & 42.297   & 2.925 \\ \cline{2-10}%\hline
& 0.0125 & 0.005  & 0.483 & 0.499 & 9922  & 10562 & 1407.787 & 630.877 & 2.232  \\ \hline

%\multicolumn{9}{|c|}{$Re=1000$} \\ \hline
1000    & 0.1    & 0.0125 & 0.451 & 0.463 & 3867  & 4106  & 5.965 & 1.845 & 3.234 \\ \cline{2-10}%\hline
&   0.05   & 0.01   & 0.518 & 0.554  & 5365  & 5563  & 36.082 & 11.211 & 3.218 \\ \cline{2-10}%\hline
&   0.025  & 0.005  & 0.346 & 0.358 & 11601 & 11599 & 358.441 & 122.940 & 2.916 \\ \cline{2-10}%\hline
& 0.0125 & 0.0025 & 0.259 & 0.263 & 23043 & 23525 & 3287.456 & 1407.136 & 2.336 \\ \hline
\end{tabular}}
\end{center}
\end{table}
To evaluate the computational efficiency of the two methods, simulations were performed for three Reynolds numbers, $Re=50, 500$, and $1000$ using different spatial discretizations ($\Delta x = \Delta y$) and time-step sizes ($\Delta t = 2\tau$). For all cases, the wave propagation was simulated up to a common final time of 2 s to ensure a consistent basis for comparison.  \tab\ref{tab:4} presents a comparison of the CPU time, iteration count, and RMS error for all test cases. To ensure a fair and unbiased assessment, both algorithms were implemented using an identical coding framework. In particular, the coefficient definitions and algebraic equation structures were kept the same for both methods. As discussed earlier, the adopted coefficient formulation can be readily applied to either the MCCNIM or RCCNIM approach, allowing the comparison to focus solely on the inherent differences between the numerical schemes. The results in \tabs\ref{tab:4} show that both methods require nearly the same number of iterations, which is expected since the overall implementation strategy and algebraic formulation are identical. Furthermore, RCCNIM achieves RMS errors comparable to those obtained with MCCNIM, demonstrating that both methods provide similar levels of accuracy. However, RCCNIM consistently requires significantly less computational time, reducing the CPU cost by approximately a factor of three across all cases considered. These results highlight the superior computational efficiency of RCCNIM while maintaining accuracy comparable to that of MCCNIM.

\begin{figure}[!htbp]%
	\centering
	%\subfigure[~]{\includegraphics[width=0.5\linewidth]{Figures/fig03_a.jpg}\label{fig03a}}%\\%
	%\subfigure[~] {\includegraphics[width=0.5\linewidth]{Figures/fig03_b.jpg}\label{fig03b}}\\
	%\subfigure[~]{\includegraphics[width=0.5\linewidth]{Figures/fig03_c.jpg}\label{fig03c}}%\\%
	%\subfigure[~] {\includegraphics[width=0.5\linewidth]{Figures/fig03_d.jpg}\label{fig03d}}\\
	\subfigure[~]{\includegraphics[width=0.5\linewidth]{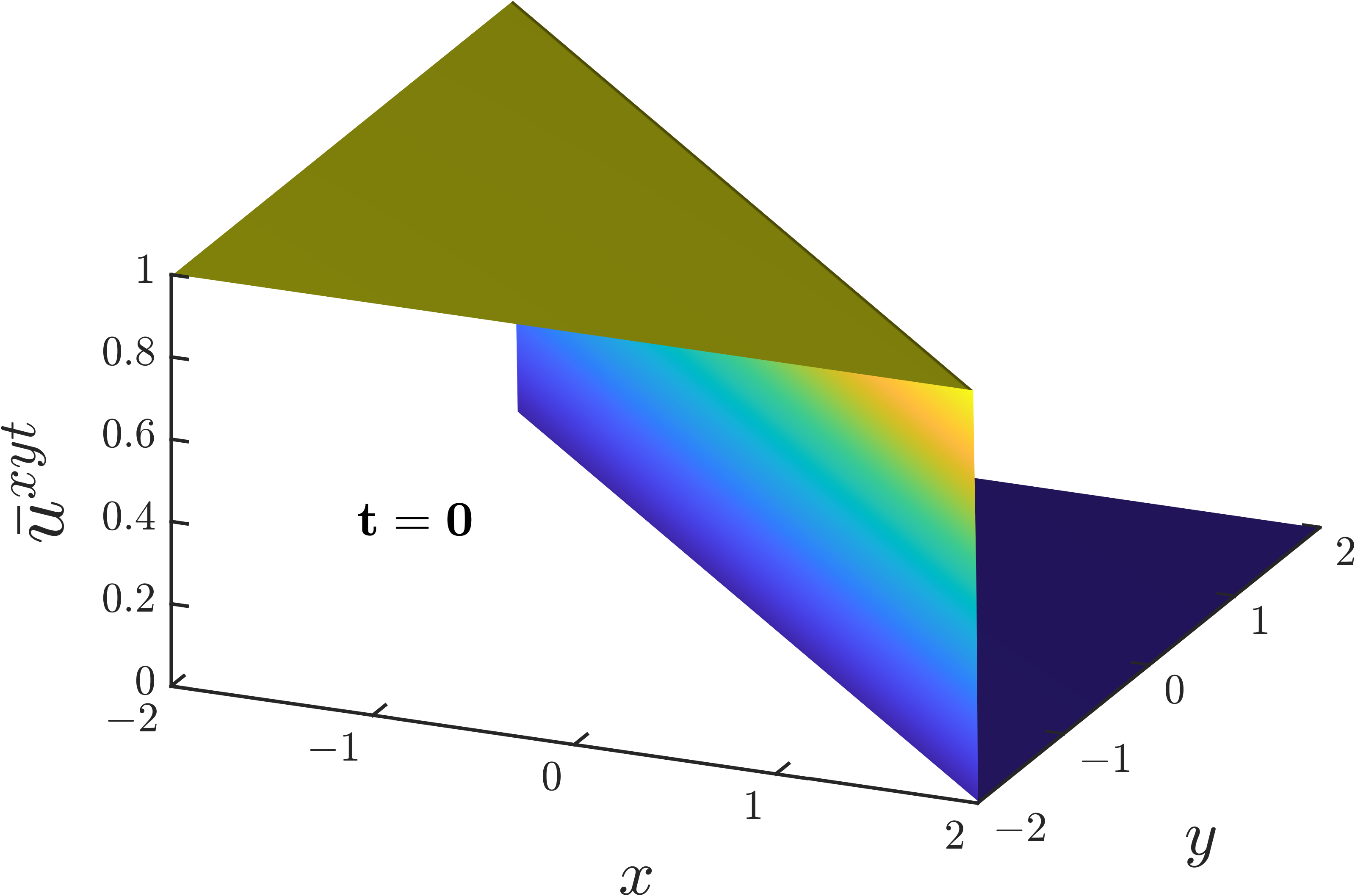}\label{fig03a}}%\\%
	\subfigure[~] {\includegraphics[width=0.5\linewidth]{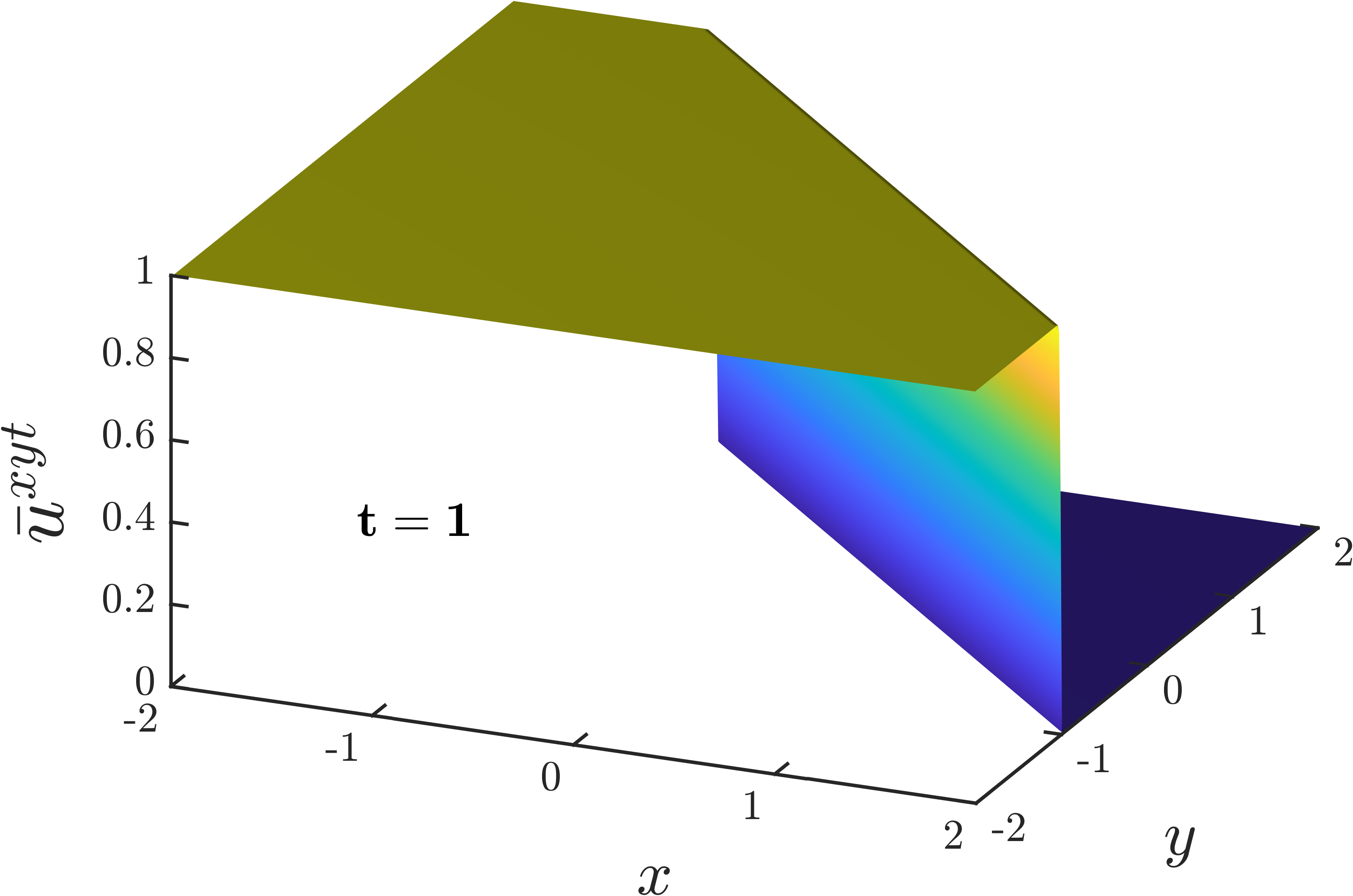}\label{fig03b}}\\
	\subfigure[~]{\includegraphics[width=0.5\linewidth]{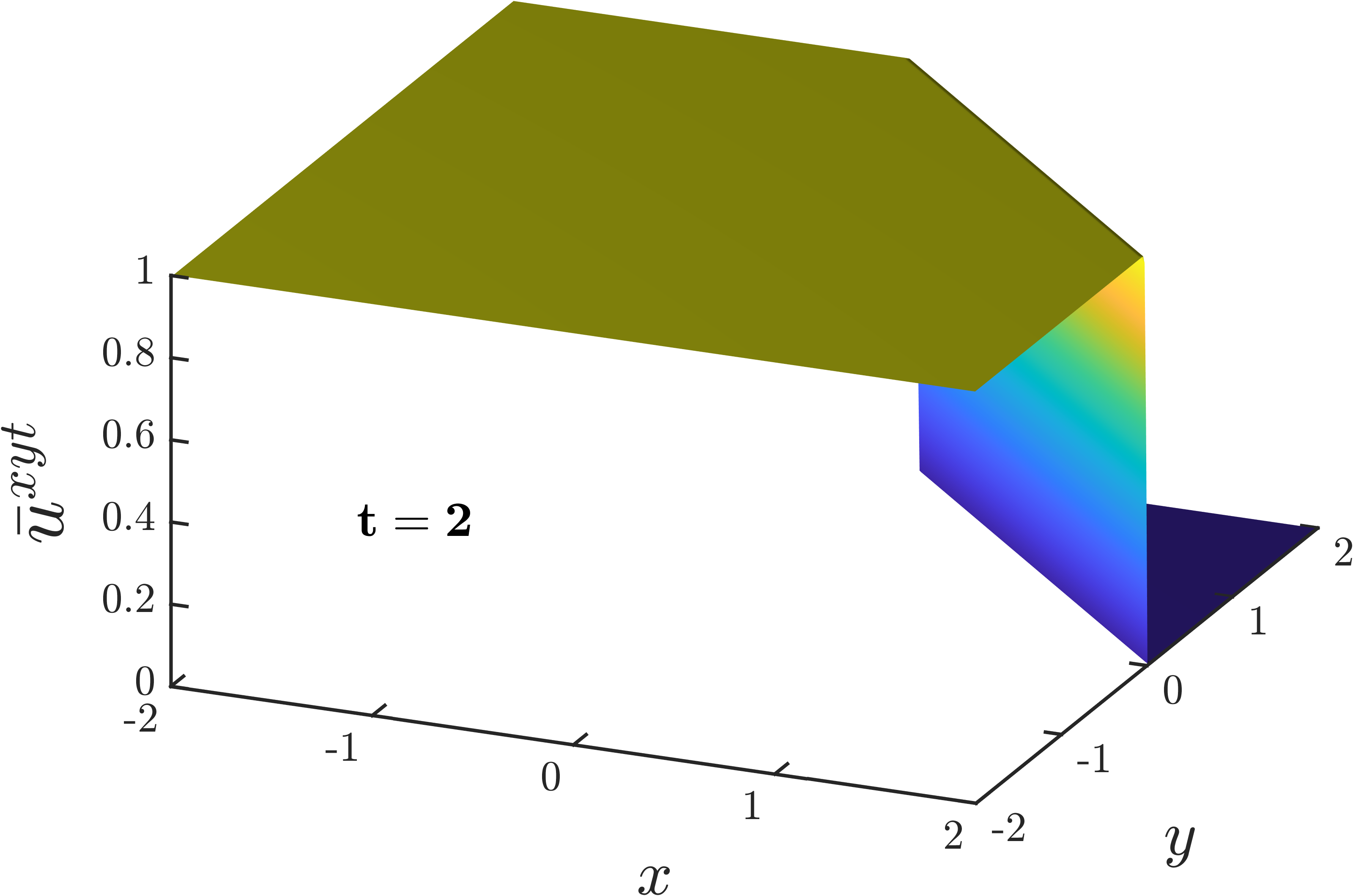}\label{fig03c}}%\\%
	\subfigure[~] {\includegraphics[width=0.5\linewidth]{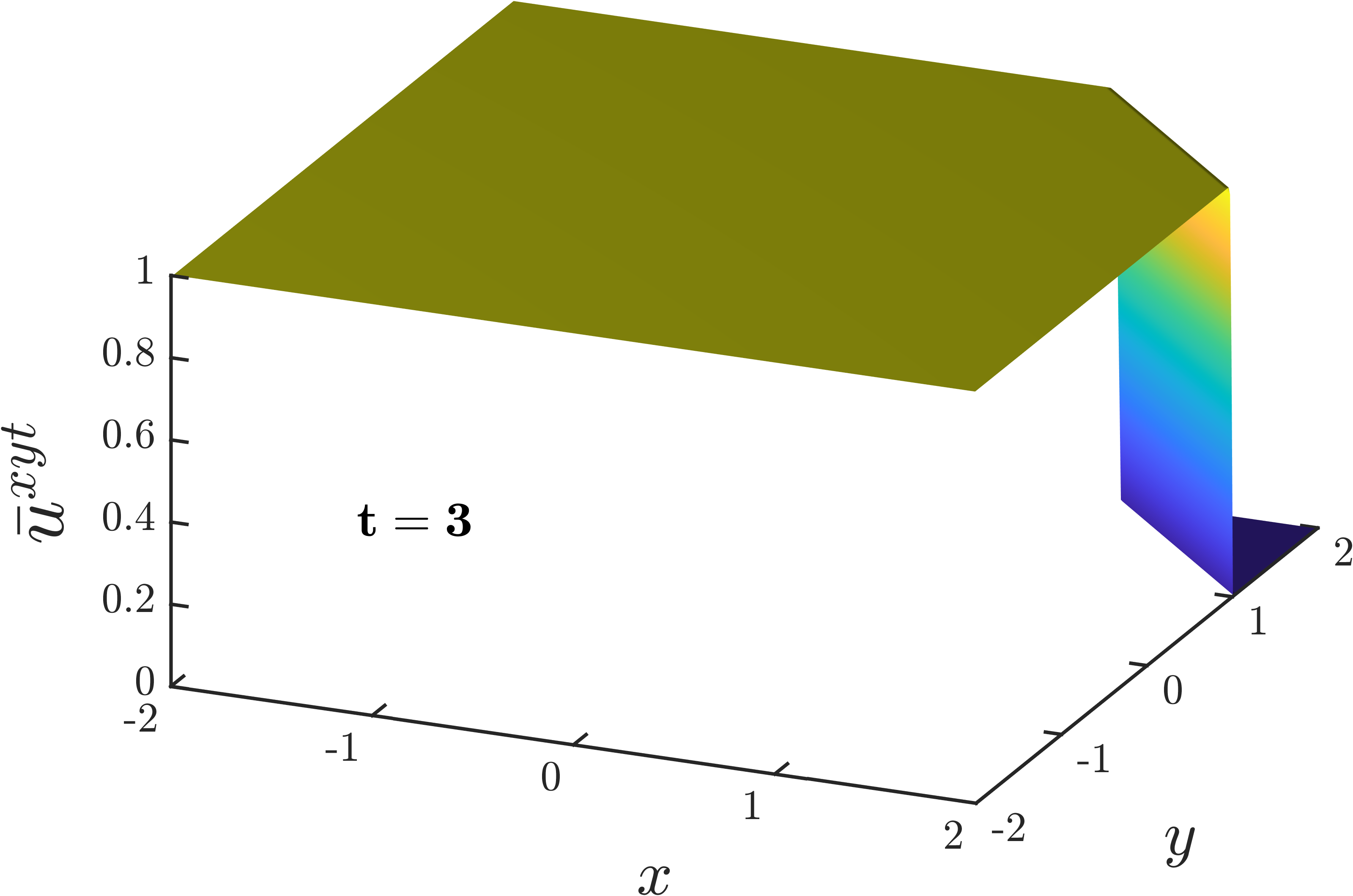}\label{fig03d}}\\
    \caption{Temporal evolution of the numerical solution for Example 2 obtained using the RCCNIM scheme at ($Re=10^5$), computed with ($\Delta t=0.0001$) and ($\Delta x=0.0167$).}%
	\label{fig:2}%
\end{figure}
To further assess the performance of the proposed RCCNIM scheme at high Reynolds numbers, the problem was solved for ($Re = 10^5$), and the corresponding results are shown in \fig\ref{fig:2}. The figure presents a sequence of three-dimensional surface plots illustrating the evolution of the propagating wave at different time instances ($t = 0-3$). These visualizations provide a comprehensive representation of the temporal development and propagation characteristics of the wave. It is evident from \fig\ref{fig:2} that the numerical solution remains smooth and free from spurious oscillations throughout the simulation. Remarkably, even at the high Reynolds number of ($Re = 10^5$), the wavefront is accurately captured on a relatively coarse computational grid ($\Delta x = 0.0167$)). The ability of the proposed scheme to maintain stability and resolution under such challenging conditions demonstrates its robustness and numerical accuracy.
\subsubsection{Example 3: Two-dimensional Burgers' equations with mixed boundary conditions}
\label{sec:3.1.3}
This example considers a two-dimensional Burgers’ problem with mixed boundary conditions, distinct from the previously studied propagating shock case. Unlike the earlier scenario, where shock propagation is the dominant feature, the present problem exhibits the formation of a discontinuity along the diagonal of the computational domain as time evolves \citep{zhang2010variational}. The governing equations are given by \eqns\eqref{eq:001} and \eqref{eq:002}, and the computational domain is defined over the unit square $[0,1]\times[0,1]$. The corresponding initial and boundary conditions are specified as follows:
\begin{gather*}
        u\left(x,y,0\right)=\sin(\pi x)\cos(\pi y), \qquad
        v\left(x,y,0\right)=\cos(\pi x)\sin(\pi y), \\
		u\left(0,y,t\right)=u\left(1,y,t\right)=0, \qquad v\left(x,0,t\right)=v\left(x,1,t\right)=0, \\
        \frac{\partial u}{\partial \mathrm{n}}\left(x,0,t\right)=\frac{\partial u}{\partial \mathrm{n}}\left(x,1,t\right)=0, \qquad
        \frac{\partial v}{\partial \mathrm{n}}\left(0,y,t\right)=\frac{\partial v}{\partial \mathrm{n}}\left(1,y,t\right)=0,
\end{gather*}
%
%\begin{flalign*}
%& u(x,y,0) = \sin(\pi x)\cos(\pi y), && \\
%& v(x,y,0) = \cos(\pi x)\sin(\pi y), && \\
%& u(0,y,t) = u(1,y,t) = 0, && \\
%& v(x,0,t) = v(x,1,t) = 0, && \\
%& \frac{\partial u}{\partial \mathrm{n}}(x,0,t) = \frac{\partial u}{\partial \mathrm{n}}(x,1,t) = 0, && \\
%& \frac{\partial v}{\partial \mathrm{n}}(0,y,t) = \frac{\partial v}{\partial \mathrm{n}}(1,y,t) = 0 && 
%\end{flalign*}
%
Achieving stable numerical solutions for this problem in highly convective regimes has long been a challenging task for numerical methods \citep{zhang2010variational,gao2017analytical,chai2020appropriate}. As the Reynolds number increases, the solution develops increasingly sharp gradients, making the computation more susceptible to numerical instabilities, particularly in the vicinity of these steep fronts. In much of the existing literature \citep{zhang2009element,zhang2010variational,gao2017analytical}, oscillation-free solutions are typically reported for relatively low Reynolds numbers ($Re = 100$). However, for higher Reynolds numbers ($Re > 100$), especially at ($Re \geq 10^3$), numerical solutions often exhibit spurious oscillations, or alternatively, very fine computational grids are required to adequately resolve these oscillations.
\begin{figure}[!htbp]%
	\centering
	%
    %\subfigure[~]{\includegraphics[width=0.5\linewidth]{Figures/fig04_a.jpg}\label{fig04a}}%\\%
	%\subfigure[~]{\includegraphics[width=0.5\linewidth]{Figures/fig04_b.jpg}\label{fig04b}}\\
	%\subfigure[~]{\includegraphics[width=0.5\linewidth]{Figures/fig04_c.jpg}\label{fig04c}}%\\%
	%\subfigure[~]{\includegraphics[width=0.5\linewidth]{Figures/fig04_d.jpg}\label{fig04d}}\\
	%\subfigure[~]{\includegraphics[width=0.5\linewidth]{Figures/fig04_e.jpg}\label{fig04e}}%\\%
	%\subfigure[~]{\includegraphics[width=0.5\linewidth]{Figures/fig04_f.jpg}\label{fig04f}}\\
	%
    \subfigure[~]{\includegraphics[width=0.5\linewidth]{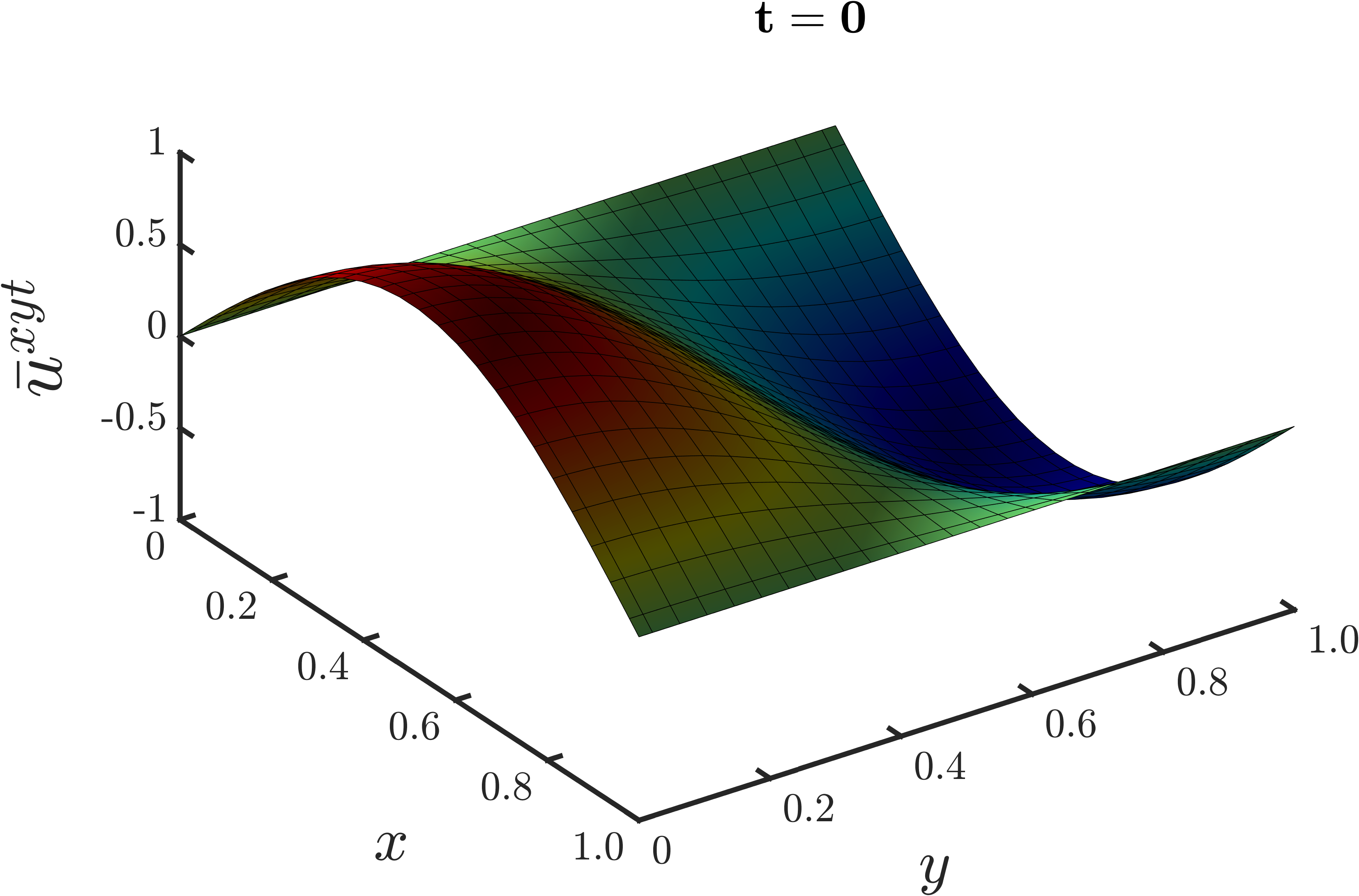}\label{fig04a}}%\\%
	\subfigure[~]{\includegraphics[width=0.5\linewidth]{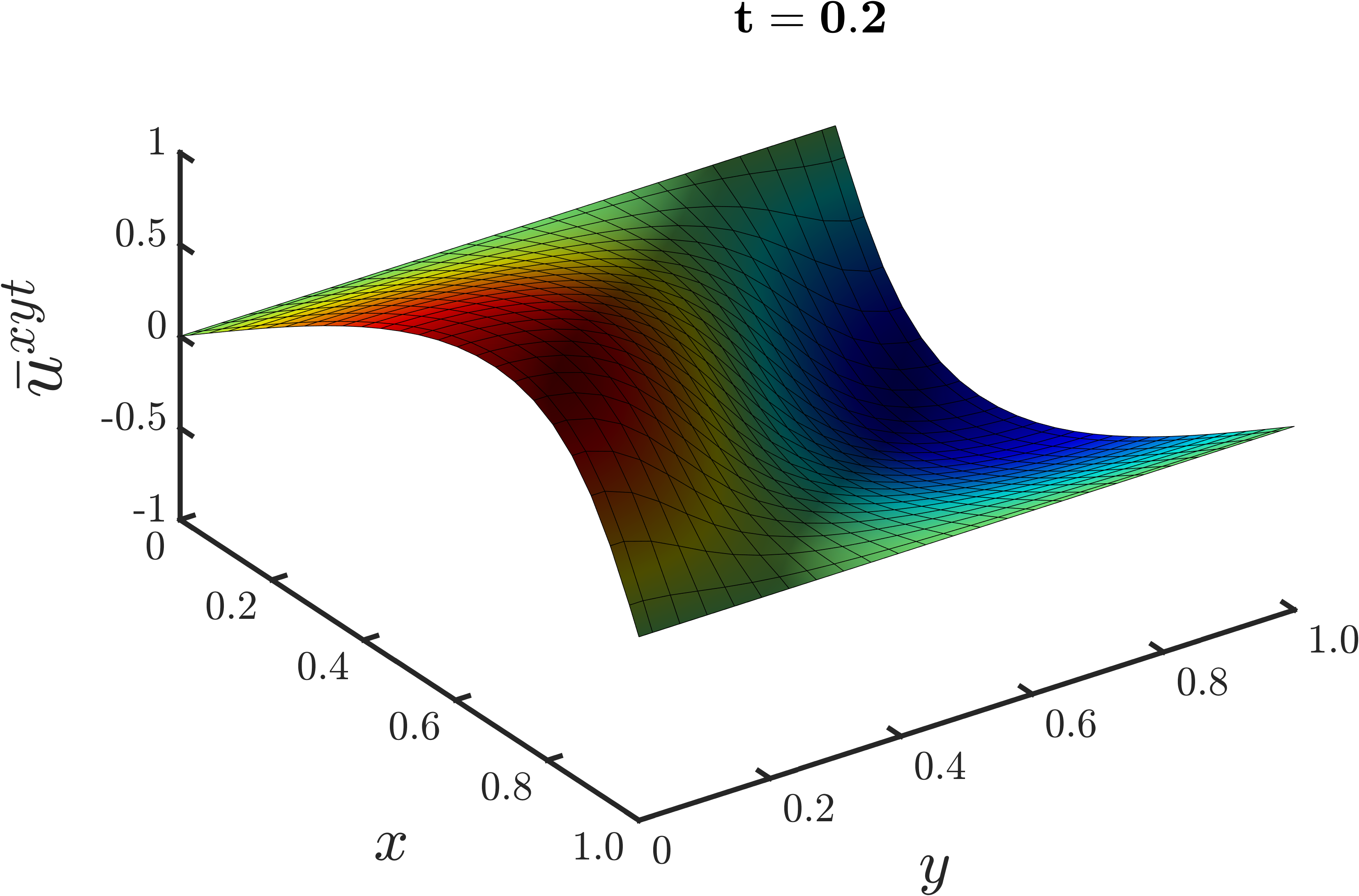}\label{fig04b}}\\
	\subfigure[~]{\includegraphics[width=0.5\linewidth]{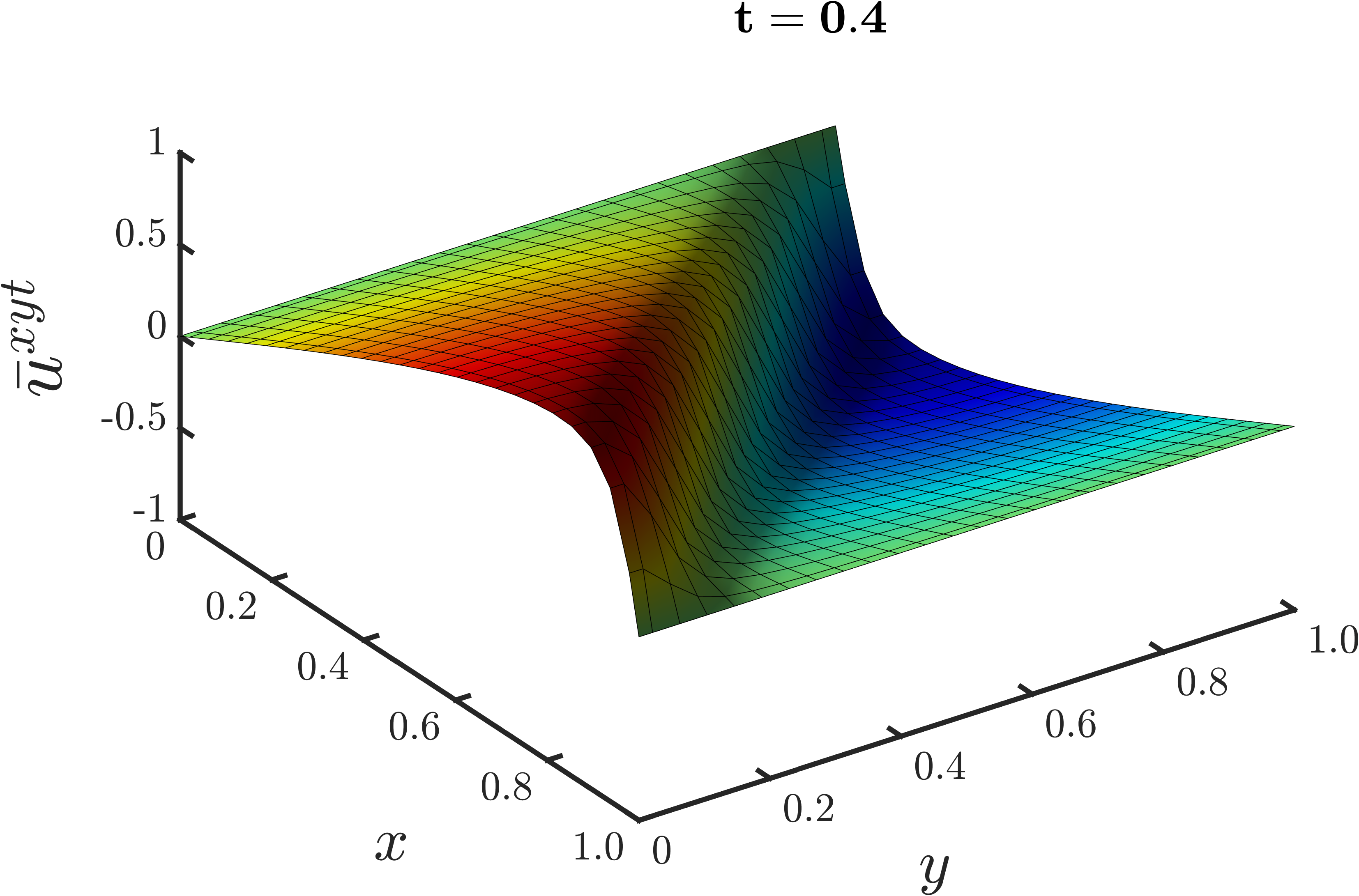}\label{fig04c}}%\\%
	\subfigure[~]{\includegraphics[width=0.5\linewidth]{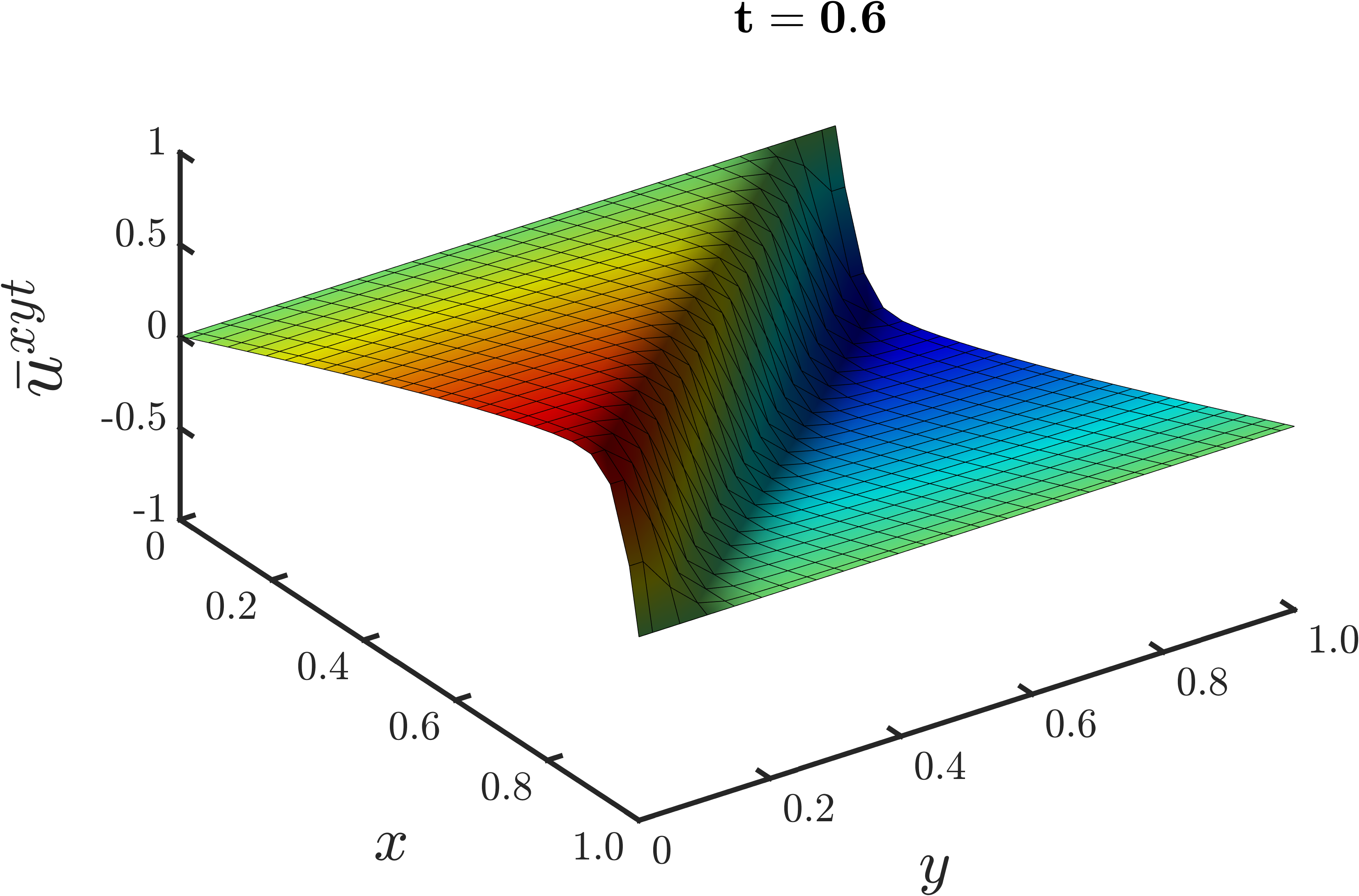}\label{fig04d}}\\
	\subfigure[~]{\includegraphics[width=0.5\linewidth]{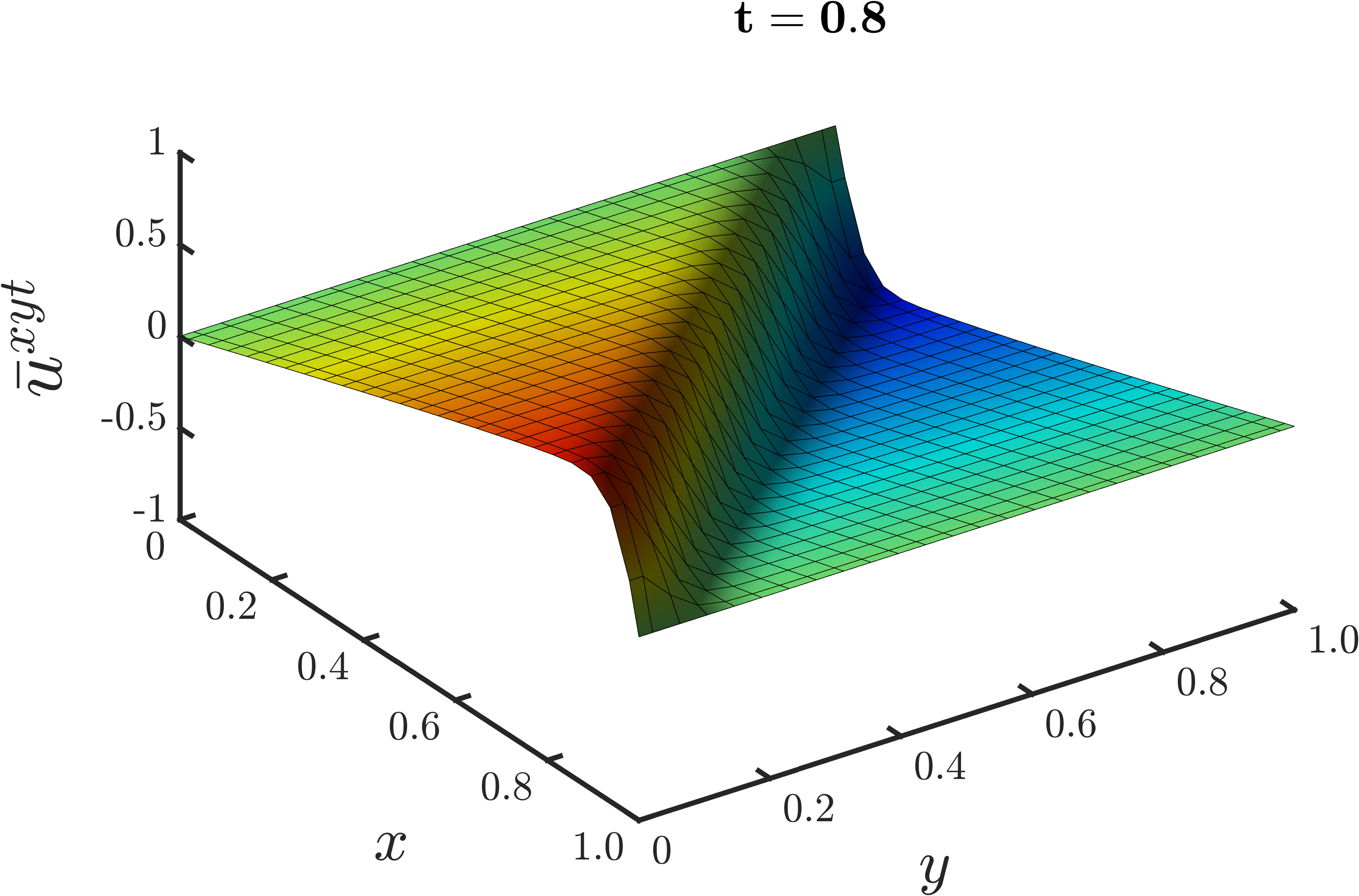}\label{fig04e}}%\\%
	\subfigure[~]{\includegraphics[width=0.5\linewidth]{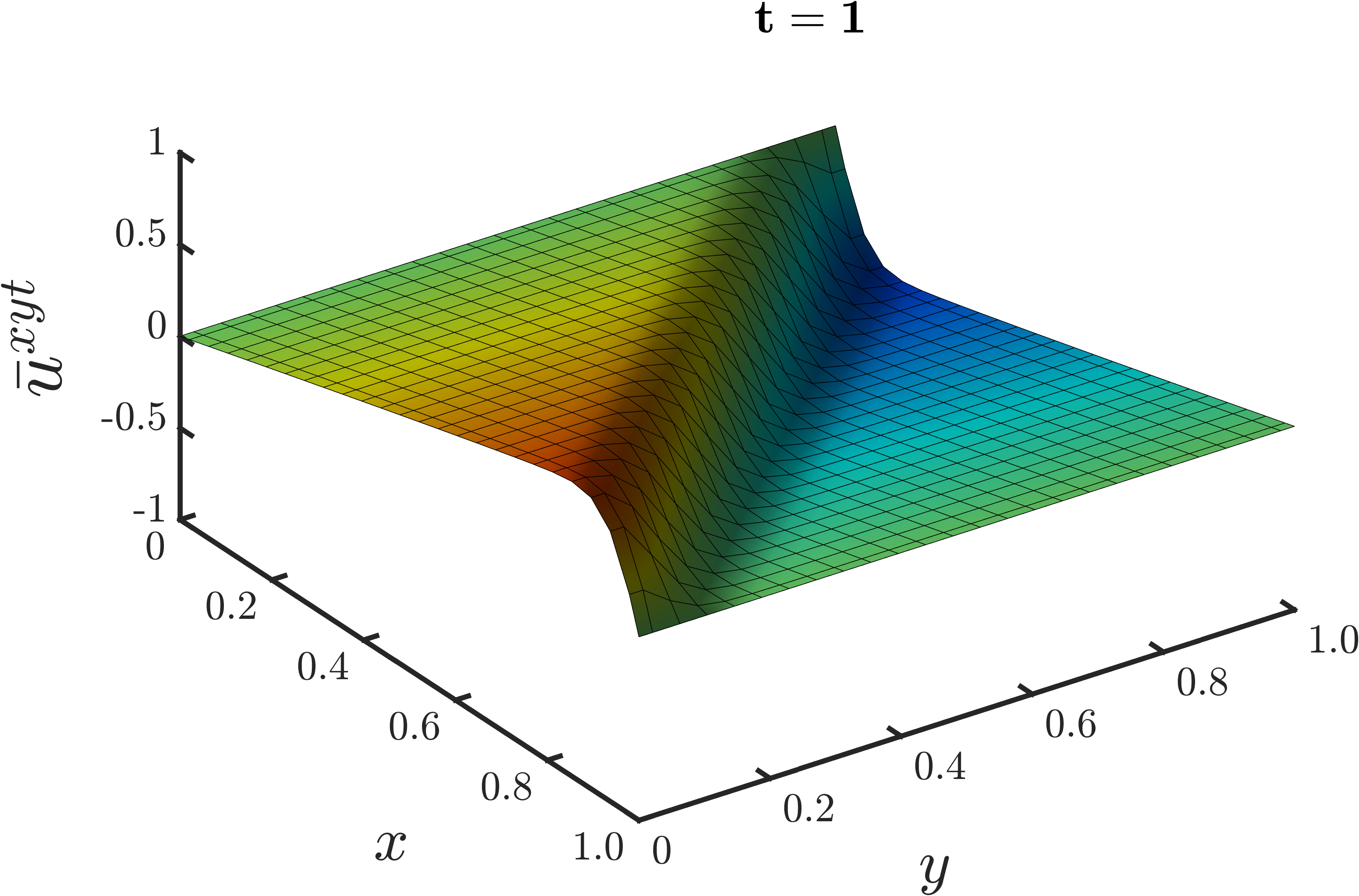}\label{fig04f}}\\
    \caption{Profiles of the velocity component $\bar{u}^{xyt}$ at different time instances for $Re=100$ on a $24 \times 24$ computational mesh.}
	\label{fig:3}%
\end{figure}
\begin{figure}[!htbp]%
	\centering
    %
	%\subfigure[~]{\includegraphics[width=0.5\linewidth]{Figures/fig05_a.jpg}\label{fig05a}}%\\%
	%\subfigure[~] {\includegraphics[width=0.5\linewidth]{Figures/fig05_b.jpg}\label{fig05b}}\\
	%\subfigure[~]{\includegraphics[width=0.5\linewidth]{Figures/fig05_c.jpg}\label{fig05c}}%\\%
	%\subfigure[~] {\includegraphics[width=0.5\linewidth]{Figures/fig05_d.jpg}\label{fig05d}}\\
	%\subfigure[~]{\includegraphics[width=0.5\linewidth]{Figures/fig05_e.jpg}\label{fig05e}}%\\%
	%\subfigure[~] {\includegraphics[width=0.5\linewidth]{Figures/fig05_f.jpg}\label{fig05f}}\\
	%
    \subfigure[~]{\includegraphics[width=0.5\linewidth]{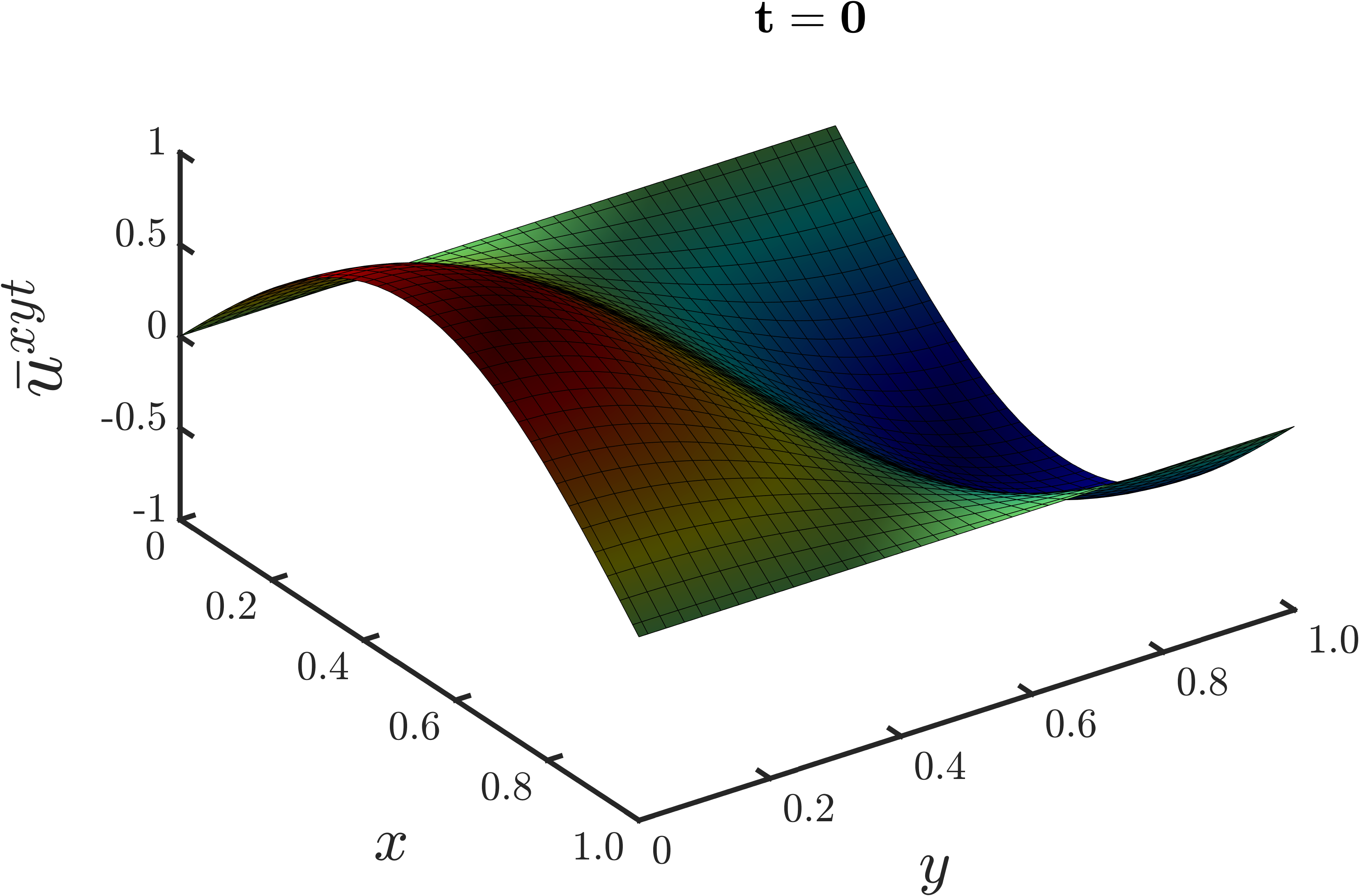}\label{fig05a}}%\\%
	\subfigure[~]{\includegraphics[width=0.5\linewidth]{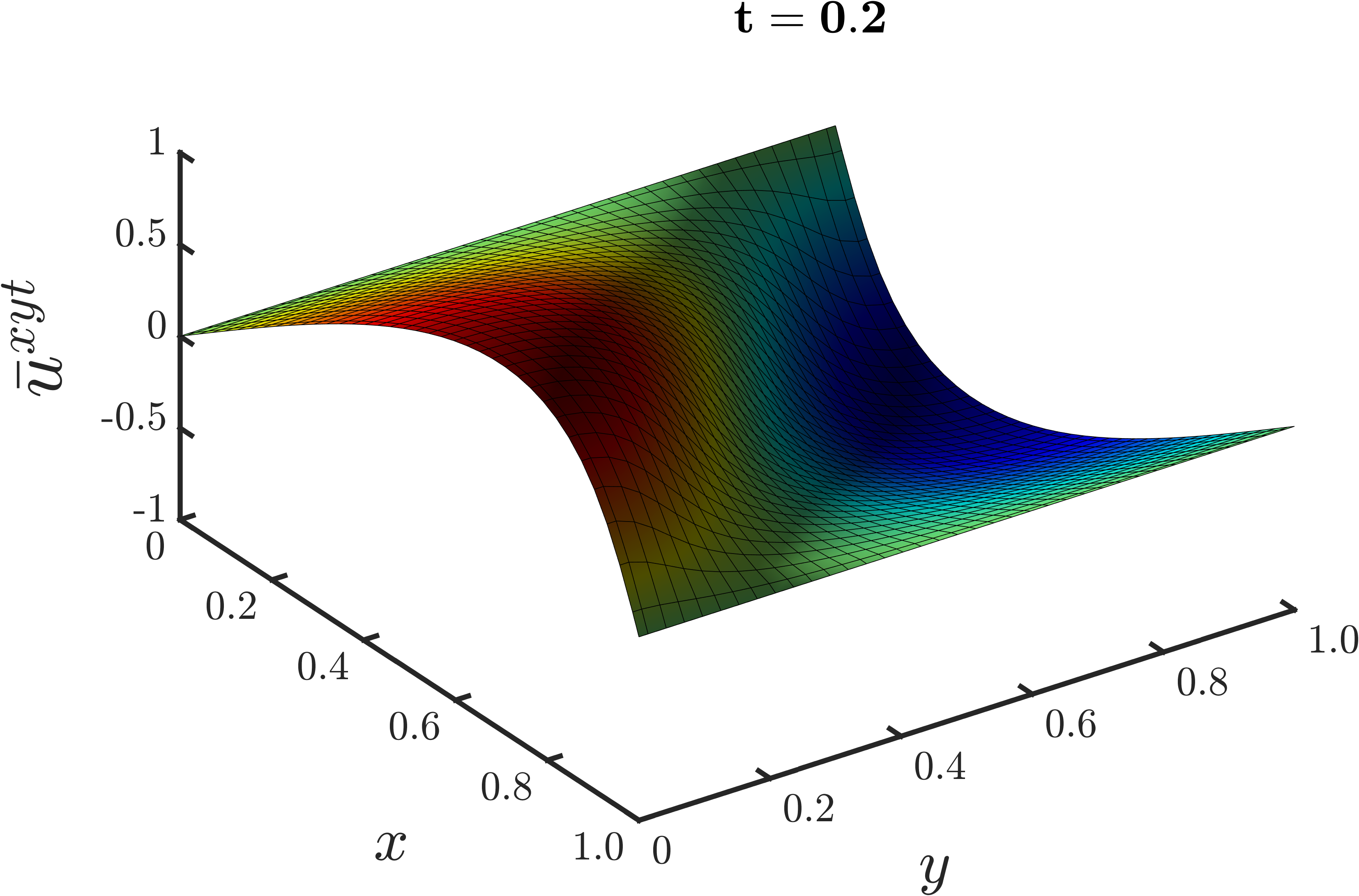}\label{fig05b}}\\
	\subfigure[~]{\includegraphics[width=0.5\linewidth]{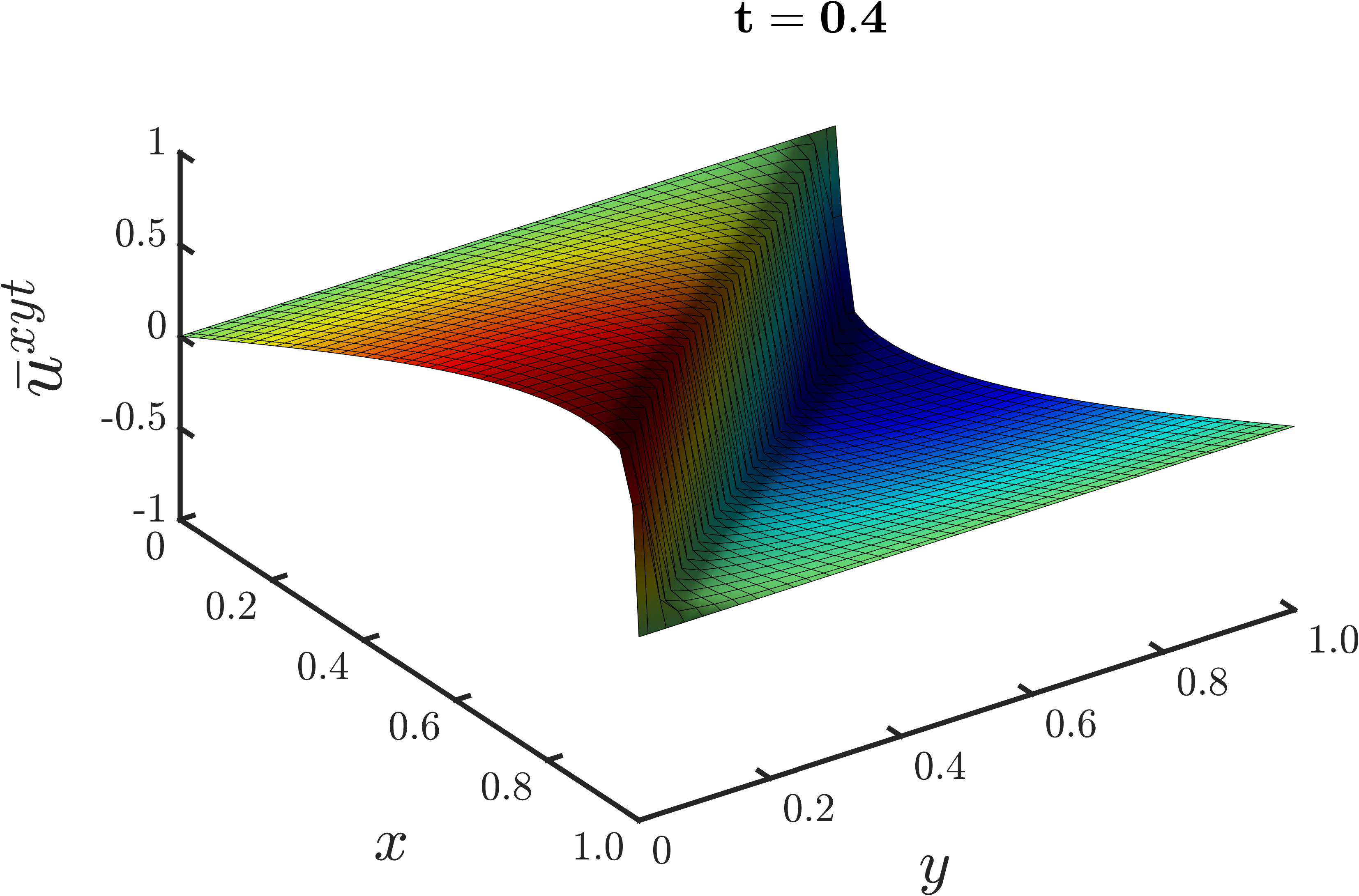}\label{fig05c}}%\\%
	\subfigure[~]{\includegraphics[width=0.5\linewidth]{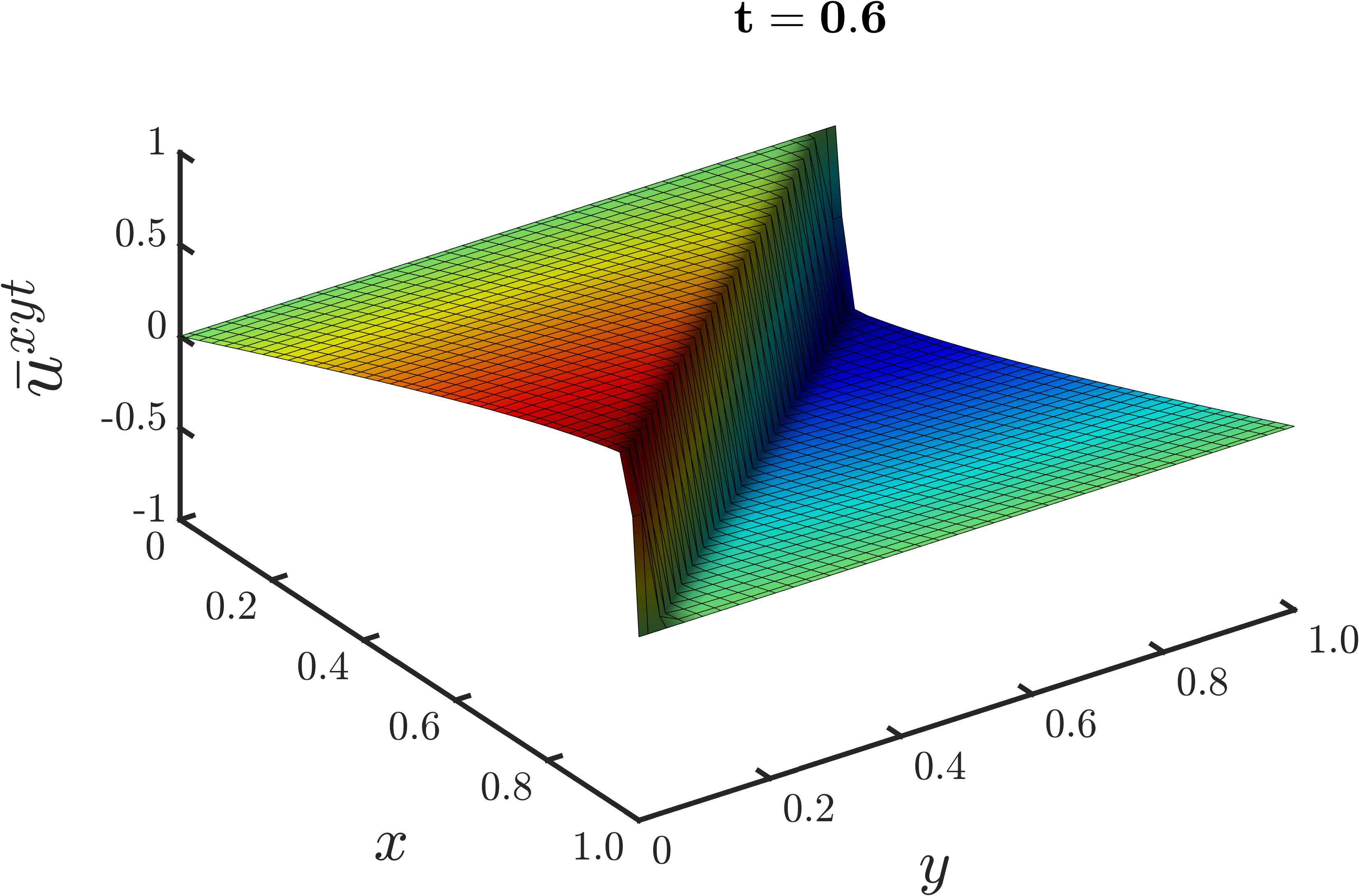}\label{fig05d}}\\
	\subfigure[~]{\includegraphics[width=0.5\linewidth]{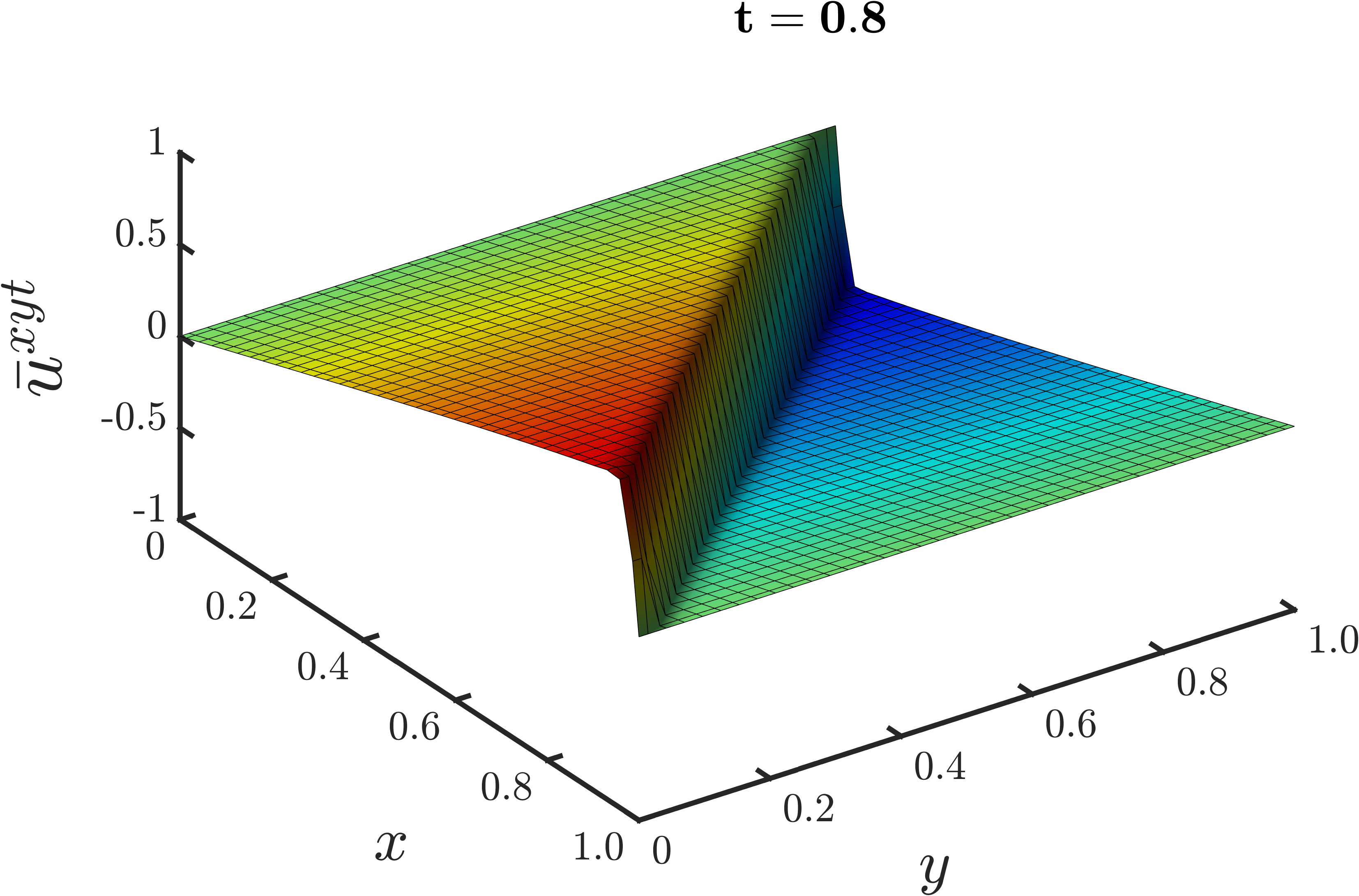}\label{fig05e}}%\\%
	\subfigure[~]{\includegraphics[width=0.5\linewidth]{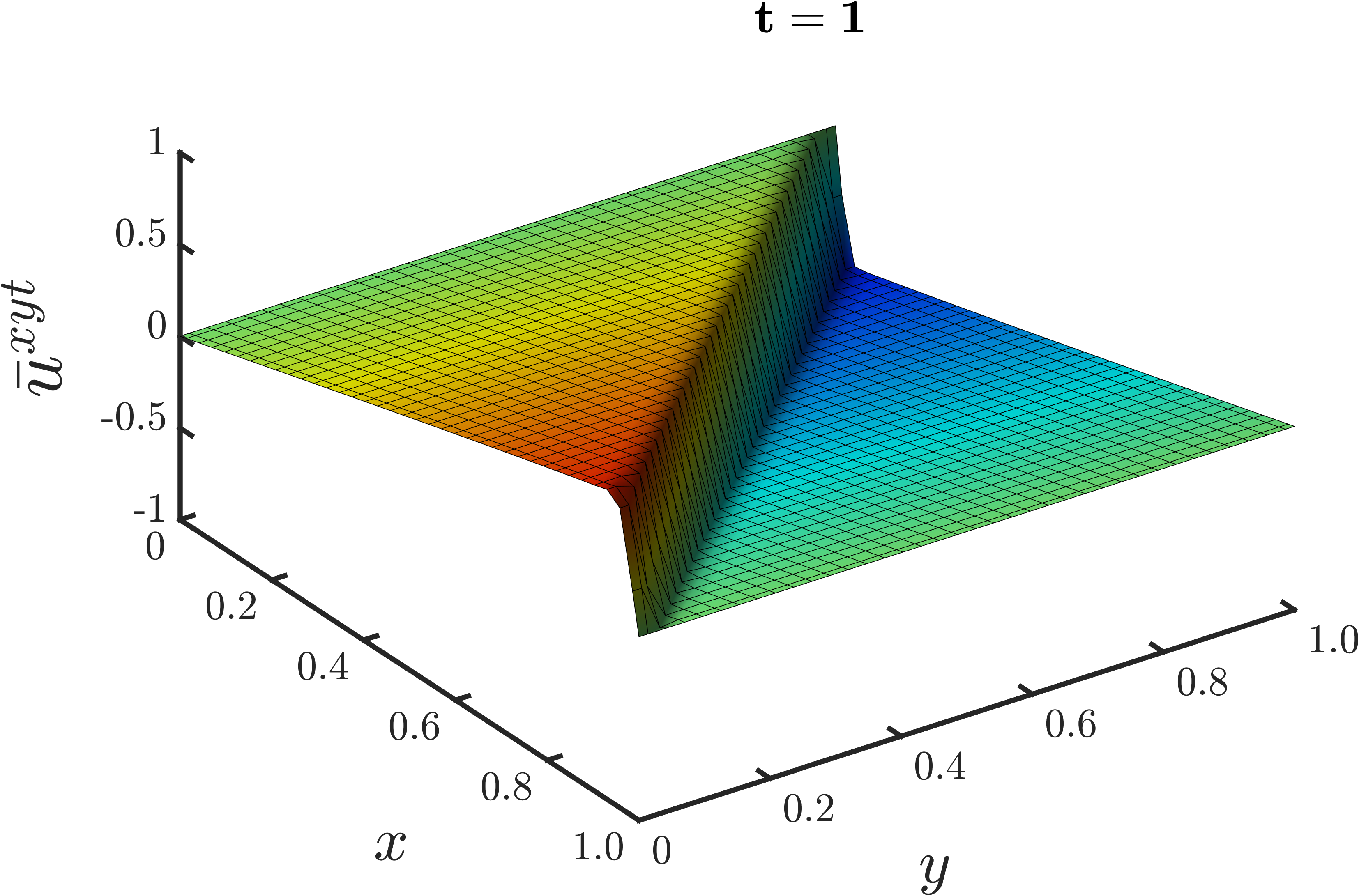}\label{fig05f}}\\
    \caption{Profiles of the velocity component $\bar{u}^{xyt}$ at different time instances for $Re=1000$ on a $36 \times 36$ computational mesh.}
	\label{fig:4}%
\end{figure}
%
%
%\begin{figure}[H]%
%	\centering
%	\subfigure[~]{\includegraphics[width=0.85\linewidth]{Figures/Re_100_plots_at_different_times.jpg}\label{fig5a}}\\%
%	\subfigure[~] {\includegraphics[width=0.85\linewidth]{Figures/Re_1000_plots_at_different_times.jpg}\label{fig5b}}\\
%    \caption{Diagonal ($x=y$) profiles of $u(x,y,t)$ extracted from the 3D RCCNIM solutions at selected time instances, comparing fine and coarse grids for (a) $Re=100$ and (b) $Re=1000$.}
	%\caption{Solutions $u(x,x,t)$ of the RCCNIM at different times for $Re=100$ and $1000$, respectively. (a) $Re=100$ and (b) $Re=1000$.}
%	\label{fig:5}%
%\end{figure}
%
\begin{figure}[!htbp]%
	\centering
	\subfigure[~]{\includegraphics[width=0.5\linewidth]{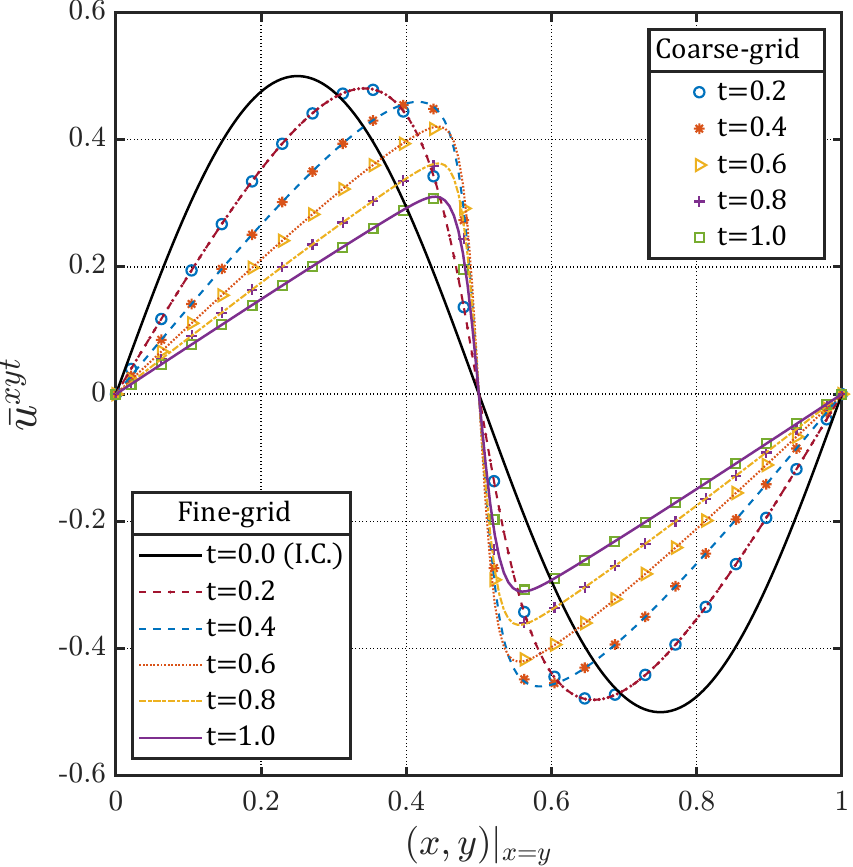}\label{fig06a}}%\\%
	\subfigure[~]{\includegraphics[width=0.5\linewidth]{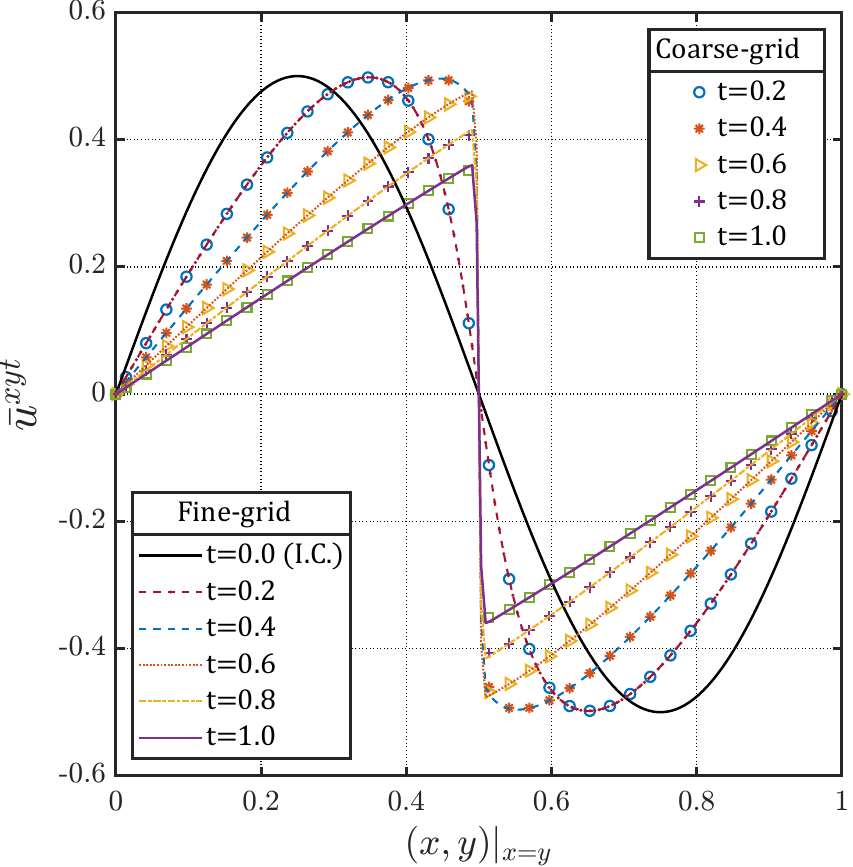}\label{fig06b}}\\
    \caption{Diagonal ($x=y$) profiles of $u(x,y,t)$ extracted from the 3D RCCNIM solutions at selected time instances, comparing fine and coarse grids for (a) $Re=100$ and (b) $Re=1000$.}
	\label{fig:5}%
\end{figure}
\begin{table}[!t]
\begin{center}\renewcommand{\arraystretch}{2}
\caption{Comparison between present RCCNIM and analytical \citep{gao2017analytical} results for $u(x,y,t)$ at selected points ($x,y$) in the computational domain for $Re=100$. The deviations of the RCCNIM results from the analytical solution are listed as $\delta$.}
\label{tab:5}
\resizebox{1\textwidth}{!}{
\begin{tabular}{|c|c|c | c | c | c | c | c | c | c|}
\hline
$x$ & $y$ & \multicolumn{2}{c|}{$t=0.25$} & \multicolumn{2}{c|}{$t=0.5$}  & \multicolumn{2}{c|}{$t=0.75$}  & \multicolumn{2}{c|}{$t=1.0$} \\ \cline{3-10}
&  & Ref \citep{gao2017analytical} & RCCNIM & Ref \citep{gao2017analytical} & RCCNIM & Ref \citep{gao2017analytical} & RCCNIM  & Ref \citep{gao2017analytical} & RCCNIM  \\ \hline

0.25 & 0.25 & 0.39355 & 0.39345 & 0.29118 & 0.29119 & 0.22738 & 0.22735 & 0.18580 & 0.18586 \\ \hline
0.50 & 0.25 & 0.68618 & 0.68532 & 0.56055 & 0.56011 & 0.44800 & 0.44784 & 0.36879 & 0.36872 \\ \hline
0.75 & 0.25 & 0.39355 & 0.39389 & 0.29118 & 0.29132 & 0.22738 & 0.22741 & 0.18580 & 0.18583 \\ \hline
0.25 & 0.50 & 0.26195 & 0.26085 & 0.26198 & 0.26172 & 0.21804 & 0.21819 & 0.18186 & 0.18177 \\ \hline
0.50 & 0.50 & 0.00000 & 0.00000 & 0.00000 & 0.00000 & 0.00000 & 0.00000 & 0.00000 & 0.00000 \\ \hline
0.75 & 0.50 & -0.26195 & -0.26085 & -0.26198 & -0.26172 & -0.21804 & -0.21819 & -0.18186 & -0.18177 \\ \hline
0.25 & 0.75 & -0.39355 & -0.39389 & -0.29118 & -0.29132 & -0.22738 & -0.22741 & -0.18580 & -0.18583 \\ \hline
0.50 & 0.75 & -0.68618 & -0.68532 & -0.56055 & -0.56011 & -0.44800 & -0.44784 & -0.36879 & -0.36872 \\ \hline
0.75 & 0.75 & -0.39355 & -0.39345 & -0.29118 & -0.29119 & -0.22738 & -0.22735 & -0.18580 & -0.18586 \\ \hline
\multicolumn{2}{|c}{$\delta_{min}$ (\%)} & \multicolumn{2}{|c}{0.0254} & \multicolumn{2}{|c}{0.0034} & \multicolumn{2}{|c}{0.0132} & \multicolumn{2}{|c|}{0.0161} \\ \hline
\multicolumn{2}{|c}{$\delta_{max}$ (\%)} & \multicolumn{2}{|c}{0.4202} & \multicolumn{2}{|c}{0.0992} & \multicolumn{2}{|c}{0.0688} & \multicolumn{2}{|c|}{0.0495} \\ \hline
\end{tabular}}
\end{center}
\end{table}
\begin{table}[!t]
\begin{center}\renewcommand{\arraystretch}{2}
\caption{Comparison between present RCCNIM and analytical \citep{gao2017analytical} results for $u(x,y,t)$ at selected points ($x,y$) in the computational domain for $Re=1000$. The deviations of the RCCNIM results from the analytical solution are listed as $\delta$.}
\label{tab:6}
\resizebox{1\textwidth}{!}{
\begin{tabular}{|c|c|c | c | c | c | c | c | c | c|}
\hline
$x$ & $y$ & \multicolumn{2}{c|}{t=0.25} & \multicolumn{2}{c|}{t=0.5}  & \multicolumn{2}{c|}{t=0.75}  & \multicolumn{2}{c|}{t=1.0} \\ \cline{3-10}
&  & Ref \citep{gao2017analytical} & RCCNIM & Ref \citep{gao2017analytical} & RCCNIM & Ref \citep{gao2017analytical} & RCCNIM  & Ref \citep{gao2017analytical} & RCCNIM  \\ \hline

0.25 & 0.25 & 0.40215 & 0.40229 & 0.29673 & 0.29682 & 0.23092 & 0.23098 & 0.18824 & 0.18828 \\ \hline
0.50 & 0.25 & 0.71299 & 0.71278 & 0.57601 & 0.57608 & 0.45667 & 0.45676 & 0.37442 & 0.37449 \\ \hline
0.75 & 0.25 & 0.40215 & 0.40239 & 0.29673 & 0.29687 & 0.23092 & 0.23087 & 0.18824 & 0.18829 \\ \hline
0.25 & 0.50 & 0.28204 & 0.28153 & 0.27288 & 0.27281 & 0.22362 & 0.22363 & 0.18529 & 0.18531 \\ \hline
0.50 & 0.50 & 0.00000 & 0.00000 & 0.00000 & 0.00000 & 0.00000 & 0.00000 & 0.00000 & 0.00000 \\ \hline
0.75 & 0.50 & -0.28204 & -0.28153 & -0.27288 & -0.27283 & -0.22362 & -0.22362 & -0.18529 & -0.18530 \\ \hline
0.25 & 0.75 & -0.40215 & -0.40239 & -0.29673 & -0.29686 & -0.23092 & -0.23087 & -0.18824 & -0.18829 \\ \hline
0.50 & 0.75 & -0.71299 & -0.71278 & -0.57601 & -0.57608 & -0.45667 & -0.45676 & -0.37442 & -0.37449 \\ \hline
0.75 & 0.75 & -0.40215 & -0.40229 & -0.29673 & -0.29683 & -0.23092 & -0.23098 & -0.18824 & -0.18828 \\ \hline
\multicolumn{2}{|c}{$\delta_{min}$ (\%)} & \multicolumn{2}{|c}{0.0348} & \multicolumn{2}{|c}{0.0257} & \multicolumn{2}{|c}{0.0045} & \multicolumn{2}{|c|}{0.0106} \\ \hline
\multicolumn{2}{|c}{$\delta_{max}$ (\%)} & \multicolumn{2}{|c}{0.1808} & \multicolumn{2}{|c}{0.0472} & \multicolumn{2}{|c}{0.0390} & \multicolumn{2}{|c|}{0.0266} \\ \hline
\end{tabular}}
\end{center}
\end{table}

The problem under consideration is symmetric with respect to the velocity components $u(x,y,t)$ and $v(x,y,t)$, i.e., $u(x,y,t) = v(x,y,t)$. Therefore, only the results corresponding to $u(x,y,t)$ are presented and discussed. The RCCNIM solutions are obtained on uniform meshes of $24 \times 24$ and $36 \times 36$ for Reynolds numbers $Re = 100$ and $Re = 1000$, respectively, with a fixed time step size of $\Delta t = 0.001$. The computed solutions are depicted in \figs\ref{fig:3} and \ref{fig:4}. For $Re = 100$ (\fig\ref{fig:3}), the solution remains smooth and does not present significant numerical difficulty. In contrast, the case of $Re = 1000$ (\fig\ref{fig:4}) is more challenging due to increased nonlinearity, and is therefore considered to evaluate the robustness of the proposed method. The results clearly indicate that the numerical solution remains stable and free from spurious oscillations throughout the simulation. Since an exact analytical solution is not available for this problem, a highly refined grid solution is treated as a reference. Specifically, a fine mesh of $160 \times 160$ is used, and the corresponding results are compared with those obtained on coarser grids in Figure~5. It is evident that the coarse grid solutions ($24 \times 24$ for $Re=100$ and $36 \times 36$ for $Re=1000$) show excellent agreement with the fine grid results. Furthermore, to strengthen the validation of the RCCNIM approach, comparisons are made with the benchmark results reported in \citep{gao2017analytical}. These comparisons are performed under identical grid resolutions and time-step conditions. The close agreement between the RCCNIM results and the benchmark data is further confirmed by the numerical values presented at selected points in \tabs\ref{tab:5} and \ref{tab:6}.   
\subsubsection{Example 4: Two-dimensional Burgers' equations with source terms: Decaying shock problem}
\label{sec:3.1.3}
The decaying shock problem is considered in this study due to its frequent use in the literature as a benchmark for demonstrating the inherent upwinding characteristics of numerical schemes \citep{25Wescott_2001}. In the present work, this problem is similarly employed to examine and confirm that the proposed scheme also possesses inherent upwinding properties.
%
%
%\begin{figure}[H]%
%	\centering
%	\subfigure[~]{\includegraphics[width=0.5\linewidth]{Figures/u_at_t_0_01.jpg}\label{fig5a}}%\\%
%	\subfigure[~] {\includegraphics[width=0.5\linewidth]{Figures/u_at_t_0_5.jpg}\label{fig5b}}\\
%	\subfigure[~] {\includegraphics[width=0.5\linewidth]{Figures/u_at_t_5.jpg}\label{fig5b}}\\
%	\caption{Time evolution of $\bar{u}^{xyt}$ obtained using RCCNIM. (a) $t=0.01$, (b) $t=0.5$ and (c) $t=5.0$.}
%	\label{fig:6}%
%\end{figure}
%
%
%\begin{figure}[H]%
%	\centering
%	\subfigure[~]{\includegraphics[width=0.5\linewidth]{Figures/v_at_t_0_01.jpg}\label{fig5a}}%\\%
%	\subfigure[~] {\includegraphics[width=0.5\linewidth]{Figures/v_at_t_0_5.jpg}\label{fig5b}}\\
%	\subfigure[~] {\includegraphics[width=0.5\linewidth]{Figures/v_at_t_5.jpg}\label{fig5b}}\\
%	\caption{Time evolution of $\bar{v}^{xyt}$ obtained using RCCNIM. (a) $t=0.01$, (b) $t=0.5$ and (c) $t=5.0$.}
%	\label{fig:7}%
%\end{figure}
%
%
%
\begin{figure}[!htbp]
\centering
\vfill

\begin{minipage}{0.9\linewidth}
	\centering
	\includegraphics[width=\linewidth]{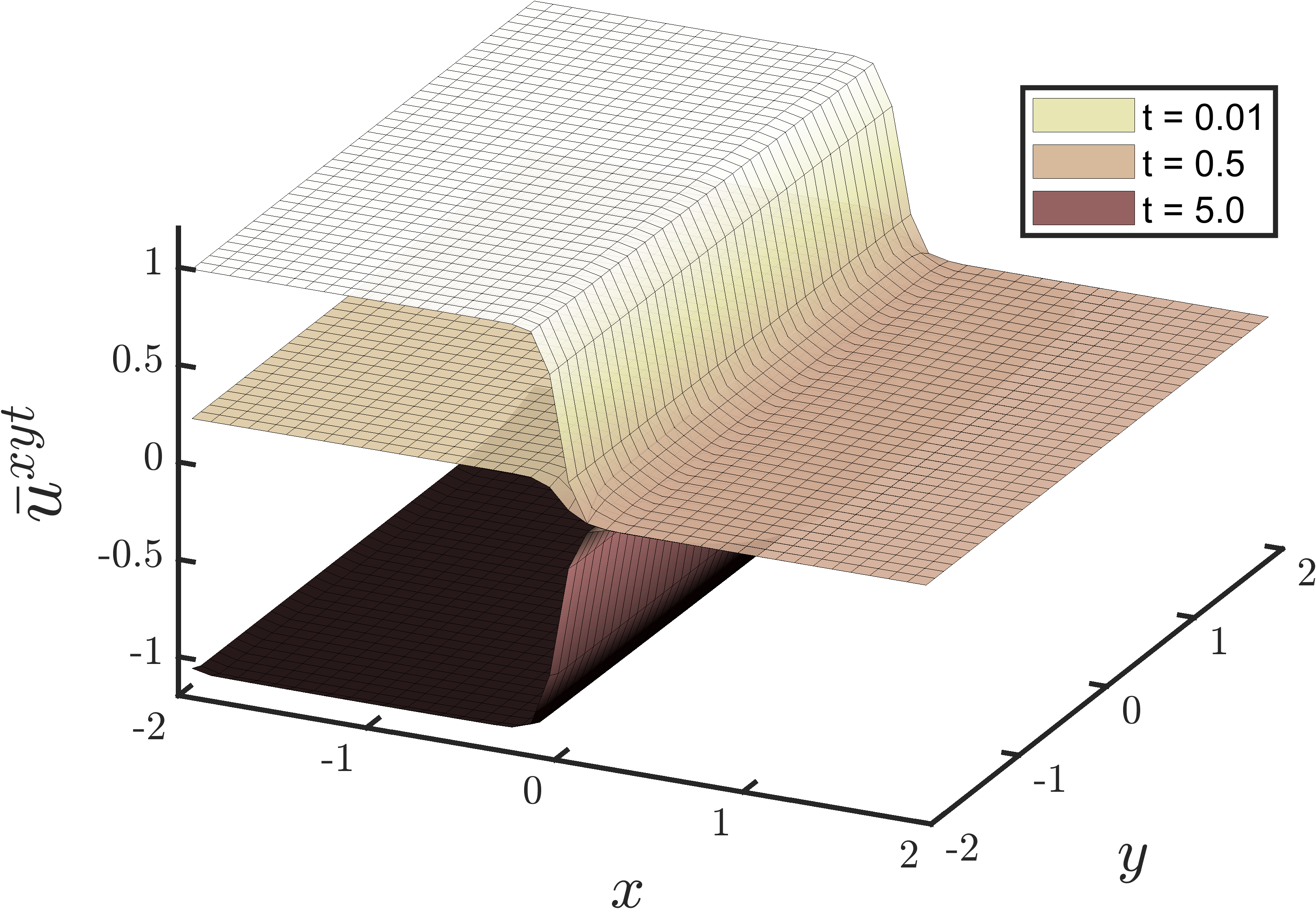}
	
	\vspace{0.2cm}
	{\Large (a)}
\end{minipage}

\vspace{0.5cm}

\begin{minipage}{0.9\linewidth}
	\centering
	\includegraphics[width=\linewidth]{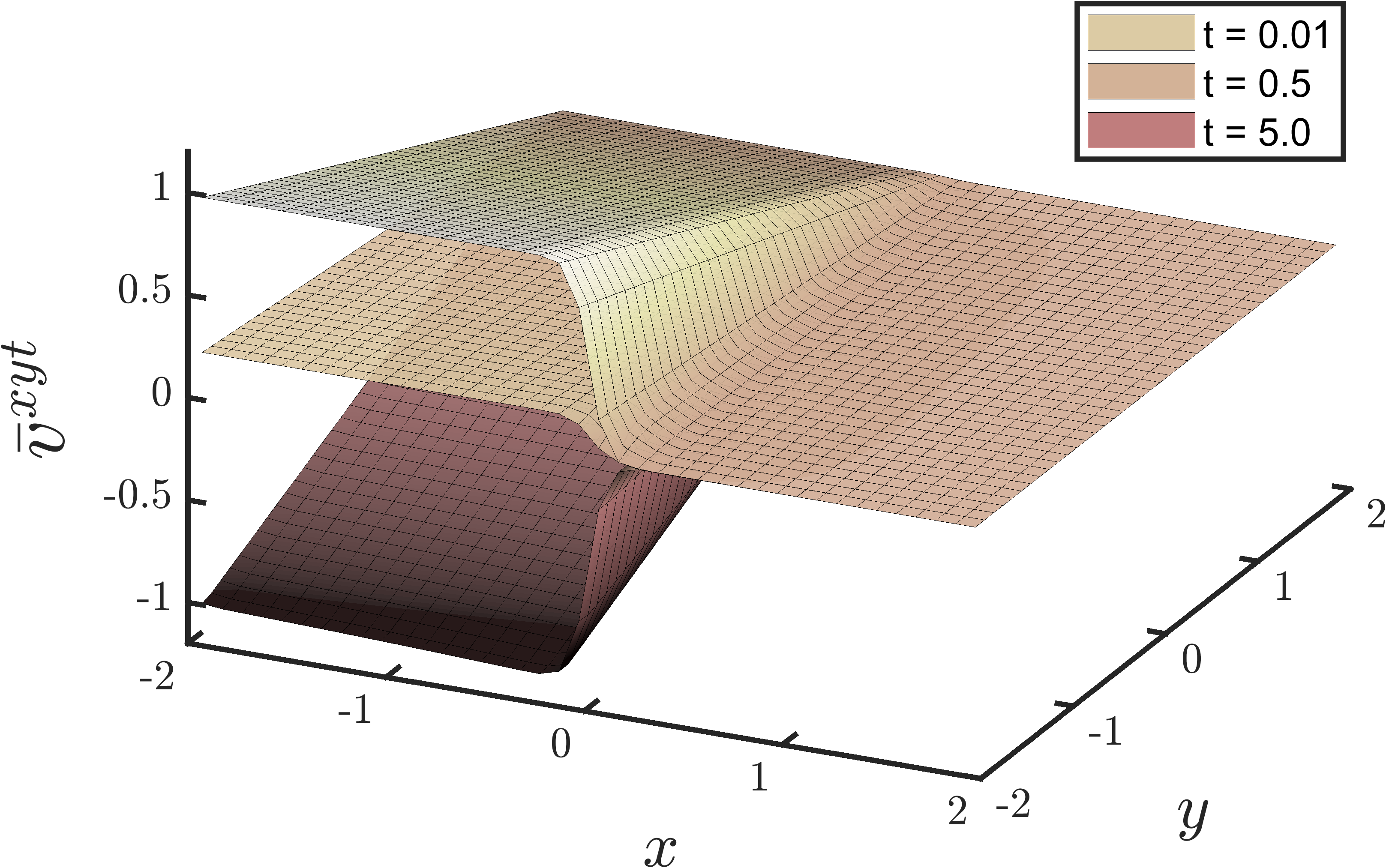}
	
	\vspace{0.2cm}
	{\Large (b)}
\end{minipage}

\vfill

\caption{Time evolution of $\bar{u}^{xyt}$ and $\bar{v}^{xyt}$ obtained using RCCNIM. (a) $\bar{u}^{xyt}$, and (b) $\bar{v}^{xyt}$.}
\label{fig:7}%

\end{figure}
The corresponding exact solution for this problem is given as follows:
\begin{gather}
	\begin{split}
		u\left(x,y,t\right)=\frac{1}{2}\left(1-\tanh\left[\frac{xRe}{4}\right]\right)\left( -1+2e^{-t}\right) 
	\end{split}
\label{eq:072}
\end{gather}
\begin{gather}
	\begin{split}
		v\left(x,y,t\right)=\frac{1}{2}\left(1-\tanh\left[\frac{xRe}{4}\right]\right)\left(\frac{1}{2}-\frac{1}{4}y\right)\left( -1+2e^{-t}\right) 
	\end{split}
\label{eq:073}
\end{gather}
The source terms corresponding to the exact solution defined in \eqn\eqref{eq:072}-\eqref{eq:073} are expressed as follows:
\begin{gather}
\begin{aligned}
f_x(x,y,t) &= -e^{-t}\left(1 - \tanh\left(\tfrac{xRe}{4}\right)\right) \\
&\quad - \frac{Re}{16}\left(-1 + 2e^{-t}\right)^2
\left(1 - \tanh\left(\tfrac{xRe}{4}\right)\right)
\operatorname{sech}^2\left(\tfrac{xRe}{4}\right) \\
&\quad - \frac{Re}{16}\left(-1 + 2e^{-t}\right)
\operatorname{sech}^2\left(\tfrac{xRe}{4}\right)
\tanh\left(\tfrac{xRe}{4}\right)
\end{aligned}
\label{eq:074}
\end{gather}
\begin{gather}
\begin{aligned}
f_y(x,y,t) &= -e^{-t}
\left(1 - \tanh\left(\tfrac{xRe}{4}\right)\right)
\left(\tfrac{1}{2} - \tfrac{y}{4}\right) \\
&\quad - \frac{Re}{16}\left(-1 + 2e^{-t}\right)^2
\left(1 - \tanh\left(\tfrac{xRe}{4}\right)\right)
\operatorname{sech}^2\left(\tfrac{xRe}{4}\right)
\left(\tfrac{1}{2} - \tfrac{y}{4}\right) \\
&\quad - \frac{1}{16}\left(-1 + 2e^{-t}\right)^2
\left(1 - \tanh\left(\tfrac{xRe}{4}\right)\right)^2
\left(\tfrac{1}{2} - \tfrac{y}{4}\right) \\
&\quad - \frac{Re}{16}\left(-1 + 2e^{-t}\right)
\operatorname{sech}^2\left(\tfrac{xRe}{4}\right)
\tanh\left(\tfrac{xRe}{4}\right)
\left(\tfrac{1}{2} - \tfrac{y}{4}\right)
\end{aligned}
\label{eq:075}
\end{gather}

In this test case, a shock is initially located at $x = 0$, which weakens over time and later reappears. The evolution of the velocity components $u$ and $v$ at $t = 0.01$, $0.5$, and $5\ \text{s}$ is presented in \fig\ref{fig:7}. For the region $x > 0$, the velocity remains nearly unchanged with time. In contrast, for $x < 0$, both velocity components decrease as time progresses, eventually becoming negative and approaching steady-state values as $t \to \infty$. The resulting steady-state solution corresponds to a reflection of the initial condition about the $x\text{-}y$ plane. The computations are carried out on a $40 \times 40$ grid with a time step $\Delta t = 0.01$ and viscosity $\nu = 0.02$ in the domain $-2 \le x \le 2$ and $-2 \le y \le 2$. A detailed analysis of this problem is presented in \citep{25Wescott_2001}, where it is demonstrated that the numerical scheme exhibits an inherent upwinding mechanism, provided the solution remains sufficiently smooth even under flow reversal. The computed solutions are illustrated in \fig\ref{fig:7}. It is evident from \fig\ref{fig:7} that the numerical solution retains its smoothness despite the reversal of flow direction; consequently, the proposed scheme continues to preserve its inherent upwinding characteristics. 

%--------------------------------------
\section{Concluding remarks}
%--------------------------------------
\label{sec:4}
\noindent 
The development of the improved RCCNIM formulation was driven by the objective of reducing the computational cost associated with MCCNIM, while preserving its superior accuracy relative to conventional nodal approaches. A key modification involves the treatment of the convective term through linearization based on the solution from the previous time level, which facilitates the construction of the transverse spatially integrated equations with explicitly known coefficients. As a result, the coefficient matrix corresponding to the discretized system needs to be assembled only once per time step, leading to enhanced computational efficiency. This improvement is reflected in a noticeable reduction in CPU time without any compromise in solution accuracy. The two-dimensional formulation of RCCNIM has been validated using multiple test cases encompassing various boundary condition types, and it consistently exhibits inherent upwinding characteristics. Furthermore, the methodology presented here provides a systematic framework that can be extended to develop efficient numerical schemes for a broader class of nonlinear, time-dependent partial differential equations.
%
%--------------- Acknowledgments
\section*{Conflict of interest}
 The authors declare that they have no known competing financial interests or personal relationships that could have appeared to influence the work reported in this paper.
%
%--------------- Acknowledgments
\section*{Acknowledgments}
 Nadeem Ahmed acknowledges the support received under the Institute Post-Doctoral Fellow (IPDF) scheme of the Indian Institute of Technology Roorkee, Roorkee, India.
%
%--------------- Nomenclature
\begin{spacing}{1.5}
\input{Nomenclature.tex}

%\renewcommand{\nompreamble}{\vspace{1em}\fontsize{10}{8pt}\selectfont}
{\printnomenclature[5em]}
%\printnomenclature
\end{spacing}
%
%\printglossaries 
%
%
%--------------- Bibliography
%
%\begin{thebibliography}{0000}
%\bibliographystyle{plainnat}
%\bibliographystyle{elsarticle/elsarticle-harv}\biboptions{authoryear}
\bibliography{references}
%
% Bibliographic references with the natbib package:
% Parenthetical: \citep{Bai92} produces (Bailyn 1992).
% Textual: \citet{Bai95} produces Bailyn et al. (1995).
% An affix and part of a reference:
%   \citep[e.g.][Ch. 2]{Bar76}
%   produces (e.g. Barnes et al. 1976, Ch. 2).
% 
% \bibitem[Names(Year)]{label} or \bibitem[Names(Year)Long names]{label}.
% (\harvarditem{Name}{Year}{label} is also supported.)
%\bibitem[]{}
%
%\input{references.tex}
%
%\end{thebibliography}
%
%=============================================================
\newpage
\appendix
\section{RCCNIM coefficients}\label{app:A}
\renewcommand{\thesection}{\Alph{section}}
\subsection{RCCNIM coefficients for 2D Burgers’ equation}
The coefficients corresponding to the averaged flux formulations, as defined in \eqn\eqref{eq:033}, \eqref{eq:034}, and \eqref{eq:039}–\eqref{eq:044}, are given below:
\begin{gather}
	\begin{split}
		A_{31}=\frac{2ae^{Ru_{i,j,\ell}}({{\bar{u}}_{i,j,\ell}^p})^2{Re}}{1-e^{Ru_{i,j,\ell}}+2ae^{Ru_{i,j,\ell}}{\bar{u}}_{i,j,\ell}^pRe} 
	\end{split}
\end{gather}
\begin{gather}
	\begin{split}
		A_{32}=\frac{-1+e^{Ru_{i,j,\ell}}-2ae^{Ru_{i,j,\ell}}{\bar{u}}_{i,j,\ell}^pRe + 2a^2e^{Ru_{i,j,\ell}}({{\bar{u}}_{i,j,\ell}^p})^2{Re}^2}{{\bar{u}}_{i,j,\ell}^pRe\left(-1+e^{Ru_{i,j,\ell}}-2ae^{Ru_{i,j,\ell}}{\bar{u}}_{i,j,\ell}^pRe\right)} 
	\end{split}
\end{gather}
\begin{gather}
	\begin{split}
		A_{51}=\frac{2a({{\bar{u}}_{i,j,\ell}^p})^2{Re}}{1-e^{Ru_{i,j,\ell}}+2a{\bar{u}}_{i,j,\ell}^pRe} 
	\end{split}
\end{gather}
\begin{gather}
	\begin{split}
		A_{52}=\frac{1-e^{Ru_{i,j,\ell}}+2a{\bar{u}}_{i,j,\ell}^pRe+2a^2({{\bar{u}}_{i,j,\ell}^p})^2{Re}^2}{{\bar{u}}_{i,j,\ell}^pRe\left(1-e^{Ru_{i,j,\ell}}+2a{\bar{u}}_{i,j,\ell}^pRe\right)} 
	\end{split}
\end{gather}
\begin{gather}
	\begin{split}
		B_{31}=\frac{2be^{Rv_{i,j,\ell}}({{\bar{v}}_{i,j,\ell}^p})^2{Re}}{1-e^{Rev_{i,j,\ell}}+2be^{Rv_{i,j,\ell}}{\bar{v}}_{i,j,\ell}^pRe} 
	\end{split}
\end{gather}
\begin{gather}
	\begin{split}
		B_{32}=\frac{-1+e^{Rv_{i,j,\ell}}-2be^{Rv_{i,j,\ell}}{\bar{v}}_{i,j,\ell}^pRe+2b^2e^{Rv_{i,j,\ell}}({{\bar{v}}_{i,j,\ell}^p})^2{Re}^2}{{\bar{v}}_{i,j,\ell}^pRe\left(-1+e^{Rv_{i,j,\ell}}-2be^{Rv_{i,j,\ell}}{\bar{v}}_{i,j,\ell}^pRe\right)} 
	\end{split}
\end{gather}
\begin{gather}
	\begin{split}
		B_{51}=\frac{2b({{\bar{v}}_{i,j,\ell}^p})^2{Re}}{1-e^{Rv_{i,j,\ell}}+2b{\bar{v}}_{i,j,\ell}^pRe} 
	\end{split}
\end{gather}
\begin{gather}
	\begin{split}
		B_{52}=\frac{1-e^{Rv_{i,j,\ell}}+2b{\bar{v}}_{i,j,\ell}^pRe+2b^2({{\bar{v}}_{i,j,\ell}^p})^2{Re}^2}{{\bar{v}}_{i,j,\ell}^pRe\left(1-e^{Rv_{i,j,\ell}}+2b{\bar{v}}_{i,j,\ell}^pRe\right)} 
	\end{split}
\end{gather}
Moreover, the final coefficients given in \eqn\eqref{eq:059}–\eqref{eq:062}, are presented below:
\begin{gather}
	\begin{split}
	 C_{31}=\frac{\left(-\frac{A_{32}A_{51,i+1}}{A_{31}-A_{51,i+1}}-\frac{A_{31,i-1}A_{52}}{A_{31,i-1}-A_{51}}\right)}{2a} 
	\end{split}
\end{gather}
\begin{gather}
	\begin{split}
	C_{32}=\frac{A_{32,i-1}A_{51}}{2a\left(A_{31,i-1}-A_{51}\right)} 
	\end{split}
\end{gather}
\begin{gather}
	\begin{split}
		C_{33}=\frac{A_{31}A_{52,i+1}}{2a\left(A_{31}-A_{51,i+1}\right)} 
	\end{split}
\end{gather}
\begin{gather}
	\begin{split}
		C_{34}=\frac{\left(-\frac{A_{31,i-1}A_{51}}{A_{31,i-1}-A_{51}}-\frac{A_{31}A_{51,i+1}}{A_{31}-A_{51,i+1}}\right)}{2a} 
	\end{split}
\end{gather}
\begin{gather}
	\begin{split}
		C_{35}=\frac{A_{31,i-1}A_{51}}{2a\left(A_{31,i-1}-A_{51}\right)} 
	\end{split}
\end{gather}
\begin{gather}
	\begin{split}
		C_{36}=\frac{A_{31}A_{51,i+1}}{2a\left(A_{31}-A_{51,i+1}\right)} 
	\end{split}
\end{gather}
\begin{gather}
	\begin{split}
	 C_{51}=\frac{\left(\frac{A_{32}}{A_{31}-A_{51,i+1}}+\frac{A_{52}}{A_{31,i-1}-A_{51}}\right)}{2a} 
	\end{split}
\end{gather}
\begin{gather}
	\begin{split}
	C_{52}=-\frac{A_{32,i-1}}{2a\left(A_{31,i-1}-A_{51}\right)} 
	\end{split}
\end{gather}
\begin{gather}
	\begin{split}
		C_{53}=-\frac{A_{52,i+1}}{2a\left(A_{31}-A_{51,i+1}\right)} 
	\end{split}
\end{gather}
\begin{gather}
	\begin{split}
		C_{54}=\frac{\left(\frac{A_{51}}{A_{31,i-1}-A_{51}}+\frac{A_{31}}{A_{31}-A_{51,i+1}}\right)}{2a} 
	\end{split}
\end{gather}
\begin{gather}
	\begin{split}
		C_{55}=-\frac{A_{31,i-1}A_{51}}{2a\left(A_{31,i-1}-A_{51}\right)} 
	\end{split}
\end{gather}
\begin{gather}
	\begin{split}
		C_{56}=-\frac{A_{51,i+1}}{2a\left(A_{31}-A_{51,i+1}\right)} 
	\end{split}
\end{gather}
\begin{gather}
	\begin{split}
	 C_{71}=\frac{\left(-\frac{B_{32}B_{51,j+1}}{B_{31}-B_{51,j+1}}-\frac{B_{31,j-1}B_{52}}{B_{31,j-1}-B_{51}}\right)}{2b} 
	\end{split}
\end{gather}
\begin{gather}
	\begin{split}
	C_{72}=\frac{B_{32,j-1}B_{51}}{2b\left(B_{31,j-1}-B_{51}\right)} 
	\end{split}
\end{gather}
\begin{gather}
	\begin{split}
		C_{73}=\frac{B_{31}B_{52,j+1}}{2b\left(B_{31}-B_{51,j+1}\right)} 
	\end{split}
\end{gather}
\begin{gather}
	\begin{split}
		C_{74}=\frac{\left(-\frac{B_{31,j-1}B_{51}}{B_{31,j-1}-B_{51}}-\frac{B_{31}B_{51,j+1}}{B_{31}-B_{51,j+1}}\right)}{2b} 
	\end{split}
\end{gather}
\begin{gather}
	\begin{split}
		C_{75}=\frac{B_{31,j-1}B_{51}}{2b\left(B_{31,j-1}-B_{51}\right)} 
	\end{split}
\end{gather}
\begin{gather}
	\begin{split}
		C_{76}=\frac{B_{31}B_{51,j+1}}{2b\left(B_{31}-B_{51,j+1}\right)} 
	\end{split}
\end{gather}
\begin{gather}
	\begin{split}
	 C_{91}=\frac{\left(\frac{B_{32}}{B_{31}-B_{51,j+1}}+\frac{B_{52}}{B_{31,j-1}-B_{51}}\right)}{2b} 
	\end{split}
\end{gather}
\begin{gather}
	\begin{split}
	C_{92}=-\frac{B_{32,j-1}}{2b\left(B_{31,j-1}-B_{51}\right)} 
	\end{split}
\end{gather}
\begin{gather}
	\begin{split}
		C_{93}=-\frac{B_{52,j+1}}{2b\left(B_{31}-B_{51,j+1}\right)} 
	\end{split}
\end{gather}
\begin{gather}
	\begin{split}
		C_{94}=\frac{\left(\frac{B_{51}}{B_{31,j-1}-B_{51}}+\frac{B_{31}}{B_{31}-B_{51,j+1}}\right)}{2b} 
	\end{split}
\end{gather}
\begin{gather}
	\begin{split}
		C_{95}=-\frac{B_{31,j-1}B_{51}}{2b\left(B_{31,j-1}-B_{51}\right)} 
	\end{split}
\end{gather}
\begin{gather}
	\begin{split}
		C_{96}=-\frac{B_{51,j+1}}{2b\left(B_{31}-B_{51,j+1}\right)} 
	\end{split}
\end{gather}
\subsection{Coefficients of the RCCNIM scheme corresponding to the approximated (node-averaged) convective velocity for 2D Burgers’ equations}
The RCCNIM coefficients corresponding to the node-averaged velocity $u_{i,j,\ell}^p$, as defined in \eqn\eqref{eq:067} for the two-dimensional Burgers’ equation, are given by:
\begin{gather}
	\begin{split}
	F_{71}=\frac{A_{32,i-1}}{2\left(A_{31,i-1}-A_{51}\right)} 
	\end{split}
\end{gather}
\begin{gather}
	\begin{split}
		F_{72}=\frac{A_{32}}{2\left(A_{31}-A_{51,i+1}\right)}-\frac{A_{52}}{2\left(A_{31,i-1}-A_{51}\right)} 
	\end{split}
\end{gather}
\begin{gather}
	\begin{split}
		F_{73}=-\frac{A_{52,i+1}}{2\left(A_{31}-A_{51,i+1}\right)} 
	\end{split}
\end{gather}
\begin{gather}
	\begin{split}
		F_{74}=\frac{A_{31,i-1}}{2\left(A_{31,i-1}-A_{51}\right)} 
	\end{split}
\end{gather}
\begin{gather}
	\begin{split}
	F_{75}=-\frac{A_{51}}{2\left(A_{31,i-1}-A_{51}\right)}+\frac{A_{31}}{2\left(A_{31}-A_{51,i+1}\right)} 
	\end{split}
\end{gather}
\begin{gather}
	\begin{split}
		F_{76}=-\frac{A_{51,i+1}}{2\left(A_{31}-A_{51,i+1}\right)} 
	\end{split}
\end{gather}
\begin{gather}
	\begin{split}
	F_{91}=\frac{B_{32,j-1}}{2\left(B_{31,j-1}-B_{51}\right)} 
	\end{split}
\end{gather}
\begin{gather}
	\begin{split}
		F_{92}=\frac{B_{32}}{2\left(B_{31}-B_{51,j+1}\right)}-\frac{B_{52}}{2\left(B_{31,j-1}-B_{51}\right)} 
	\end{split}
\end{gather}
\begin{gather}
	\begin{split}
		F_{93}=-\frac{B_{52,j+1}}{2\left(B_{31}-B_{51,j+1}\right)} 
	\end{split}
\end{gather}
\begin{gather}
	\begin{split}
	F_{94}=\frac{B_{31,j-1}}{2\left(B_{31,j-1}-B_{51}\right)} 
	\end{split}
\end{gather}
\begin{gather}
	\begin{split}
		F_{95}=-\frac{B_{51}}{2\left(B_{31,j-1}-B_{51}\right)}+\frac{B_{31}}{2\left(B_{31}-B_{51,j+1}\right)} 
	\end{split}
\end{gather}
\begin{gather}
	\begin{split}
		F_{96}=-\frac{B_{51,j+1}}{2\left(B_{31}-B_{51,j+1}\right)} 
	\end{split}
\end{gather}
The coefficients corresponding to the remaining node-averaged velocity component, namely $v_{i,j,\ell}^p$, can be obtained through an analogous derivation procedure. It is important to emphasize that, for problems with fully Dirichlet boundary conditions, the same coefficient structure remains applicable for $y$-momentum equation. However, in the presence of mixed boundary conditions, as considered in Example 3, the coefficients associated with the $y$-momentum equation must be reformulated to consistently incorporate the prescribed boundary constraints.
%=============================================================
%
%\addcontentsline{toc}{section}{References} 
%\clearpage\listoftables
%
%\clearpage\listoffigures
%
%\clearpage
%
%\input{tables.tex}
%
%\clearpage
%
%\renewcommand{\thesubfigure}{(\roman{subfigure})}
%
%\input{figures.tex}
%
\end{document}

%% file: nomenclature.tex
\fontsize{10}{10pt}\selectfont
%\fontsize{10}{10pt}\selectfont
%-------------------------------------
%\section*{Nomenclature}
%\noindent 
%------------------------------------------------------------------------
\renewcommand{\nomgroup}[1]{%
	\ifthenelse{\equal{#1}{0}}{~}{%{\item[\textbf{Variables}]}
		\ifthenelse{\equal{#1}{A}}{\item[\textbf{Abbreviations}]}{%
			\ifthenelse{\equal{#1}{G}}{\item[\textbf{Greek Symbols}]}{%
				\ifthenelse{\equal{#1}{D}}{\item[\textbf{Dimensionless Groups / Constants}]}{%
					\ifthenelse{\equal{#1}{S}}{\item[\textbf{Subcripts / Superscripts}]}
					{}}}}}
}
%------------------------------------------------------------------------
%
%%%%%%%%%%% Abbreviations
\nomenclature[A]{CCNIM}{cell-centered nodal integral method}
\nomenclature[A]{CN-4PU}{Crank-Nicolson 4-point upwind}
\nomenclature[A]{NIM}{nodal integral method}
\nomenclature[A]{MCCNIM}{modified cell-centered nodal integral method}
\nomenclature[A]{M$^2$NIM}{modified-modified nodal integral method}
\nomenclature[A]{DAE}{differential-algebraic equation}
\nomenclature[A]{ODE}{ordinary differential equation}
\nomenclature[A]{PDE}{partial differential equation}
\nomenclature[A]{TIP}{transverse integration process}
%
%%%%%%%%%%% List of Symbols
\nomenclature[0a]{$a$}{half of the node size in the $x$-direction}
\nomenclature[0b]{$b$}{half of the node size in the $y$-direction}
\nomenclature[0j]{$\bar{J}$}{surface-averaged flux}	
\nomenclature[0s]{$\bar{S}$}{pseudo-source term}	
\nomenclature[0u]{${\bar{u}}^{xy}$}{$x$ and $y$-space averaged velocity in the $t$-direction}
\nomenclature[0u]{${\bar{u}}^{xt}$}{$x$-space and time-averaged velocity in the $y$-direction}
\nomenclature[0u]{${\bar{u}}^{yt}$}{$y$-space and time-averaged velocity in the $x$-direction}
\nomenclature[0u]{${\bar{u}}^0$}{approximate convective velocity in the $x$-direction at current timestep}
\nomenclature[0v]{${\bar{v}}^{xy}$}{$x$ and $y$-space averaged velocity in the $t$-direction}
\nomenclature[0v]{${\bar{v}}^{xt}$}{$x$-space and time-averaged velocity in the $y$-direction}
\nomenclature[0v]{${\bar{v}}^{yt}$}{$y$-space and time-averaged velocity in the $x$-direction}
\nomenclature[0v]{${\bar{v}}^0$}{approximate convective velocity in the $y$-direction at current timestep}
%
%------- {Greek symbols}
%
\nomenclature[G]{$\Delta t$}{time step}
\nomenclature[G]{$\Delta x$}{node size in the $x$-direction}
\nomenclature[G]{$\Delta y$}{node size in the $y$-direction}
\nomenclature[G]{$\tau$}{half of the time step}
%
%------- {Dimensionless group}
%
\nomenclature[D]{$Re$}{Reynolds number}
\nomenclature[D]{${Reu}_{i,j}$}{local Reynolds number}
%
%-------  subscripts and superscripts
%\nomenclature[sz]{$0$}{without electroviscous effects}
\nomenclature[s]{$i$}{spatial index in the x-direction }
\nomenclature[s]{$j$}{spatial index in the y-direction}
\nomenclature[s]{$\ell$}{temporal index}
\nomenclature[s]{$t$}{transverse averaging in time}
\nomenclature[s]{$x$}{transverse averaging in the x-direction}
\nomenclature[s]{$y$}{transverse averaging in the y-direction}

%
%\printnomenclature

%% file: references.bib
@Article{25Wescott_2001,
	author  = {Wescott, Bradley L. and Rizwan-uddin},
	journal = {Nuclear Science and Engineering},
	title   = {An efficient formulation of the modified nodal integral method and application to the two-dimensional Burgers’ equation},
	year    = {2001},
	note    = {doi:10.13182/nse01-a2239},
	number  = {3},
	pages   = {293--305},
	volume  = {139},
}

@InProceedings{decker1993block,
  author    = {Decker, William J and Dorning, JJ},
  booktitle = {International topical meeting on mathematical methods and supercomputing in nuclear applications},
  title     = {A block iterative nodal integral method for fluid dynamics problems},
  year      = {1993},
  month     = {Apr},
  pages     = {208--223},
  place     = {Germany},
}

@Article{23Rizwan_uddin_1997,
	author  = {Rizwan-uddin},
	journal = {Numerical Methods for Partial Differential Equations},
	title   = {An improved coarse-mesh nodal integral method for partial differential equations},
	year    = {1997},
	note    = {doi:10.1002/(SICI)1098-2426(199703)13:2\textless 113::AID-NUM1\textgreater 3.0.CO;2-S},
	number  = {2},
	pages   = {113--145},
	volume  = {13},
}

@Article{22Rizwan_uddin_1997,
	author  = {Rizwan-uddin},
	journal = {Computers \& Fluids},
	title   = {A second-order space and time nodal method for the one-dimensional convection-diffusion equation},
	year    = {1997},
	note    = {doi:10.1016/s0045-7930(96)00039-4},
	number  = {3},
	pages   = {233--247},
	volume  = {26},
}

@Article{20Elnawawy_1990,
	author  = {Elnawawy, Osman A. and Valocchi, Albert J. and Ougouag, Abderrafi M.},
	journal = {Water Resources Research},
	title   = {The Cell analytical‐numerical method for solution of the advection‐dispersion equation: Two‐dimensional problems},
	year    = {1990},
	note    = {doi:10.1029/wr026i011p02705},
	number  = {11},
	pages   = {2705--2716},
	volume  = {26},
}

@InProceedings{19Ahmed_2021,
  author     = {Ahmed, N. and Kumar, N. and Singh, S.},
  booktitle  = {14th WCCM-ECCOMAS Congress},
  title      = {Node averaged nodal integral method},
  year       = {2021},
  note       = {doi: https://doi.org/10.23967/wccm-eccomas.2020.219},
  pages      = {1-8},
  publisher  = {CIMNE},
  series     = {WCCM-ECCOMAS 2020},
  volume     = {700},
  collection = {WCCM-ECCOMAS 2020},
}

@Article{18Ahmed_2023,
	author  = {Ahmed, Nadeem and Maurya, Govind and Singh, Suneet},
	journal = {Annals of Nuclear Energy},
	title   = {A novel cell-centered nodal integral method for the convection-diffusion equation},
	year    = {2023},
	note    = {doi:10.1016/j.anucene.2023.109858},
	pages   = {109858},
	volume  = {189},
}

@Article{17Ahmed_2024,
	author  = {Ahmed, Nadeem and Singh, Suneet},
	journal = {Journal of Computational Science},
	title   = {A modified cell-centered nodal integral scheme for the convection-diffusion equation},
	year    = {2024},
	note    = {doi:10.1016/j.jocs.2024.102320},
	pages   = {102320},
	volume  = {80},
}

@Article{kumar2012pressure,
  title={Pressure correction--based iterative scheme for Navier-stokes equations using nodal integral method},
  author={Kumar, Neeraj and Singh, Suneet and Doshi, JB},
  journal={Numerical Heat Transfer, Part B: Fundamentals},
  volume={62},
  number={4},
  pages={264--288},
  year={2012},
  publisher={Taylor \& Francis},
  note = {doi: 10.1080/10407790.2012.709169}
}

@TechReport{13Shober1978,
	author      = {Shober, R. A.},
	institution = {Applied Physics Division. Argonne National Laboratory. Argonne, Illinois},
	title       = {A nodal method for solving transient fewgroup neutron diffusion equations},
	year        = {1978},
	address     = {Argonne, Illinois},
	month       = {jun},
	note        = {https://digital.library.unt.edu/ark:/67531/metadc303752/},
	number      = {ANL-78-51},
}

@InProceedings{11Azmy_1983,
	author    = {Y. Y. Azmy},
	booktitle = {Advances in Reactor Computations, Vol II},
	title     = {A nodal integral approach to the numerical solution of partial differential equations},
	year      = {1983},
	pages     = {893-909},
	publisher = {American Nuclear Society. LaGrange Park, IL},
}

@Article{09Shober_1977,
	author  = {Shober, R. A. and Sims, R. N. and Henry, A. F.},
	journal = {Nuclear Science and Engineering},
	title   = {Two nodal methods for solving time-dependent group diffusion equations},
	year    = {1977},
	note    = {doi:10.13182/nse77-a27392},
	number  = {2},
	pages   = {582--592},
	volume  = {64},
}

@Article{07Lawrence_1986,
	author  = {Lawrence, R. D.},
	journal = {Progress in Nuclear Energy},
	title   = {Progress in nodal methods for the solution of the neutron diffusion and transport equations},
	year    = {1986},
	note    = {doi:10.1016/0149-1970(86)90034-x},
	number  = {3},
	pages   = {271--301},
	volume  = {17},
}

@Article{06Hennart_1986,
	author  = {Hennart, J. P.},
	journal = {SIAM Journal on Scientific and Statistical Computing},
	title   = {A general family of nodal schemes},
	year    = {1986},
	note    = {doi:10.1137/0907018},
	number  = {1},
	pages   = {264--287},
	volume  = {7},
}

@InProceedings{05Dorning_1979,
	author    = {Dorning, J. J.},
	booktitle = {Computational Methods in Nuclear Engineering, Vol. 1},
	title     = {Modern coarse-mesh methods - A development of the 70's},
	year      = {1979},
	pages     = {3.1–3.31},
	publisher = {American Nuclear Society, Williamsburg, Virginia},
}

@Article{chai2020appropriate,
  title={Appropriate stabilized Galerkin approaches for solving two-dimensional coupled Burgers’ equations at high Reynolds numbers},
  author={Chai, Yong and Ouyang, Jie},
  journal={Computers \& Mathematics with Applications},
  volume={79},
  number={5},
  pages={1287--1301},
  year={2020},
  publisher={Elsevier},
  note = {doi: 10.1016/j.camwa.2019.08.036}
}

@Article{zhang2010variational,
  title={Variational multiscale element-free Galerkin method for 2D Burgers’ equation},
  author={Zhang, Lin and Ouyang, Jie and Wang, Xiaoxia and Zhang, Xiaohua},
  journal={Journal of computational physics},
  volume={229},
  number={19},
  pages={7147--7161},
  year={2010},
  publisher={Elsevier},
  note = {doi: 10.1016/j.jcp.2010.06.004}
}

@article{zhang2009element,
  title={Element-free characteristic Galerkin method for Burgers’ equation},
  author={Zhang, Xiao Hua and Ouyang, Jie and Zhang, Lin},
  journal={Engineering Analysis with Boundary Elements},
  volume={33},
  number={3},
  pages={356--362},
  year={2009},
  publisher={Elsevier},
  note = {doi: 10.1016/j.enganabound.2008.07.001}
}

@Article{gao2017analytical,
  title={An analytical solution for two and three dimensional nonlinear Burgers' equation},
  author={Gao, Q and Zou, MY},
  journal={Applied Mathematical Modelling},
  volume={45},
  pages={255--270},
  year={2017},
  publisher={Elsevier},
  note = {doi: 10.1016/j.apm.2016.12.018}
}

@Article{ahmed2025efficient,
  title={Efficient cell-centered nodal integral method for multi-dimensional Burgers’ equations},
  author={Ahmed, Nadeem and Bharti, Ram Prakash and Singh, Suneet},
  journal={Engineering with Computers},
  pages={3611--3649},
  volume={41},  
  year={2025},
  publisher={Springer},
  note = {doi: 10.1007/s00366-025-02177-1}
}

@Article{Ahmed_Diffusion_2024,
	author  = {Nadeem Ahmed and Uma Siva Sankar Kompalli and Suneet Singh},
	journal = {Sadhana},
	title   = {An improved cell-centered nodal integral approach for the transient diffusion equation},
	pages   = {1--14},
        year    = {2024},
        publisher = {Springer Nature},
	note = {doi: 10.1007/s12046-025-02728-8}
}

@Article{jarrah2024nodal,
	title={Nodal integral method to solve the two-dimensional, time-dependent, incompressible Navier-Stokes equations in curvilinear coordinates},
	author={Jarrah, Ibrahim and Rizwan-Uddin},
	journal={Computers \& Mathematics with Applications},
	volume={158},
	pages={219-243},
	year={2024},
	note = {doi: 10.1016/j.camwa.2024.02.009},
	publisher={Elsevier}
}

@Inproceedings{gander2022new,
	title={A new nodal integration method for Helmholtz problems based on domain decomposition techniques},
	author={Gander, Martin J and Kumar, Niteen},
	booktitle={International Conference on Domain Decomposition Methods},
	pages={199-206},
	year={2022},
	organization={Springer}
}

@Article{zhou2016general,
	title={General nodal expansion method for multi-dimensional steady and transient convection-diffusion equation},
	author={Zhou, Xiafeng and Guo, Jiong and Li, Fu},
	journal={Annals of Nuclear Energy},
	volume={88},
	pages={118-125},
	year={2016},
	note = {doi: 10.1016/j.anucene.2015.10.023},
	publisher={Elsevier}
}

@Article{lee2011nodal,
  title={Nodal integral expansion method for one-dimensional time-dependent linear convection--diffusion equation},
  author={Lee, Kyu Bok},
  journal={Nuclear engineering and design},
  volume={241},
  number={3},
  pages={767--774},
  year={2011},
  note = {doi: 10.1016/j.nucengdes.2010.11.016},
  publisher={Elsevier}
}

@Article{Niteen2022nodal,
	title={A nodal integral scheme for acoustic wavefield simulation},
	author={Kumar, Niteen and Shekar, Bharath and Singh, Suneet},
	journal={Pure and applied geophysics},
	volume={179},
	number={10},
	pages={3677-3691},
	year={2022},
	note = {doi: 10.1007/s00024-022-03160-3},
	publisher={Springer}
}

@Article{esser1993upwind,
	author  = {Esser, Peter D. and Witt, Robert J.},
	journal = {Nuclear Science and Engineering},
	title   = {An upwind nodal integral method for incompressible fluid flow},
	year    = {1993},
	note    = {doi: 10.13182/NSE93-A24011},
	number  = {1},
	pages   = {20--35},
	volume  = {114},
}

@Article{Neeraj2024coupled,
	author  = {Singh, Amritpal and Kumar, Neeraj},
	journal = {Journal of Computational Physics},
	title   = {Coupled nodal integral-immersed boundary method ({NI-IBM}) for simulating convection-diffusion physics},
	year    = {2024},
	note    = {doi:10.1016/j.jcp.2024.113394},
	pages   = {113394},
	volume  = {519},
}

@Article{wang2005modified,
	author  = {Wang, Fei and Rizwan-uddin},
	journal = {Nuclear Science and Engineering},
	title   = {Modified nodal integral method for the three-dimensional, time-dependent, incompressible Navier-Stokes equations},
	year    = {2005},
	note    = {doi:10.13182/nse149-107},
	number  = {1},
	pages   = {107--114},
	volume  = {149},
}
